%% file: pde-5-CR.tex
\documentclass[12pt,twoside,reqno,openany]{amsart}

\pdfoutput=1

\usepackage{amsbsy,amscd,amsfonts,amsmath,amssymb,amsthm,color,
fancybox,fancyhdr,footmisc,graphics,graphicx,ifthen,mathrsfs,
multicol,pdfpages,rotating,times,wasysym}
\usepackage[dvipsnames,svgnames,x11names]{xcolor}

\usepackage{pifont}

\usepackage[all]{xy}
\usepackage[utf8]{inputenc}
\usepackage[T1]{fontenc}
\usepackage{mathtools}
\newtagform{EngelLie}[\scriptstyle]{$}{$}
\makeatletter\newcommand{\leqnomode}{\tagsleft@true}
\newcommand{\reqnomode}{\tagsleft@false}\makeatother

\input macros.tex


\begin{document}

\setcounter{section}{0}

$\:$

\bigskip\bigskip\bigskip\bigskip\bigskip

\begin{center}

{\large\bf 
Equivalences\footnotemark[1] of PDE Systems}
\label{equivalences-pde-systems}

\medskip

{\large\bf Associated to Degenerate 
Para-CR Structures\footnotemark[2]:}

\medskip

{\large\bf
Foundational Aspects}

\end{center}

\footnotetext[1]{\,
This work was supported
in part by the Polish National Science Centre (NCN) 
via the grant number 2018/29/B/ST1/02583, and also by 
the GRIEG research project
{\sl Symmetry, Curvature Reduction, and
EquivAlence Methods} (with the registration
number 2019/34/H/ST1/00636), which is
funded by the Norwegian Financial Mechanism 2014-2021,
and is directed by Boris Kruglikov 
(PI, The Arctic University of Norway) 
and by Pawe{\l} Nurowski (PI, 
Centrum Fizyki Teoretycznej Polskiej Akademii Nauk).
}

\bigskip\bigskip

\begin{center}

Joël~{\sc Merker}\footnotemark[3]

\end{center}

\bigskip

\footnotetext[3]{\,\,Institut de Math\'ematique d'Orsay,
CNRS, Universit\'e Paris-Saclay, 91405 Orsay Cedex, France,
{\bf joel.merker@universite-paris-saclay.fr}}

\begin{center}
\begin{minipage}[t]{12.5cm}
\parindent 0.53cm
\scriptsize
\noindent
{\sc Abstract}.
Let $\K = \R$ or $\C$. We study basic invariants of
submanifolds of solutions $\mathcal{M} = \{
y = Q(x,a,b)\} = \{b = P(a,x,y)\}$ in coordinates
$x \in \K^{n\geqslant 1}$, $y \in \K$, $a \in \K^{m\geqslant 1}$,
$b \in \K$ under split-diffeomorphisms
$(x,y,a,b)
\,\longmapsto\,
\big(
f(x,y),\,g(x,y),\,\varphi(a,b),\,\psi(a,b)
\big)$. 
Two Levi forms exist, and have the same rank $r \leqslant \min
(n,m)$. If $\mathcal{M}$ is $k$-nondegenerate
with respect to parameters and $l$-nondegenerate
with respect to variables, $\Aut(\mathcal{M})$ 
is a local Lie group of dimension:
\[
\dim\,
\Aut
(\mathcal{M})
\,\,\leqslant\,\,
(n+1)\,
\tbinom{n+1+2k+2l}{n+1}
+
(m+1)\,
\tbinom{m+1+2k+2l}{m+1}.
\]

Mainly, our goal is to set up foundational material addressed to CR
geometers. We focus on $n = m = 2$, assuming $r = 1$. In coordinates
$(x,y,z, a,b,c)$, a local equation is:
\[
z
\,=\,
c
+
xa
+
\beta\,xxb
+
\underline{\beta}\,yaa
+
c\,{\rm O}_{x,y,a,b}(2)
+
{\rm O}_{x,y,a,b,c}(4),
\]
with $\beta$ and $\underline{\beta}$ 
representing the two $2$-nondegeneracy invariants at $0$.
The associated para-CR \pde~system:
\[
z_y
\,=\,
F\big(x,y,z,z_x,z_{xx}\big)
\ \ \ \ \ \ \ \ \ \ \ \ \
\&
\ \ \ \ \ \ \ \ \ \ \ \ \
z_{xxx}
\,=\,
H\big(x,y,z,z_x,z_{xx}\big),
\]
satisfies $F_{z_{xx}} \equiv 0$ from Levi degeneracy. 

We show in details that the
hypothesis of $2$-nondegeneracy with respect to variables
is equivalent to $F_{z_x z_x} \neq 0$. 
This gives a CR-geometric meaning to
the first two para-CR relative differential
invariants encountered independently in another
paper, joint with Pawe{\l} Nurowski.
\end{minipage}
\end{center}

\Section{\bf Introduction}
\label{introduction-pde-5-CR}
\HEAD{{\ref{introduction-pde-5-CR}}.~{\sf Introduction}
}{
Joël {\sc Merker} (Orsay)}

The first goal of this article is to develop some 
foundational aspects of para-CR structures, 
neither touched
in the first two memoirs
on the subject~{\cite{Merker-2008, Hill-Nurowski-2010}}, 
nor published somewhere, since then.
The second goal is to provide a detailed entrance
to an advanced paper~{\cite{Merker-Nurowski-2020-a}},
joint with Pawe{\l} Nurowski,
devoted to the classification
of Levi degenerate $5$-dimensional para-CR structures.
Thus, the focus is mainly on accessible aspects of para-CR structures,
at the threshold of more advanced equivalence results
about equivalences of Partial Diffenrential Equations,
which are highly demanding, computationally speaking. 

This introductory Section~{\ref{introduction-pde-5-CR}} presents, in a
self-contained and condensed way, some of the concepts and results of
the paper, emphasizing three new theorems taken from its core.

\smallskip

We work over $\K := \R$ or $\K := \C$. We study codimension $1$ {\sl
para-CR structures}~{\cite{Hill-Nurowski-2010}}, called {\sl
submanifolds of solutions} in~{\cite{Merker-2008}}. Let therefore
$x \in \K^n$, with $n \geqslant 1$, let $y \in \K$, let $a \in \K^m$,
with $m \geqslant 1$, and let $b \in \K$, be coordinates. In the CR
setting $x := z \in \C^n$, $y := w \in \C$, $m = n$ (due to complex
conjugation), $a := \overline{z}$, $b := \overline{w}$. In fact, CR
objects will be only seldom mentioned in this paper, since they are
less general than the para-CR ones.

A (local) {\sl submanifold of solutions} $\mathcal{M} \subset \K^{n+1}
\times \K^{m+1}$ is an analytic submanifold 
(hypersurface) passing through the
origin:
\[
\mathcal{M}
\,=\,
\big\{
(x,y,a,b)
\colon\,
R(x,y,a,b)
=
0
\big\}
\eqno
{\scriptstyle{(0\,\in\,\mathcal{M})}},
\] 
with $R_y \neq 0 \neq R_b$. By solving for $y$ or for $b$, we may
represent $\mathcal{M} = \{ y = Q(x,a,b)\}$ with $Q_b \neq 0$ or
$\mathcal{M} = \{ b = P(a,x,y)\}$ with $P_y \neq 0$. Such two
graphing functions satisfy identically:
\leqnomode\usetagform{default}
\begin{align}
\label{y-Q-P-b-identical}
y
\,\equiv\,
Q\big(x,a,P(a,x,y)\big)
\ \ \ \ \ \ \
\text{and}
\ \ \ \ \ \ \
b
\,\equiv\,
P\big(a,x,Q(x,a,b)\big).
\end{align}

The goal is to understand in a precise
manner para-CR invariants under local maps:
\leqnomode\usetagform{default}
\begin{align}
\label{para-CR-equivalence-intro}
\ \ \ \ \ \ \
(F,\Phi)
\colon
(x,y,a,b)
\,\longmapsto\,
\big(
f(x,y),\,g(x,y),\,\varphi(a,b),\,\psi(a,b)
\big)
\,=:\,
\big(x',y',a',b'\big),
\end{align}
which send $\mathcal{M}$ diffeomorphicallo to another
submanifold of solutions $\mathcal{M}'$
represented by similar equivalent equations:
\[
y'
\,=\,
Q'(x',a',b')
\ \ \ \ \ \ \
\text{and}
\ \ \ \ \ \ \
b'
\,=\,
P'(a',x',y').
\] 
In the CR setting, $m = n$, and 
equivalences are (local) biholomorphisms 
$\C^{n+1} \longrightarrow {\C'}^{n+1}$.
Equivalently, the following power series identities hold:
\reqnomode\usetagform{EngelLie}
\begin{align}
g\big(x,Q(x,a,b)\big)
&
\,\equiv\,
Q'
\big(
f(x,Q(x,a,b)),\varphi(a,b),\psi(a,b)
\big)
\tag{(\text{\rm in}\,\C\{x,a,b\}),}
\\
\psi\big(a,P(a,x,y)\big)
&
\,\equiv\,
P'
\big(
\varphi(a,P(a,x,y)),f(x,y),g(x,y)
\big)
\tag{(\text{\rm in}\,\C\{a,x,y\}).}
\end{align}

The link with completely integrable \pde~systems
whose general solution depends on a finite number of constants
is explained
in~{\cite{Merker-2008, Hill-Nurowski-2010}}.
Accordingly, $(x,y)$ are called {\sl variables}, "$\variables$"
(of the concerned \pde~system), while $(a,b)$ are
called {\sl parameters}, "$\parameters$"
(namely "constants of integration").
Thus $y := Q(x,a,b)$ can be thought of as being the general
solution to some \pde~system, and this explains
the terming "{\sl submanifold of solutions}".

Fixing either $(a,b)$ or $(x,y)$, we introduce two kinds of 
leaves which foliate $\mathcal{M}$:
\[
\mathcal{Q}_{a,b}
\,:=\,
\big\{
(x,y)
\colon\,
y
=
Q(x,a,b)
\big\}
\ \ \ \ \ \ \
\text{and}
\ \ \ \ \ \ \
\mathcal{P}_{x,y}
\,:=\,
\big\{
(a,b)
\colon\,
b
=
P(a,x,y)
\big\}.
\]
They are preserved through any para-CR equivalence:
\[
(F,\Phi)
\big(\mathcal{Q}_{a,b}\big)
\,\subset\,
\mathcal{Q}_{\varphi(a,b),\psi(a,b)}'
\ \ \ \ \ \ \
\text{and}
\ \ \ \ \ \ \
(F,\Phi)
\big(\mathcal{P}_{x,y}\big)
\,\subset\,
\mathcal{P}_{f(x,y),g(x,y)}'.
\]

In Section~{\ref{nonholonomy-form-and-dual}}, 
we introduce {\em two}
kinds of {\sl Levi forms}, whose two matrices
in some natural local frames are
the $n \times m$ matrix and the $m \times n$ matrix:
\[
\footnotesize
\aligned
\Levi_\parameterssmall(Q)
&
\,:=\,
\left(\!
\begin{array}{ccc}
\frac{-Q_bQ_{x_1a_1}+Q_{a_1}Q_{x_1b}}{Q_b\,Q_b}
& \cdots &
\frac{-Q_bQ_{x_1a_m}+Q_{a_m}Q_{x_1b}}{Q_b\,Q_b}
\\
\vdots & \ddots & \vdots
\\
\frac{-Q_bQ_{x_na_1}+Q_{a_1}Q_{x_nb}}{Q_b\,Q_b}
& \cdots &
\frac{-Q_bQ_{x_na_m}+Q_{a_m}Q_{x_nb}}{Q_b\,Q_b}
\end{array}
\!\right),
\\
\Levi_\variablessmall(P)
&
\,:=\,
\left(\!
\begin{array}{ccc}
\frac{-P_yP_{a_1x_1}+P_{x_1}P_{a_1y}}{P_y\,P_y}
& \cdots &
\frac{-P_yP_{a_1x_n}+P_{x_n}P_{a_1y}}{P_y\,P_y}
\\
\vdots & \ddots & \vdots
\\
\frac{-P_yP_{a_mx_1}+P_{x_1}P_{a_my}}{P_y\,P_y}
& \cdots &
\frac{-P_yP_{a_mx_n}+P_{x_n}P_{a_my}}{P_y\,P_y}
\end{array}
\!\right).
\endaligned
\]

By differentiating~({\ref{y-Q-P-b-identical}}),
we see that $P_y = \frac{1}{Q_b}$, and 
we show in Lemma~{\ref{Lemma-transfer-from-Q-to-P}} that
for all indices $1 \leqslant i \leqslant n$ and $1 \leqslant j
\leqslant m$, it holds identically on $\mathcal{M}$:
\[
\frac{-Q_bQ_{x_ia_j}+Q_{a_j}Q_{x_ib}}{
Q_b\,Q_b}
\,\,\equiv\,\,
-\,P_y\,
\bigg(
\frac{-P_yP_{x_ia_j}+P_{x_i}P_{a_jy}}{
P_y\,P_y}
\bigg).
\] 
From this, it follows that the two Levi forms 
have the same rank, since:
\[
\Levi_\parameterssmall(Q)
\,=\,
-\,P_y\,
{}^\TT
\Levi_\variablessmall(P)
\ \ \ \ \ \ \ \ \ \ \ \ \
\Longleftrightarrow
\ \ \ \ \ \ \ \ \ \ \ \ \
-\,Q_b\,
{}^\TT
\Levi_\parameterssmall(Q)
\,=\,
\Levi_\variablessmall(P).
\]

Denoting this common rank by $r \leqslant \min\, (n,m)$, 
we deduce that $\mathcal{M}$ 
can always be represented,
in suitably normalized coordinates $(x,y,a,b)$, as:
\[
\aligned
y
&
\,=\,
b
+
x_1a_1
+\cdots+
x_ra_r
+
{\rm O}_{x,a}(3)
+
b\,{\rm O}_{x,a,b}(2),
\\
\Longleftrightarrow
\ \ \ \ \ 
b
&
\,=\,
y
-
x_1a_1
-\cdots-
x_ra_r
+
{\rm O}_{a,x}(3)
+
y\,{\rm O}_{a,x,y}(2).
\endaligned
\]

In Section~{\ref{Levi-form-order-1-jets-Segre}}, we show that 
the two 
Levi forms are related to jets of order $1$
of $\mathcal{Q}_{a,b}$ and of $\mathcal{P}_{x,y}$. 
More generally, we study jet-maps of any orders
$k \geqslant 1$ and $l \geqslant 1$: 
\[
J_x^k\mathcal{Q}_{a,b}
\,:=\,
\Big(
x,\,\,
\big(
\partial_x^\beta Q(x,a,b)
\big)_{\vert\beta\vert\leqslant k}
\Big)
\ \ \ \ \ \ \
\text{and}
\ \ \ \ \ \ \
J_a^l\mathcal{P}_{x,y}
\,:=\,
\Big(
a,\,\,
\big(
\partial_a^\gamma
P(a,x,y)
\big)_{\vert\gamma\vert\leqslant l}
\Big).
\]
 
Assume a para-CR equivalence~({\ref{para-CR-equivalence-intro}})
is given.
In the less general CR context, the following
theorem already appeared
in~{\cite{Merker-2002, Merker-Porten-2006,
Merker-Pocchiola-Sabzevari-2013-5-CR-II}},
and in Section~{\ref{higher-order-jets-invariant-leaves}}, 
Theorem~{\ref{Thm-transfer-jets}} will provide
more details.

\begin{Theorem}
For any two integers $k \geqslant 1$ and $l \geqslant 1$, there are
two maps $R_{f,g}^k$ and $S_{\varphi,\psi}^l$ of the form:
\[
R_{f,g}^k\,
\Big(
x,y,\,
\big(
y_{x^\beta}
\big)_{1\leqslant\vert\beta\vert\leqslant k}
\Big)
\ \ \ \ \ \ \ \ \
\text{and}
\ \ \ \ \ \ \ \ \
S_{\varphi,\psi}^l\,
\Big(
a,b,\,
\big(
b_{a^\gamma}
\big)_{1\leqslant\vert\gamma\vert\leqslant l}
\Big),
\]
which make commutative the two diagrams:
\[
\aligned
\xymatrix{
\mathcal{M}
\ar[rr]^{(f,g,\varphi,\psi)}
\ar[d]_{J_x^k\mathcal{Q}_{a,b}}
&
&
\mathcal{M}'
\ar[d]^{J_{x'}^k\mathcal{Q}_{a',b'}'}
\\
\K^{n+\frac{(n+k)!}{n!\,k!}}
\ar[rr]_{R_{f,g}^k}
&
&
{\K'}^{n+\frac{(n+k)!}{n!\,k!}},
}
\ \ \ \ \ \ \ \ \ \ \ \ \
\ \ \ \ \ \ \ \ \ \ \ \ \
\xymatrix{
\mathcal{M}
\ar[rr]^{(\varphi,\psi,f,g)}
\ar[d]_{J_a^l\mathcal{P}_{x,y}}
&
&
\mathcal{M}'
\ar[d]^{J_{a'}^l\mathcal{P}_{x',y'}'}
\\
\K^{m+\frac{(m+l)!}{m!\,l!}}
\ar[rr]_{S_{\varphi,\psi}^l}
&
&
{\K'}^{m+\frac{(m+l)!}{m!\,l!}}.
}
\endaligned
\]
\end{Theorem}

A submanifold of solutions $\mathcal{M}$ is called:

\smallskip\noindent$\bullet$\, 
{\sl $k$-nondegenerate with respect to parameters} at a point $(x_0,
a_0, b_0)$ if there exists an integer $k \geqslant 1$ (smallest
possible) such that the invariant jet map $J_x^k \mathcal{Q}_{a,b}$:
\[
(x,a,b)
\,\,\longmapsto\,\,
\Big(
x,\,\,
\big(
\partial_x^\beta Q(x,a,b)
\big)_{\vert\beta\vert\leqslant k}
\Big)
\]
is of maximal possible rank $n + 1 + m$ at $(x_0, a_0, b_0)$;

\smallskip\noindent$\bullet$\, 
{\sl 
$l$-nondegenerate with respect to variables} at a point $(a_0, x_0,
y_0)$ if there exists an integer $l \geqslant 1$ (smallest possible)
such that the invariant jet map $J_a^l \mathcal{P}_{x,y}$:
\[
(a,x,y)
\,\,\longmapsto\,\,
\Big(
a,\,\,
\big(
\partial_a^\gamma P(a,x,y)
\big)_{\vert\gamma\vert\leqslant l}
\Big)
\]
is of maximal possible rank $m + 1 + n$ at $(a_0, x_0, y_0)$.

\smallskip

A general result not contained in~{\cite{Merker-2008}}
can be formally written by taking inspiration 
from~{\cite{Gaussier-Merker-2004}}. Its dimension bounds 
are general, but not sharp. 

\begin{Theorem}
If $\mathcal{M}$ is $k$-nondegenerate with respect to parameters
and $l$-nondegenerate with respect to variables,
then:
\[
\Aut(\mathcal{M})
\,:=\,
\big\{
(F,\Phi)
\in
\Diff_\variablessmall\times\Diff_\parameterssmall
\colon\,
(F,\Phi)(\mathcal{M})
\subset
\mathcal{M}
\big\},
\]
is a local Lie group of dimension:
\[
\dim\,
\Aut
(\mathcal{M})
\,\,\leqslant\,\,
(n+1)\,
\binom{n+1+2k+2l}{n+1}
+
(m+1)\,
\binom{m+1+2k+2l}{m+1}.
\]
\end{Theorem}

As is usual with any general result, there is a price to pay: its
hidden coarseness. For instance, with $n = m = 2$, and 
with $k = l = 2$, the bound is:
\[
\dim\,
\Aut
(\mathcal{M})
\,\,\leqslant\,\,
(2+1)\,
\tbinom{2+1+2\cdot 2+2\cdot2}{2+1}
+
(2+1)\,
\tbinom{2+1+2\cdot 2+2\cdot2}{2+1}
\,\,=\,\,
990,
\]
but the right bound, attained, 
is~{\cite{Merker-Nurowski-2020-a}}:
\[
\dim\,
\Aut
(\mathcal{M})
\,\,\leqslant\,\,
10.
\]

The present paper, preliminary to the 
articles~{\cite{Merker-Nurowski-2020-a, Merker-Nurowski-2020-b}},
examines mainly the dimensions $n = m = 2$.
More suitable notations for the two graphed equations of
$\mathcal{M}$ then are:
\[
z
\,=\,
Q(x,y,a,b,c)
\ \ \ \ \ \ \ \ \ \ \ \ \
\Longleftrightarrow
\ \ \ \ \ \ \ \ \ \ \ \ \
c
\,=\,
P(a,b,x,y,z).
\]
Also, we assume throughout that the Levi rank $r$ is constant
$r = 1 < 2$, so that both Levi forms are {\em degenerate}.
Thus, we assume that: 
\[
Q_b
\,\neq\, 
0
\,\neq\, 
P_y,
\ \ \ \ \ \ \ \ \ \ \ \ \ \
-\,Q_c\,Q_{xa}+Q_a\,Q_{xc}
\,\neq\,
0
\,\neq\,
-\,P_z\,P_{ax}+Px\,P_{az},
\]
and we assume the two (equivalent) identical
vanishings:
\[
\footnotesize
\aligned
0
&
\,\equiv\,
\det\,
\Levi_\parameterssmall(Q)
\,=\,
\left\vert\!
\begin{array}{ccc}
\frac{-Q_cQ_{xa}+Q_aQ_{xc}}{Q_c\,Q_c}
&
\frac{-Q_cQ_{xb}+Q_bQ_{xc}}{Q_c\,Q_c}
\\
\frac{-Q_cQ_{ya}+Q_aQ_{yc}}{Q_c\,Q_c}
&
\frac{-Q_cQ_{yb}+Q_bQ_{yc}}{Q_c\,Q_c}
\end{array}
\!\right\vert
\\
\ \ \ \ \
\Longleftrightarrow
\ \ \ \ \ 
0
&
\,\equiv\,
\det\,
\Levi_\variablessmall(P)
\,=\,
\left\vert\!
\begin{array}{ccc}
\frac{-P_zP_{ax}+P_xP_{az}}{P_z\,P_z}
&
\frac{-P_zP_{ay}+P_xP_{ay}}{P_z\,P_z}
\\
\frac{-P_zP_{bx}+P_xP_{bz}}{P_z\,P_z}
&
\frac{-P_zP_{by}+P_yP_{bz}}{P_z\,P_z}
\end{array}
\!\right\vert.
\endaligned
\]

In Section~{\ref{local-graphs-twin-2-nondegenerate-M}}, we show
that there exist normalized coordinates $(x,y,z,a,b,c)$ in which
these two equations read:
\[
\aligned
z
&
\,=\,
c
+
xa
+
\beta\,xxb
+
\underline{\beta}\,yaa
+
c\,{\rm O}_{x,y,a,b}(2)
+
{\rm O}_{x,y,a,b,c}(4),
\\
c
&
\,=\,
z
-
ax
-
\underline{\beta}\,aay
-
\beta\,bxx
+
z\,{\rm O}_{a,b,x,y}(2)
+
{\rm O}_{a,b,x,y,z}(4).
\endaligned
\]
Also, we verify that, at the origin $(x_0, y_0, z_0, a_0, b_0, c_0) = 
(0, 0, 0, 0, 0, 0)$, one has: 
\[
\aligned
\mathcal{M}\,\,
\text{is $2$-nondegenerate with respect to parameters}\,\,
\ \ \ \ \
&
\Longleftrightarrow
\ \ \ \ \
\beta
\,\neq\,
0
\ \ \ \ \
\Longleftrightarrow
\ \ \ \ \
\Delta(Q)
\,\neq\,
0,
\\
\mathcal{M}\,\,
\text{is $2$-nondegenerate with respect to variables}\,\,
\ \ \ \ \
&
\Longleftrightarrow
\ \ \ \ \
\underline{\beta}
\,\neq\,
0
\ \ \ \ \
\Longleftrightarrow
\ \ \ \ \
\Box(P)
\,\neq\,
0,
\endaligned
\]
where:
\[
\Delta(Q)
\,:=\,
\left\vert\!
\begin{array}{ccc}
Q_a & Q_b & Q_c
\\
Q_{xa} & Q_{xb} & Q_{xc}
\\
Q_{xxa} & Q_{xxb} & Q_{xxc}
\end{array}
\!\right\vert,
\ \ \ \ \ \ \ \ \ \ \ \ \
\ \ \ \ \ \ \ \ \ \ \ \ \
\Box(P)
\,:=\,
\left\vert\!
\begin{array}{ccc}
P_x & P_y & P_z
\\
P_{ax} & P_{ay} & P_{az}
\\
P_{aax} & P_{aay} & P_{aaz}
\end{array}
\!\right\vert.
\]
In the CR setting, $\underline{\beta} = \overline{\beta}$ ,
hence both $2$-nondegeneracy conditions are equivalent.
From now on, we shall assume that $\mathcal{M}$
is $2$-nondegenerate with respect to
parameters.
 
\smallskip

Next, in terms of the two equivalent choices of fundamental vector
fields:
\[
\aligned
\mathcal{K}_x
&
\,:=\,
\partial_x,
&
\ \ \ \ \
\mathcal{K}_y
&
\,:=\,
\partial_y,
&
\ \ \ \ \
\mathcal{L}_a
&
\,:=\,
\partial_a
-
\tfrac{Q_a}{Q_c}\,
\partial_c,
&
\ \ \ \ \
\mathcal{L}_b
&
\,:=\,
\partial_b
-
\tfrac{Q_b}{Q_c}\,
\partial_c,
\\
\mathcal{K}_x
&
\,:=\,
\partial_x
-
\tfrac{P_x}{P_z}\,
\partial_z,
&
\ \ \ \ \
\mathcal{K}_y
&
\,:=\,
\partial_y
-
\tfrac{P_y}{P_z}\,
\partial_z,
&
\ \ \ \ \
\mathcal{L}_a
&
\,:=\,
\partial_a,
&
\ \ \ \ \
\mathcal{L}_b
&
\,:=\,
\partial_b.
\endaligned
\]
we may introduce the two kernels of the two Levi forms:
\[
\aligned
\Ker\,\Levi_\variablessmall
&
\,:=\,
\Big\{
\mathcal{K}
\in
\Gamma\big(T^\variablessmall\mathcal{M}\big)
\colon\,
0
=
\Levi_\variablessmall
\big(
\mathcal{K},
\mathcal{L}
\big),
\ \ 
\forall\,
\mathcal{L}
\in
\Gamma
\big(
T^\parameterssmall\mathcal{M}
\big)
\Big\}, 
\\
\Ker\,\Levi_\parameterssmall
&
\,:=\,
\Big\{
\mathcal{L}
\in
\Gamma\big(T^\parameterssmall\mathcal{M}\big)
\colon\,
0
=
\Levi_\parameterssmall
\big(
\mathcal{L},
\mathcal{K}
\big),
\ \ 
\forall\,
\mathcal{K}
\in
\Gamma
\big(
T^\variablessmall\mathcal{M}
\big)
\Big\}.
\endaligned
\]
Abbreviating:
\[
\kaux
\,:=\,
-\,
\frac{-P_zP_{ay}+P_yP_{az}}{-P_zP_{ax}+P_xP_{az}}\,
\ \ \ \ \ \ \ \ \ \ \ \ \
\text{and}
\ \ \ \ \ \ \ \ \ \ \ \ \
\laux
\,:=\,
-\,
\frac{-Q_cQ_{xb}+Q_bQ_{xc}}{-Q_cQ_{xa}+Q_aQ_{xc}},
\]
we see that these two kernels are {\em line bundles}
generated by:
\[
\mathcal{K}_\kersmall^\variablessmall
\,=\,
\kaux\,
\mathcal{K}_x
+
\mathcal{K}_y
\ \ \ \ \ \ \ \ \ \ \ \ \
\text{and}
\ \ \ \ \ \ \ \ \ \ \ \ \
\mathcal{L}_\kersmall^\parameterssmall
\,=\,
\laux\,
\mathcal{L}_a
+
\mathcal{L}_b.
\] 
We will come back to these two Levi kernel bundles in a while.

\smallskip

For the moment, 
by eliminating the three parameters $(a, b, c)$ from the three
equations $y = Q$, $y_x = Q_x$, $y_{xx} = Q_{xx}$, 
and by substituting the solved $(a,b,c)$
in $z_y = Q_y$, $z_{xxx} = Q_{xxx}$, 
we reconstitute the
associated \pde~system which is the main object of study
in~{\cite{Merker-Nurowski-2020-a, Merker-Nurowski-2020-b}}:
\leqnomode\usetagform{default}
\begin{align}
\label{pde-system-intro}
z_y
\,=\,
F\big(x,y,z,z_x,z_{xx}\big)
\ \ \ \ \ \ \ \ \ \ \ \ \
\&
\ \ \ \ \ \ \ \ \ \ \ \ \
z_{xxx}
\,=\,
H\big(x,y,z,z_x,z_{xx}\big).
\end{align}
From $\big( Q_y\big)_{xxx} = \big( Q_{xxx} \big)_y$, we deduce a
compatibility constraint on $F$ and $H$:
\[
\Dop_x
\big(
\Dop_x
(\Dop_x(F))
\big)
\,=\,
\Dop_y
\big(H\big),
\] 
in terms of the two {\sl total differentiation operators:}
\[
\Dop_x
\,:=\,
\partial_x
+
z_x\,
\partial_z
+
z_{xx}\,
\partial_{z_x}
+
H\,
\partial_{z_{xx}},
\ \ \ \ \
\Dop_y
\,:=\,
\partial_y
+
F\,
\partial_z
+
\Dop_x(F)\,
\partial_{z_x}
+
\Dop_x\big(\Dop_x(F)\big)\,
\partial_{z_{xx}}.
\]
The converse holds true, and is an elementary consequence of the 
Frobenius theorem~{\cite{Merker-2008}}.

\smallskip

{\sl As a preliminary to the 
articles~{\cite{Merker-Nurowski-2020-a,
Merker-Nurowski-2020-b}},
we can now explain how the geometric (invariant) hypotheses
on submanifolds of solutions
do transfer to 
the para-CR \pde~systems}. At first, we realize in 
Proposition~{\ref{Prp-F-zxx-zero}} the following nice fact.

\begin{Proposition}
If $\mathcal{M}$ has (degenerate)
Levi form of constant rank $1$ and if it is $2$-nondegenerate
with respect to parameters, then $F$ is independent of $z_{xx}$:
\[
0
\,\equiv\,
F_{z_{xx}}.
\]
\end{Proposition}

Next, due to the {\em independency} between $\beta \neq 0$ and
$\underline{\beta} \neq 0$, there comes a

\begin{Question}
{\sl How can one view the hypothesis of $2$-nondegeneracy 
with respect to {\em variables} in
the \pde~system~({\ref{pde-system-intro}}){\bf ?}}
\end{Question}

By some (non-straightforward) differential-algebraic computations,
we show in Lemmas~{\ref{Lm-F-zx}} 
and in Proposition~{\ref{Proposition-k-equal-F-zx}} 
that:
\[
F_{z_x}
\,=\,
\frac{-Q_cQ_{ya}+Q_aQ_{yc}}{-Q_cQ_{xa}+Q_aQ_{xc}}
\,=\,
\frac{-P_zP_{ay}+P_yP_{az}}{-P_zP_{ax}+P_xP_{az}}
\,=\,
-\,
\kaux.
\]
Further, in Lemma~{\ref{Lm-da-F-zx}}, we show
the nowhere vanishing invariant expression:
\[
\frac{\partial}{\partial a}\,
\big(
F_{z_x}
\big)
\,=\,
P_z\,
\frac{
\left\vert\!
\begin{smallmatrix}
P_x & P_y & P_z
\\
P_{ax} & P_{ay} & P_{az}
\\
P_{aax} & P_{aay} & P_{aaz}
\end{smallmatrix}
\!\right\vert
}{
\big(
-\,P_z\,P_{ax}+P_x\,P_{az}
\big)^2
}
\,\,\neq\,\,
0,
\]
from which we deduce:

\begin{Theorem}
Let~({\ref{pde-system-intro}}) be a \pde~system
coming from a submanifold of solutions
$\mathcal{M}$ which is $2$-nondegenerate with respect
to parameters. Then 
$\mathcal{M}$ is also $2$-nondegenerate with respect
to variables if and only if:
\[
F_{z_xz_x}
\,\neq\,
0.
\]
\end{Theorem}

This theorem therefore explains, in standard CR terms
or para-CR language, why
in~{\cite{Merker-Nurowski-2020-a, Merker-Nurowski-2020-b}},
we encountered the natural differential
invariant $F_{z_x z_x}$, which we assumed $\neq 0$. 

In order to present a bit more of foundational
material preparatory to~{\cite{Merker-Nurowski-2020-a,
Merker-Nurowski-2020-b}},
let us describe how to launch Cartan's equivalence
method, both from the point of view of $y = Q(x,a,b)$ 
and from the point of view of $z_y = F$, $z_{xxx} = H$.

Coming back to our two Levi kernel bundles
$T_\kersmall^\variablessmall\mathcal{M}$ and
$T_\kersmall^\parameterssmall\mathcal{M}$, 
we observe that,
through any equivalence $(F, \Phi) \colon \mathcal{M}
\longrightarrow \mathcal{M}'$:
\[
(F,\Phi)_*
\big(
T_\kersmall^\variablessmall\mathcal{M}
\big)
\,=\,
T_\kersmall^\variablessmall\mathcal{M}'
\ \ \ \ \ \ \ \ 
\text{and}
\ \ \ \ \ \ \ \ 
(F,\Phi)_*
\big(
T_\kersmall^\parameterssmall\mathcal{M}
\big)
\,=\,
T_\kersmall^\parameterssmall\mathcal{M}'.
\]
It follows that there exist $11$ functions on 
$\mathcal{M}$ such that:
\[
\left(\!
\begin{array}{c}
\mathcal{T}_{c'}'
\\
\mathcal{L}_{a'}'
\\
\laux'\mathcal{L}_{a'}'+\mathcal{L}_{b'}'
\\
\mathcal{K}_{x'}'
\\
\kaux'\mathcal{K}_{x'}'+\mathcal{K}_{y'}'
\end{array}
\!\right)
\,=\,
\left(\!
\begin{array}{ccccc}
{\sf a} & {\sf b}_1 & {\sf b}_2 & {\sf c}_1 & {\sf c}_2
\\
0 & {\sf f}_1 & {\sf f}_2 & 0 & 0
\\
0 & 0 & {\sf f}_3 & 0 & 0
\\
0 & 0 & 0 & {\sf h}_1 & {\sf h}_2
\\
0 & 0 & 0 & 0 & {\sf h}_3
\end{array}
\!\right)\,
\left(\!
\begin{array}{c}
\mathcal{T}_{c}
\\
\mathcal{L}_{a}
\\
\laux\mathcal{L}_{a}+\mathcal{L}_{b}
\\
\mathcal{K}_{x}
\\
\kaux\mathcal{K}_{x}+\mathcal{K}_{y}
\end{array}
\!\right).
\]
Hence passing to dual $1$-forms, 
a plain transposition yields
(with new functions):
\[
\left(\!
\begin{array}{c}
Q_{a'}'da'+Q_{b'}'db'+Q_{c'}'dc'
\\
da'-\laux'db'
\\
db'
\\
dx'-\kaux'dy'
\\
dy'
\end{array}
\!\right)
\,=\,
\left(\!
\begin{array}{ccccc}
{\sf a} & 0 & 0 & 0 & 0
\\
{\sf b}_1 & {\sf f}_1 & 0 & 0 & 0
\\
{\sf b}_2 & {\sf f}_2 & {\sf f}_3 & 0 & 0
\\
{\sf c}_1 & 0 & 0 & {\sf h}_1 & 0
\\
{\sf c}_2 & 0 & 0 & {\sf h}_2 & {\sf h}_3
\end{array}
\!\right)\,
\left(\!
\begin{array}{c}
Q_ada+Q_bdb+dc
\\
da-\laux db
\\
db
\\
dx-\kaux dy
\\
dy
\end{array}
\!\right).
\]

From the point of view of the 
\pde~system~({\ref{pde-system-intro}}),
the natural coframe to start with when
running Cartan's equivalence 
method is:
\[
\aligned
\lambda
&
\,:=\,
dz
-
z_x\,dx
-
F\,dy,
\\
\mu_1
&
\,:=\,
dz_x
-
z_{xx}\,dx
-
\Dop(F)\,dy,
\\
\mu_2
&
\,:=\,
dz_{xx}
-
H\,dx
-
\Dop_x(\Dop_x(F))\,dy,
\\
\nu_1
&
\,:=\,
dx
-
F_{z_x}\,dy,
\\
\nu_2
&
\,:=\,
dy,
\endaligned
\] 
{\sl This is under
the assumption that $F_{z_xz_x} \neq 0$ only!}
In conclusion, this explains why the {\sl lifted coframe}
we start with in~{\cite{Merker-Nurowski-2020-a}} is:
\[
\left(\!
\begin{array}{c}
\pmb{\lambda}
\\
\pmb{\mu_1}
\\
\pmb{\mu_2}
\\
\pmb{\nu_1}
\\
\pmb{\nu_2}
\end{array}
\!\right)
\,:=\,
\left(\!
\begin{array}{ccccc}
{\sf a} & 0 & 0 & 0 & 0
\\
{\sf b}_1 & {\sf f}_1 & 0 & 0 & 0
\\
{\sf b}_2 & {\sf f}_2 & {\sf f}_3 & 0 & 0
\\
{\sf c}_1 & 0 & 0 & {\sf h}_1 & {\sf h}_4
\\
{\sf c}_2 & 0 & 0 & {\sf h}_2 & {\sf h}_3
\end{array}
\!\right)
\left(\!
\begin{array}{c}
\lambda
\\
\mu_1
\\
\mu_2
\\
\nu_1
\\
\nu_2
\end{array}
\!\right).
\]
Beyond these elementary foundational considerations,
readers interested in {\em advanced} Cartan-type results
are referred 
to~{\cite{Merker-Nurowski-2020-a,
Merker-Nurowski-2020-b, Merker-Nurowski-2021}}.

\medskip\noindent{\bf Acknowledgments.}
The author would like to thank two referees 
for deep understanding and for wise advices.

\Section{\bf Motivation From CR Geometry}
\label{invariance-Omega}
\HEAD{{\ref{invariance-Omega}}.~{\sf Motivation From CR Geometry}
}{
Joël {\sc Merker} (Orsay)}

Already in 1907, Poincar\'e observed that there are many more
real-analytic ($\mathcal{C}^\omega$) 
hypersurfaces $M^3 \subset \C^2$, represented in coordinates
$(z,w) = (x+\isqrt\, y,\, u+\isqrt\, v)$ as graphs
by means of converging real power series:
\[
u
\,=\,
\varphi\big(x,y,v\big)
\,=\,
\smallsum{i,j,k}\,\,
\varphi_{i,j,k}\,
x^iy^jv^k
\eqno
{\scriptstyle{(\varphi_{i,j,k}\,\in\,\R)}},
\] 
than there are local biholomorphic transformations:
\[
(z,w)
\,\longmapsto\,
\big(f(z,w),\,g(z,w)\big)
\,\,=\,\,
\Big(
\smallsum{i,j}\,
f_{i,j}\,
z^i\,w^j,\,\,
\smallsum{i,j}\,
g_{i,j}\,
z^i\,w^j
\Big)
\eqno
{\scriptstyle{(f_{i,j},\,g_{i,j}\,\in\,\C)}},
\]
since for instance if all objects are polynomials of a certain large 
degree $d\gg 1$:
\[
\dim_\R\,
\big\{
\varphi_{i,j,k}\in\R
\colon\,
i+j+k\leqslant d
\big\}
\,=\,
{\textstyle{\binom{d+3}{3}}}
\,\,\gg\,\,
2\cdot 2\cdot
{\textstyle{\binom{d+2}{2}}}
\,=\,
\dim_\R\,
\big\{
f_{i,j},g_{i,j}
\in\C
\colon\,
i+j\leqslant d
\big\}.
\]

Hence there {\em is} a problem of classification, which can be
formulated generally in any dimension $\NN \geqslant 1$, and in a
local setting: given two real submanifolds $M \subset \C^\NN$ and $M'
\subset {\C'}^\NN$, determine
whether there exist local biholomorphic transformations $h
\colon \C^\NN \longrightarrow {\C'}^\NN$ with $h(M) = M'$ (up to
restriction to open subsets),
and classify such real submanifolds of
$\C^\NN$ up to local biholomorphisms.

Here, to avoid dwelling into the analysis of embedability problems and
to benefit of Zariski-generic properties, only real-analytic
($\mathcal{C}^\omega$) objects will be dealt with. Even in this
context, the problem is too ample to be solved completely in any
dimension.

To start with, two basic facts must be reminded, {\em see e.g.}
{\cite[Section~1]{Merker-Pocchiola-Sabzevari-2013-5-CR-II}} 
for details. Let $\NN \geqslant 1$ and let:
\[
\big(\zaux_1,\dots,\zaux_\NN\big)
\,\in\,
\C^\NN.
\]
Let $J \colon T\C^\NN \longrightarrow T\C^\NN$ be the
standard complex structure which expresses the multiplication
by $\isqrt$ of real vector fields:
\[
J
\big(
{\textstyle{\frac{\partial}{\partial x_k}}}
\big)
\,:=\,
{\textstyle{\frac{\partial}{\partial y_k}}}
\ \ \ \ \ \ \ \ \ \ \ \ \
\text{and}
\ \ \ \ \ \ \ \ \ \ \ \ \
J
\big(
{\textstyle{\frac{\partial}{\partial y_k}}}
\big)
\,:=\,
-\,{\textstyle{\frac{\partial}{\partial x_k}}}
\eqno
{\scriptstyle{(1\,\leqslant\,k\,\leqslant\,n)}}.
\]
A $\mathcal{C}^\omega$ 
real submanifold $M \subset \C^\NN$ is called {\sl Cauchy-Riemann}
({\sl CR}) when the smallest $J$-invariant subspaces of the tangent
spaces $T_pM$:
\[
\dim_\R\,
\big(T_pM\cap J(T_pM)\big)
\,=:\,
2\,n,
\]
are of constant (necessarily even) dimension when $p \in M$ varies.
It is called {\sl CR-generic} when its tangent spaces $T_pM$ 
generate the whole ambient tangent spaces:
\[
T_pM
+
J(T_pM)
\,=\,
T_p\C^\NN
\eqno
{\scriptstyle{(\forall\,p\,\in\,M)}}.
\]
Elementary linear algebra shows that CR-generic implies CR. 
The {\sl complex tangent spaces} are then:
\[
T_p^cM
\,:=\,
T_pM\cap J(T_pM)
\eqno
{\scriptstyle{(p\,\in\,M)}}.
\]

The first fact is that any $\mathcal{C}^\omega$ real-analytic
submanifold $M \subset \C^\NN$ (even every real-analytic subset) is CR
at a Zariski-generic point, that is, outside some proper real-analytic
subset. Agreeing to study 
problems only at generic points like Lie and Cartan did,
one {\em disregards} such exceptional sets.

The second basic fact is that every $\mathcal{C}^\omega$ CR
submanifold $M \subset \C^\NN$ is locally straightenable, by means of
some appropriate biholomorphism, to a {\sl CR-generic} one times zero:
\[
M
\,=\,
M_1
\times
\{0\}^{\NN-\NN_1}
\,\,\subset\,\,
\C^{\NN_1}
\times
\C^{\NN-\NN_1},
\]
with $M_1 \subset \C^{\NN_1}$ being CR-generic. Hence the problem
reduces to equivalences of CR-generic submanifolds, and in
this context, a deep link exists with equivalences of
certain systems of partial differential equations,
as we will see. 

Before coming to {\sc pde}'s, consider a CR-generic $M \subset
\C^\NN$ of real codimension: 
\[
c
\,:=\,
\codim_\R\,
M,
\]
represented locally
in some neighborhood of some `central' 
point $p_0 \in M$ as the zero-set:
\[
M
\,=\,
\big\{
\zaux\in\C^\NN
\colon\,
\rho_1(\zaux,\overline{\zaux})
=\cdots=
\rho_c(\zaux,\overline{\zaux})
=
0
\big\}
\eqno
{\scriptstyle{(\rho(p_0)\,=\,0)}},
\]
of $c$ real-analytic functions, which, to
guarantee smoothness, have independent differentials:
\[
0
\,\neq\,
d\rho_1\wedge\cdots\wedge d\rho_c(p)
\eqno
{\scriptstyle{(\forall\,p\,\in\,M)}}.
\]
The reality $\overline{\rho} = \rho$
of this $\R^c$-valued function
$\rho = \big(\rho_1, \dots, \rho_c\big)$ means,
when expanding it in converging power series:
\[
\smallsum{\alpha,\beta}\,
\rho_{\alpha,\beta}\,
\zaux^\alpha\,\overline{\zaux}^\beta
\,=\,
\rho\big(\zaux,\overline{\zaux}\big)
\,=\,
\overline{\rho(\zaux,\overline{\zaux})}
\,=\,
\smallsum{\alpha,\beta}\,
\overline{\rho_{\alpha,\beta}}\,
\overline{\zaux}^\alpha\,
\zaux^\beta,
\]
that its coefficients satisfy by identification:
\[
\overline{\rho_{\beta,\alpha}}
\,=\,
\rho_{\alpha,\beta}
\eqno
{\scriptstyle{(\forall\,\alpha\,\in\,\N^\NN,\,\,\forall\,\beta\,\in\,
\N^\NN)}},
\]
so that if one defines, by conjugating {\em only} coefficients:
\[
\overline{\rho}
\big(\zaux,\overline{\zaux}\big)
\,:=\,
\smallsum{\alpha,\beta}\,
\overline{\rho_{\alpha,\beta}}\,
\zaux^\alpha\,\overline{\zaux}^\beta,
\]
one indeed has an identity which expresses $\rho = 
\overline{\rho}$ with arguments:
\[
\rho
\big(\zaux,\overline{\zaux}\big)
\,\equiv\,
\overline{\rho}
\big(\overline{\zaux},\zaux\big)
\eqno
{\scriptstyle{(\text{\rm in}\,\C\{\zaux,\,\overline{\zaux}\}^c)}}.
\]

As is known, the CR-genericity of $M$ is then equivalent to:
\[
0
\,\neq\,
\partial\rho_1\wedge\cdots\wedge\partial\rho_c(p)
\eqno
{\scriptstyle{(\forall\,p\,\in\,M)}},
\]
where, for any $\mathcal{C}^\omega$ function $\chi$, the
holomorphic $\partial$ and antiholomorphic
$\overline{\partial}$ 
differential operators are defined by:
\[
\partial\chi
\,:=\,
\smallsum{k}\,
{\textstyle{\frac{\partial\chi}{\partial\zaux_k}}}\,
d\zaux_k,
\ \ \ \ \ \ \ \ \ \ \ \ \
\overline{\partial}\chi
\,:=\,
\smallsum{k}\,
{\textstyle{\frac{\partial\chi}{\partial\overline{\zaux}_k}}}\,
d\overline{\zaux}_k,
\ \ \ \ \ \ \ \ \ \ \ \ \
d\chi
\,=\,
\partial\chi
+
\overline{\partial}\chi,
\]
and have sum equal to the standard real differential.
Then the integer $\NN - c =: 2\,n$ is even
{\cite{Merker-Pocchiola-Sabzevari-2013-5-CR-II}}, and its half:
\[
n
\,:=\,
\CRdim\,
M
\]
is called the {\sl CR dimension} of the CR-generic submanifold $M
\subset \C^{n+c}$. Always, to avoid purely real and
purely complex geometries, we will assume:
\[
c
\,\geqslant\,
1
\ \ \ \ \ \ \ \ \ \ \ \ \ \ \ \ \ \ \
\text{and}
\ \ \ \ \ \ \ \ \ \ \ \ \ \ \ \ \ \ \
n
\,\geqslant\,
1,
\]
hence:
\[
\NN
\,=\,
n+c
\,\geqslant\,
2.
\]

With any $c \times c$ invertible matrix $(\alpha_{j,k})$ of
$\mathcal{C}^\omega$ functions, the functions $\rho_j' := \sum_k\,
\alpha_{j,k}\,\rho_k$ are still defining $M = \{ \rho_1' = \cdots =
\rho_c' = 0\}$, which shows a (known)
lack of adequate correspondence
$M \longleftrightarrow \rho$. 

Actually, most ot the time, it is more appropriate to work with
graphed representations of submanifolds, since the link
with intrinsic geometric properties becomes one-to-one.

\begin{Lemma}
{\cite{Merker-Pocchiola-Sabzevari-2013-5-CR-II}}
Given a $\mathcal{C}^\omega$ CR-generic $M \subset \C^{n+c}$
with $\codim_\R\, M = c \geqslant 1$ and
$\CRdim\, M = n \geqslant 1$, at each point $p_0 \in M$, 
there exist centered holomorphic coordinates:
\[
\big(z,w\big)
\,=\,
\big(z_1,\dots,z_n,\,w_1,\dots,w_c\big)
\,\in\,
\C^n\times\C^c,
\]
and a complex-analytic $\C^c$-valued function $\Theta\big(z,
\overline{z}, \overline{w}\big)$ such that, near $p$:
\[
M
\,=\,
\big\{
(z,w)\in\C^\NN
\colon\,
w_1=\Theta_1\big(z,\overline{z},\overline{w}\big),
\,\dots\dots,\,
w_c=\Theta_c\big(z,\overline{z},\overline{w}\big)
\big\}.
\eqno\qed
\]
\end{Lemma}

This vector-valued graphing function 
$\Theta = (\Theta_1, \dots, \Theta_c)$ 
is constructed by first performing a complex-linear
transformation in order that:
\[
0
\,\neq\,
\det\,
\big(
{\textstyle{\frac{\partial\rho_{j_1}}{\partial w_{j_2}}}}
(p_0)
\big),
\]
and then by applying the analytic implicit function theorem to solve
for $w$ in the $c$ equations:
\[
0
\,=\,
\rho_1\big(z,w,\overline{z},\overline{w}\big)
\,=\,\cdots\,=\,
\rho_c\big(z,w,\overline{z},\overline{w}\big).
\]

With the
convention, seen above, 
that the conjugation of power series distributes
bars over coefficients and over variables as well:
\[
\aligned
\overline{\Theta(z,\overline{z},\overline{w})}
\,=\,
\overline{
\smallsum{\alpha,\beta,\gamma}\,
\Theta_{\alpha,\beta,\gamma}\,
z^\alpha\,
\overline{z}^\beta\,
\overline{w}^\gamma}
&
\,=\,
\smallsum{\alpha,\beta,\gamma}\,
\overline{\Theta_{\alpha,\beta,\gamma}}\
\overline{z}^\alpha\,
z^\beta\,
w^\gamma
\\
&
\,=:\,
\overline{\Theta}
\big(\overline{z},z,w\big),
\endaligned
\]
this graphed $\C^c$-valued function representing:
\[
M
\,=\,
\big\{
w
=
\Theta\big(z,\overline{z},\overline{w}\big)
\big\},
\]
then satisfies identically by construction:
\[
0
\,\equiv\,
\rho
\big(
z,
\Theta(z,\overline{z},\overline{w}),
\overline{z},
\overline{w}
\big).
\]

Conjugating this and remembering $\overline{\rho}
(\overline{\zaux},\zaux) \equiv 
\rho (\zaux, \overline{\zaux})$ yields:
\[
\aligned
0
&
\,\equiv\,
\rho
\big(
z,w,\overline{z},
\overline{\Theta}(\overline{z},z,w)
\big),
\endaligned
\]
which means that the equations $\rho = 0$ can also be solved
with respect to $\overline{w}$ as:
\[
\overline{w}
=
\overline{\Theta}
\big(
\overline{z},z,w
\big).
\]
A more accurate and complete analysis conducts to the basic

\begin{Lemma}
The complex $\C^c$-valued graphing function $\Theta (z, 
\overline{z}, \overline{w})$ for a CR-generic submanifold
$M \subset \C^{n+c}$ satisfies the two (equivalent by conjugation)
functional equations:
\[
\aligned
w
&
\,\equiv\,
\Theta
\big(
z,\overline{z},
\overline{\Theta}(\overline{z},z,w)
\big),
\\
\overline{w}
&
\,\equiv\,
\overline{\Theta}
\big(
\overline{z},z,
\Theta(z,\overline{z},\overline{w})
\big).
\endaligned
\]

Conversely, given a $\C^c$-valued power series $\Theta \in \C \{z, 
\overline{z}, \overline{w} \}^c$ vanishing at
the origin which satisfies these identities,
the two zero-sets:
\[
\big\{
w
=
\Theta(z,\overline{z},\overline{w})
\big\}
\,=\,
\big\{
\overline{w}
=
\overline{\Theta}
(\overline{z},z,w)
\big\},
\]
coincide and define a local 
real-analytic CR-generic submanifold $M \subset
\C^{n+c}$.\qed
\end{Lemma}

Next, assume that two CR-generic 
submanifolds passing through the origin
$0 \in M \subset \C^{n+c}$ and $0' \in M' 
\subset {\C'}^{n+c}$ having the same codimension $c$ and
the same CR dimension $n$ represented as:
\[
M
\,=\,
\big\{
w 
= 
\Theta(z,\overline{z},\overline{w})
\big\}
\ \ \ \ \ \ \ \ \ \ \ \ \
\text{and}
\ \ \ \ \ \ \ \ \ \ \ \ \
M'
\,=\,
\big\{
w'
= 
\Theta'(z',\overline{z}',\overline{w}')
\big\},
\]
are {\sl equivalent} through a local biholomorphism:
\[
h
\colon\ \ \
(z,w)
\,\longmapsto\,
\big(
f(z,w),g(z,w)
\big)
\,=:\,
(z',w'),
\]
that is to say, assume that $h(M) \subset M'$. This
inclusion can be expressed
as two (equivalent) relations holding for $(z,w) \in M$:
\[
\aligned
g(z,w)
&
\,=\,
\Theta'
\big(
f(z,w),\overline{f}(\overline{z},\overline{w}),
\overline{g}(\overline{z},\overline{w})
\big),
\\
\overline{g}(\overline{z},\overline{w})
&
\,=\,
\overline{\Theta}'
\big(
\overline{f}(\overline{z},\overline{w}),
f(z,w),g(z,w)
\big).
\endaligned
\]
As is visible in these relations, simultaneously with 
the biholomorphic equivalence $h$, there comes its antiholomorphic
conjugation:
\[
\overline{h}
\colon\ \ \
(\overline{z},\overline{w})
\,\longmapsto\,
\big(
\overline{f}(\overline{z},\overline{w}),\,
\overline{g}(\overline{z},\overline{w})
\big),
\]
and consequently, it is natural to consider simutaneously both maps:
\[
\big(z,w,\overline{z},\overline{w}\big)
\,=\,
\big(\zaux,\overline{\zaux}\big)
\,\,\longmapsto\,\,
\big(h(\zaux),\,
\overline{h}(\overline{\zaux})\big)
\,=\,
\big(
f(z,w),g(z,w),\,
\overline{f}(\overline{z},\overline{w}),
\overline{g}(\overline{z},\overline{w})
\big).
\]
These aspects of thought will now become more transparent when passing
to equivalences of {\sc pde} systems.

\Section{\bf Equivalences of Submanifolds of Solutions
for {\sc pde} Systems}
\label{equivalences-mathcal-M}
\HEAD{{\ref{equivalences-mathcal-M}}.~{\sf Equivalences of 
Submanifolds of Solutions for {\sc pde} Systems}
}{
Joël {\sc Merker} (Orsay)}

As explained in~{\cite{Merker-2008, Hill-Nurowski-2010}}, 
point equivalences
of completely integrable systems of partial differentials equations 
whose solutions depend on a finite number of parameters
can be understood as point equivalences of their
respective {\sl submanifolds of solutions}, which share remarkable
common properties with CR-generic submanifolds.
The theory works equally well 
over the fields $\K = \R$ or $\K = \C$.

The general set up is a local product space of {\sl variables}
$(x,y)$ times {\sl parameters} $(a,b)$, in any dimension:
\[
x
\,=\,
\big(x_1,\dots,x_n\big),
\ \ \ \ \ \ \ \ \ \ \ \ \
y
\,=\,
\big(y_1,\dots,y_c\big),
\ \ \ \ \ \ \ \ \ \ \ \ \
a
\,=\,
\big(a_1,\dots,a_m\big),
\ \ \ \ \ \ \ \ \ \ \ \ \
b
\,=\,
\big(b_1,\dots,b_c\big),
\]
with $n \geqslant 1$, with $m \geqslant 1$ and with $c \geqslant 1$
being arbitrary 
integers, the number of parameters $b$'s being the same as
that of variables $y$'s, because we will understand $y(x) = Q(x,a,b)$
as a function of the variables $x$ depending on the parameters
$(a,b)$, with $b = Q(0,a,b) = y(0)$ being its value
at the origin. Abbreviate:
\[
\NN
\,:=\,
n+c
\ \ \ \ \ \ \ \ \ \ \ \ \
\text{and}
\ \ \ \ \ \ \ \ \ \ \ \ \
\MM
\,:=\,
m+c,
\]

For CR-generic $M \subset \C^{n+c}$, one has the correspondence:
\[
x
\,\longleftrightarrow\,
z,
\ \ \ \ \ \ \ \ \ \ \ \ \
y
\,\longleftrightarrow\,
w,
\ \ \ \ \ \ \ \ \ \ \ \ \
a
\,\longleftrightarrow\,
\overline{z},
\ \ \ \ \ \ \ \ \ \ \ \ \
b
\,\longleftrightarrow\,
\overline{w},
\]
which imposes $m = n$, a restriction that will hence be removed from
now on. Equivalences of CR manifolds will neverheless remain
inspirational for these more general structures, 
termed {\sl para-CR} in~{\cite{Hill-Nurowski-2010}}.

By definition, 
point equivalences of {\sc pde} systems transform graphs of solutions
into graphs of solutions, and as these graphs are marked out 
by parameters, such equivalences naturally induce equivalences
in the parameter space. Consequently, the
natural infinite pseudogroup of equivalences
is the product:
\[
\Diff_\variablessmall
\times
\Diff_\parameterssmall
\]
of diffeomorphisms in the space of variables times diffeomorphisms
in the space of parameters, namely invertible transformations of the 
form:
\[
\big(x,y,a,b)
\,\,\longmapsto\,\,
\big(
f(x,y),g(x,y),\,\varphi(a,b),\psi(a,b)
\big).
\]
In the CR case, we have:
\[
\varphi
\,\longleftrightarrow\,
\overline{f},
\ \ \ \ \ \ \ \ \ \ \ \ \
\psi
\,\longleftrightarrow\,
\overline{g}.
\]
Denoting symbolically:
\[
\variables
\,:=\,
(x,y)
\,\in\,
\K^\NN
\ \ \ \ \ \ \ \ \ \ \ \ \
\text{and}
\ \ \ \ \ \ \ \ \ \ \ \ \
\parameters
\,:=\,
(a,b)
\,\in\,
\K^\MM,
\]
the infinite pseudogroup 
$\Diff_\variablessmall
\times
\Diff_\parameterssmall$ can be thought of as 
the group of {\sl split diffeomorphisms} in
the total space $\K^\NN \times \K^\MM \ni (x,y,a,b)$. 

\begin{Definition}
A (local) 
{\sl submanifold of solutions} $\mathcal{M} \subset \K^{n+c}
\times \K^{m+c}$ is a $c$-codimensional
$\K$-analytic submanifold passing through the origin:
\[
\mathcal{M}
\,=\,
\big\{
(x,y,a,b)
\colon\,
R(x,y,a,b)
=
0
\big\}
\eqno
{\scriptstyle{(0\,\in\,\mathcal{M})}},
\] 
with its defining $\K^c$-valued function $R \in 
\K\{x, y, a, b\}^c$ satisfying:
\[
\det\,
\big(
{\textstyle{\frac{\partial R_{j_1}}{\partial y_{j_2}}}}
(0)
\big)
\,\neq\,
0
\,\neq\,
\det\,
\big(
{\textstyle{\frac{\partial R_{j_1}}{\partial b_{j_2}}}}
(0)
\big)
\eqno
{\scriptstyle{(1\,\leqslant\,j_1,\,j_2\,\leqslant\,c)}}.
\]
\end{Definition}

The analytic implicit function theorem then yields two graphing
functions $Q = Q(x,a,b)$ and $P = P(a,x,y)$ so that 
\[
\aligned
R(x,y,a,b)
\,=\,
0
\ \ \ \ \ \ \ \ \ \ \ \ \
&
\Longleftrightarrow
\ \ \ \ \ \ \ \ \ \ \ \ \
y
\,=\,
Q(x,a,b)
\\
&
\Longleftrightarrow
\ \ \ \ \ \ \ \ \ \ \ \ \
b
\,=\,
P(a,x,y),
\endaligned
\]
and this can be expressed as two functional equations holding
identically:
\leqnomode\usetagform{default}
\begin{align}
\label{0-R-Q-identically}
0
&
\,\equiv\,
R\big(x,Q(x,a,b),a,b\big),
\\
\label{0-R-P-identically}
0
&
\,\equiv\,
R\big(x,y,a,P(a,x,y)\big),
\end{align}
in $\C\{x,a,b\}^c$ and in $\C\{a,x,y\}^c$.

The problem of equivalence considered here is the one of submanifolds
of solutions considered up to split diffeomorphisms belonging to the
infinite pseudogroup $\Diff_\variablessmall \times
\Diff_\parameterssmall$.

\Section{\bf Normal Coordinates for Submanifolds of Solutions}
\label{normal-coordinates-M}
\HEAD{{\ref{normal-coordinates-M}}.~{\sf Normal Coordinates for 
Submanifolds of Solutions}
}{
Joël {\sc Merker} (Orsay)}

A first normalization, inspired from CR geometry, provides (for free)
useful normalizing coordinates in which, later, some invariant
concepts of the theory will be more transparently seen.

\begin{Proposition}
\label{Proposition-normal-coordinates}
After a $\K$-analytic change of coordinates belonging to
$\Diff_\variablessmall \times \Diff_\parameterssmall$ 
having the shape:
\[
x
\,=\,
x',
\ \ \ \ \ \ \ \ \ \ \ \ \
y
\,=\,
{\tt Y}'(x',y'),
\ \ \ \ \ \ \ \ \ \ \ \ \
a
\,=\,
a',
\ \ \ \ \ \ \ \ \ \ \ \ \
b
\,=\,
{\tt B}'(a',b'),
\]
$\mathcal{M}$ is transformed
into an new $\mathcal{M}'$
whose new $\K^c$-valued implicit defining function:
\[
R'\big(x',y',a',b'\big)
\,:=\,
R
\big(
x',{\tt Y}'(x',y'),\,a',{\tt B}'(a',b')
\big),
\]
satisfies the normalization:
\[
R'\big(0,b',a',b')
\,\equiv\,
0
\ \ \ \ \ \ \ \ \ \ \ \ \
\text{and}
\ \ \ \ \ \ \ \ \ \ \ \ \
0
\,\equiv\,
R'
\big(x',y',0,y'\big).
\]
\end{Proposition}

Formally in the proof below and in some later ones,
thinking of $\K^c$-valued functions or
assuming for simplicity that $c = 1$ works equally well.

\proof
After a linear change of coordinates belonging
to $\Diff_\variablessmall \times \Diff_\parameterssmall$
which stabilizes $x' = x$ and $a' = a$, 
we can assume that:
\leqnomode\usetagform{default}
\begin{align}
\label{linear-normalization-R}
R
\,=\,
y
-
b
+
{\rm O}_2(x,y,a,b).
\end{align}
Set $x := 0$ and solve for $y$ in the equation: 
\[
0
\,=\,
R\big(0,y,a,b\big)
\ \ \ \ \ \ \ \ \ \ \ \ \
\Longleftrightarrow
\ \ \ \ \ \ \ \ \ \ \ \ \
y
\,=\,
{\tt U}(a,b)
\,=\,
b
+
{\rm O}_2(a,b),
\]
with a certain unique analytic map ${\tt U}$ hence satisfying:
\leqnomode\usetagform{default}
\begin{align}
\label{R-U-a-b}
0
\,\equiv\,
R\big(
0,{\tt U}(a,b),a,b\big).
\end{align}

Next, solve:
\[
{\tt U}(a,b)
\,=\,
y
\ \ \ \ \ \ \ \ \ \ \ \ \
\Longleftrightarrow
\ \ \ \ \ \ \ \ \ \ \ \ \
b
\,=\,
{\tt B}(a,y),
\]
and then take as defining equation for the transformed 
$\mathcal{M}'$:
\[
R'
\big(x,y,a,b'\big)
\,:=\,
R
\big(x,y,a,{\tt B}(a,b')\big).
\]

Using ${\tt U}\big( a, {\tt B}(a,b')\big)
\equiv b'$ when replacing $b := {\tt B}(a,b')$ 
in~({\ref{R-U-a-b}}) conducts to the first claimed
normalization:
\[
\aligned
0
&
\,\equiv\,
R
\big(
0,{\tt U}(a,{\tt B}(a,b')),a,{\tt B}(a,b')
\big)
\\
&
\,\equiv\,
R
\big(
0,b',a,{\tt B}(a,b')
\big)
\\
&
\,\equiv\,
R'
\big(0,b',a,b'\big),
\endaligned
\]
which we rewrite, dropping primes, simply as:
\leqnomode\usetagform{default}
\begin{align}
\label{first-normalization-R}
0
\,\equiv\,
R
\big(0,b,a,b\big).
\end{align}

To reach the second normalization, observe at first that 
setting $a := 0$ gives:
\leqnomode\usetagform{default}
\begin{align}
\label{R-0-b-0}
0
\,\equiv\,
R
\big(0,b,0,b\big).
\end{align}
Now, solve for $b$ in the equation:
\[
0
\,=\,
R
\big(x,y,0,b\big)
\ \ \ \ \ \ \ \ \ \ \ \ \
\Longleftrightarrow
\ \ \ \ \ \ \ \ \ \ \ \ \
b
\,=\,
{\tt V}(x,y)
\,=\,
y
+
{\rm O}_2(x,y),
\]
to get a certain unique analytic map ${\tt V}$ 
with ${\tt V}(0,y) \equiv y$ then satisfying:
\leqnomode\usetagform{default}
\begin{align}
\label{R-V-x-y}
0
\,\equiv\,
R\big(
x,y,0,{\tt V}(x,y)\big).
\end{align}
Next, solve:
\[
{\tt V}(x,y)
\,=\,
b
\ \ \ \ \ \ \ \ \ \ \ \ \
\Longleftrightarrow
\ \ \ \ \ \ \ \ \ \ \ \ \
y
\,=\,
{\tt Y}(x,b),
\]
and observe that ${\tt Y}(0,b) \equiv b$ by uniqueness
in the implicit function theorem, since
${\tt V}(0,y) \equiv y$. 

Change defining equation as:
\[
R'
\big(
x,y',a,b
\big)
\,:=\,
R
\big(
x,{\tt Y}(x,y'),a,b
\big),
\]
and, importantly, observe that the first normalization is preserved:
\[
\aligned
R'
\big(0,b,a,b\big)
&
\,=\,
R
\big(
0,{\tt Y}(0,b),a,b
\big)
\\
&
\,=\,
R
\big(0,b,a,b\big)
\\
&
\,\equiv\,
0.
\endaligned
\]

Lastly, verify that the second normalization is also
attained, by coming back to~({\ref{R-V-x-y}}),
in which $y$ is replaced by $y := {\tt Y}(x,y')$:
\[
\aligned
0
&
\,\equiv\,
R
\big(
x,{\tt Y}(x,y'),0,{\tt V}(x,{\tt Y}(x,y'))
\big)
\\
&
\,\equiv\,
R
\big(
x,{\tt Y}(x,y'),0,y'
\big)
\\
&
\,\equiv\,
R'
\big(
x,y',0,y'
\big).
\endaligned
\]
In conclusion, the composition of the two changes of 
coordinates employed above is
indeed of the announced form stabilizing
$x' = x$ and $a' = a$.
\endproof

Any system of coordinates $(x,y,a,b)$ in which the considered
submanifold of solutions $\mathcal{M} \subset 
\K_\variablessmall^\NN
\times \K_\parameterssmall^\MM$ has a defining 
function $R$ enjoying:
\leqnomode\usetagform{default}
\begin{align}
\label{R-Q-P-equiv-0}
0
\,\equiv\,
R
\big(
0,b,a,b
\big)
\,\equiv\,
R
\big(
x,y,0,y
\big),
\end{align}
will be called {\sl normalized coordinates}.
What precedes shows that they always exist, for free.

\begin{Corollary}
In normalized coordinates $(x,y,a,b)$ at the
origin, the two graphing functions
$Q$ and $P$ for:
\[
\mathcal{M}
\,=\,
\big\{
y=Q(x,a,b)
\big\}
\,=\,
\big\{
b
\,=\,
P(a,x,y)
\big\},
\]
enjoy the normalization conditions:
\[
Q(0,a,b)
\,\equiv\,
b
\,\equiv\,
Q(x,0,b)
\ \ \ \ \ \ \ \ \ \ \ \ \
\text{and}
\ \ \ \ \ \ \ \ \ \ \ \ \
P(0,x,y)
\,\equiv\,
y
\,\equiv\,
P(a,0,y).
\]
\end{Corollary}

\proof
Treat only the first pair of normalizations for $Q$, the one
for $P$ being similar. In~({\ref{0-R-Q-identically}}),
setting $x := 0$ and then $a := 0$,
we have:
\[
\aligned
0
&
\,\equiv\,
R
\big(
0,Q(0,a,b),a,b
\big)
\ \ \ \ \ \ \ \ \ \ \ \ \
\text{to be compared with}
\ \ \ \ \ \ \ \ \ \ \ \ \
0
\overset{{\eqref{R-Q-P-equiv-0}}}{\,\,\equiv\,\,}
R
\big(
0,b,a,b
\big),
\\
0
&
\,\equiv\,
R
\big(
x,Q(x,0,b),0,b
\big)
\ \ \ \ \ \ \ \ \ \ \ \ \
\text{to be compared with}
\ \ \ \ \ \ \ \ \ \ \ \ \
0
\overset{{\eqref{R-Q-P-equiv-0}}}{\,\,\equiv\,\,}
R
\big(
x,b,0,b
\big),
\endaligned
\]
and the uniqueness property in the implicit function theorem 
when solving $0 = R(0,y,a,b)$ for $y$ and when solving
$0 = R(x,y,0,b)$ also for $y$ yields as wanted:
\[
Q(0,a,b)
\,\equiv\,
b
\ \ \ \ \ \ \ \ \ \ \ \ \
\text{and}
\ \ \ \ \ \ \ \ \ \ \ \ \
Q(x,0,b)
\,\equiv\,
b.
\qedhere
\]
\endproof

\Section{\bf Two Fundamental Foliations 
$\mathcal{F}_\variablessmall$
and $\mathcal{F}_\parameterssmall$ on $\mathcal{M}$}
\label{two-fundamental-foliations}
\HEAD{{\ref{two-fundamental-foliations}}.~{\sf Two Fundamental 
Foliations $\mathcal{F}_\variablessmall$
and $\mathcal{F}_\parameterssmall$ on $\mathcal{M}$}
}{
Joël {\sc Merker} (Orsay)}

As before, working locally, we adopt the convention of almost never
mentioning open (sub)sets, in order to lighten formalism, preserve
clarity, and save symbols.

Equivalences belonging to $\Diff_\variablessmall \times
\Diff_\parameterssmall$ are of the shape:
\[
\big(x,y,a,b)
\,\,\longmapsto\,\,
\big(
f(x,y),g(x,y),\,
\varphi(a,b),\psi(a,b)
\big)
\,=:\,
\big(
x',y',a',b'
\big),
\]
hence they stabilize vertical sets $\{ \variables = \constant\}$
and horizontal sets $\{ \parameters = \constant\}$.
The tangent bundle $T \big( \K_\variablessmall^\NN \times
\K_\parameterssmall^\MM \big)$ to the 
ambient space $\K_\variablessmall^\NN \times
\K_\parameterssmall^\MM$ is therefore equipped with
two invariant subbundles:
\[
\mathcal{T}^\variablessmall
\,:=\,
\Span\,
\frac{\partial}{\partial x}
+
\Span\,
\frac{\partial}{\partial y}
\ \ \ \ \ \ \ \ \ \ \ \ \
\text{and}
\ \ \ \ \ \ \ \ \ \ \ \ \
\mathcal{T}^\parameterssmall
\,:=\,
\Span\,
\frac{\partial}{\partial a}
+
\Span\,
\frac{\partial}{\partial b},
\]
which are trivially (Frobenius) integrable: 
\[
\big[
\Gamma\big(\mathcal{T}^\variablessmall\big),\,
\Gamma\big(\mathcal{T}^\variablessmall\big)
\big]
\,\subset\,
\Gamma\big(\mathcal{T}^\variablessmall\big)
\ \ \ \ \ \ \ \ \ \ \ \ \
\text{and}
\ \ \ \ \ \ \ \ \ \ \ \ \
\big[
\Gamma\big(\mathcal{T}^\parameterssmall\big),\,
\Gamma\big(\mathcal{T}^\parameterssmall\big)
\big]
\,\subset\,
\Gamma\big(\mathcal{T}^\parameterssmall\big),
\]
hence define two foliations:
\[
\mathcal{F}_\variablessmall
\,:=\,
\bigcup_{a,b}\,
\big\{(x,y,a,b)\big\}
\ \ \ \ \ \ \ \ \ \ \ \ \
\text{and}
\ \ \ \ \ \ \ \ \ \ \ \ \
\mathcal{F}_\parameterssmall
\,:=\,
\bigcup_{x,y}\,
\big\{(x,y,a,b)\big\}.
\]

\begin{center}
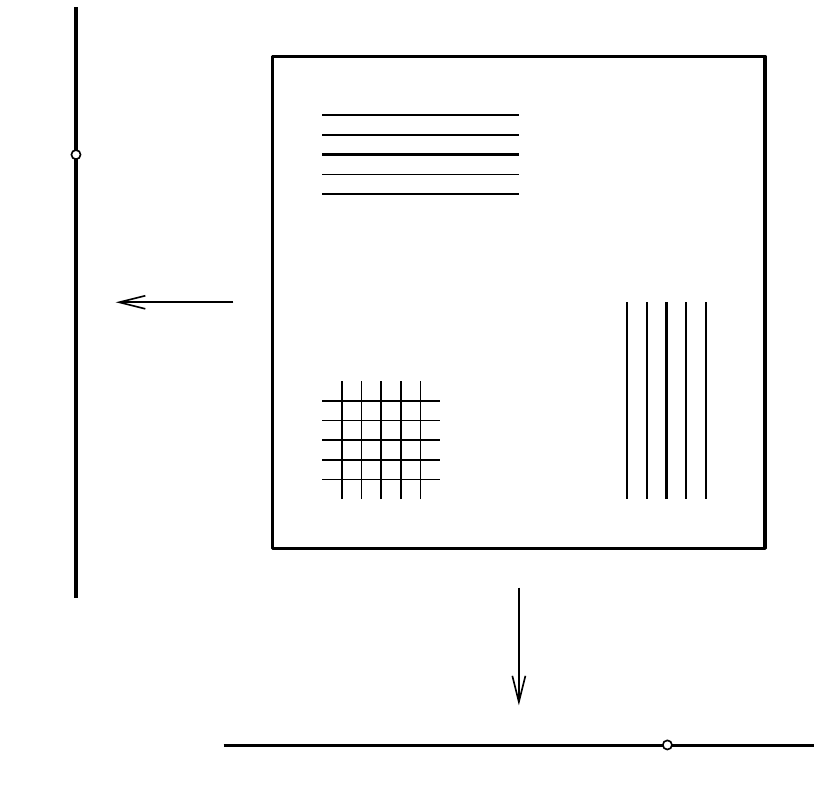
\end{center}

Thanks to the assumption that $\mathcal{M} = \{ y = Q(x,a,b)\} = 
\{ b = P(a,x,y)\}$, these two foliations intersect $\mathcal{M}$
transversely and define two foliations on $\mathcal{M}$,
still denoted:
\[
\mathcal{F}_\variablessmall
\,:=\,
\bigcup_{a,b}\,
\mathcal{Q}_{a,b}
\ \ \ \ \ \ \ \ \ \ \ \ \
\text{and}
\ \ \ \ \ \ \ \ \ \ \ \ \
\mathcal{F}_\parameterssmall
\,=\,
\bigcup_{x,y}\,
\mathcal{P}_{x,y},
\]
whose leaves are defined as:
\[
\aligned
\mathcal{Q}_{a,b}
\,:=\,
\big\{
(x,y)
\colon\,
y
=
Q(x,a,b)
\big\},
\\
\mathcal{P}_{x,y}
\,:=\,
\big\{
(a,b)
\colon\,
b
=
P(a,x,y)
\big\}.
\endaligned
\]
At the level of vector bundles, these $n$-dimensional leaves
$\mathcal{Q}_{a,b}$ and these $m$-dimensional leaves
$\mathcal{P}_{x,y}$ are just integral manifolds of the integrable
subbundles:
\[
T^\variablessmall\mathcal{M}
\,:=\,
T\mathcal{M}
\cap
\mathcal{T}^\variablessmall
\ \ \ \ \ \ \ \ \ \ \ \ \
\text{and}
\ \ \ \ \ \ \ \ \ \ \ \ \
T^\parameterssmall\mathcal{M}
\,:=\,
T\mathcal{M}
\cap
\mathcal{T}^\parameterssmall.
\]

Next, assume that two submanifolds of solutions 
$0 \in \mathcal{M} \subset \K^\NN \times \K^\MM$
and $0' \in \mathcal{M}' \subset {\K'}^\NN \times
{\K'}^\MM$ having the same type $(n,c,m,c)$
and passing through the origin,
represented as:
\[
\aligned
\mathcal{M}
&
\,=\,
\big\{
y
=
Q(x,a,b)
\big\}
\\
&
\,=\,
\big\{
b
=
P(a,x,y)
\big\}
\endaligned
\ \ \ \ \ \ \ \ \ \ \ \ \
\text{and}
\ \ \ \ \ \ \ \ \ \ \ \ \
\aligned
\mathcal{M}'
&
\,=\,
\big\{
y'
=
Q'(x',a',b')
\big\}
\\
&
\,=\,
\big\{
b'
=
P'(a',x',y')
\big\}
\endaligned
\]
are equivalent under an allowed transformation
$(f,g,\varphi,\psi)$. Such an equivalence expresses
as two (equivalent) relations holding for $(x,y,a,b) \in 
\mathcal{M}$:
\[
\aligned
g(x,y)
&
\,=\,
Q'
\big(
f(x,y),\varphi(a,b),\psi(a,b)
\big),
\\
\psi(a,b)
&
\,=\,
P'
\big(
\varphi(a,b),f(x,y),g(x,y)
\big),
\endaligned
\]
that is as power series identities:
\reqnomode\usetagform{EngelLie}
\begin{align}
g\big(x,Q(x,a,b)\big)
&
\,\equiv\,
Q'
\big(
f(x,Q(x,a,b)),\varphi(a,b),\psi(a,b)
\big)
\tag{(\text{\rm in}\,\C\{x,a,b\}^c),}
\\
\psi\big(a,P(a,x,y)\big)
&
\,\equiv\,
P'
\big(
\varphi(a,P(a,x,y)),f(x,y),g(x,y)
\big)
\tag{(\text{\rm in}\,\C\{a,x,y\}^c).}
\end{align}

\begin{center}
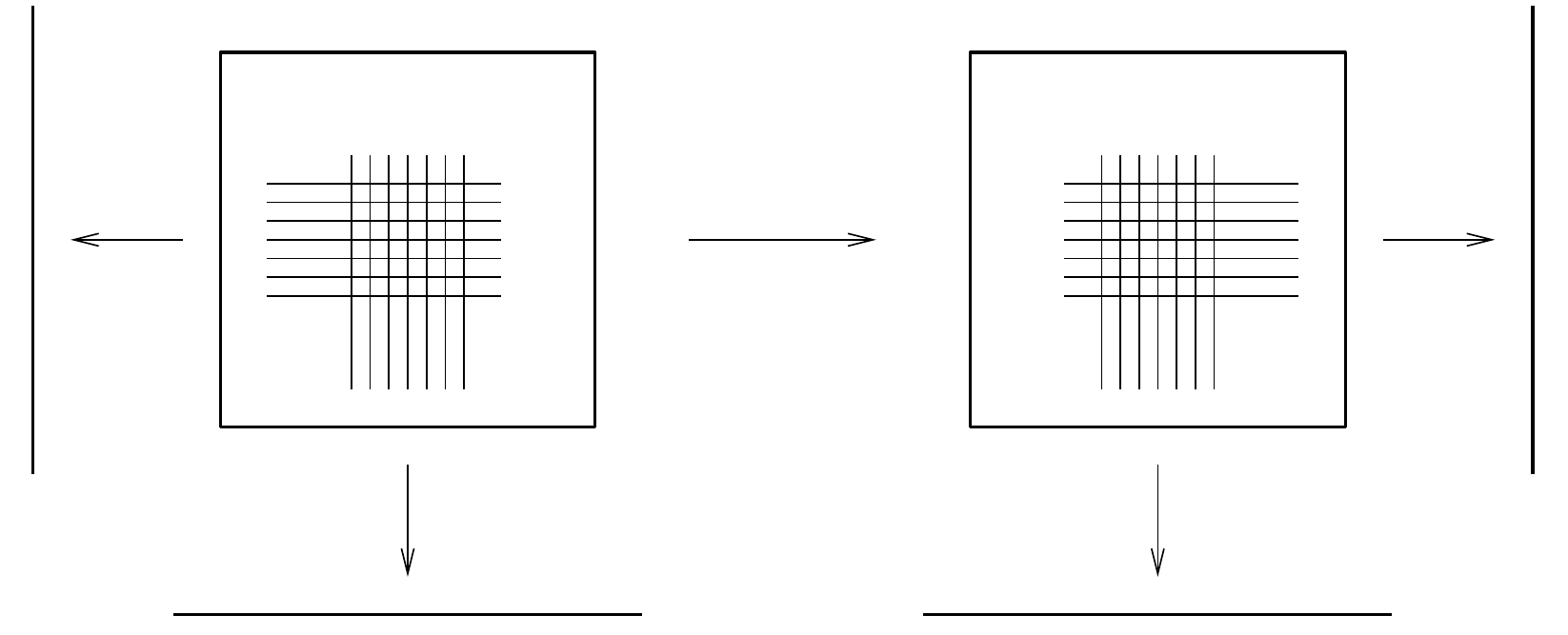
\end{center}

Abbreviating:
\[
F
\,:=\,
(f,g)
\ \ \ \ \ \ \ \ \ \ \ \ \
\text{and}
\ \ \ \ \ \ \ \ \ \ \ \ \
\Phi
\,:=\,
(\varphi,\psi),
\]
these identities can be read off as expressing that:
\[
\aligned
(F,\Phi)
\big(
\mathcal{Q}_{a,b}
\big)
&
\,=\,
\mathcal{Q}_{\Phi(a,b)}',
\\
(F,\Phi)
\big(
\mathcal{P}_{x,y}
\big)
&
\,=\,
\mathcal{P}_{F(x,y)}',
\endaligned
\]
which means that the allowed transformations $(F, \Phi) \in
\Diff_\variablessmall \times \Diff_\parameterssmall$
stabilize the pairs of foliations:
\[
\aligned
(F,\Phi)
\big(\mathcal{F}_\variablessmall\big)
&
\,\subset\,
\mathcal{F}_\variablessmall',
\\
(F,\Phi)
\big(\mathcal{F}_\parameterssmall\big)
&
\,\subset\,
\mathcal{F}_\parameterssmall'.
\endaligned
\]

\Section{\bf Two Systems of Coordinates and Functional Relations}
\label{two-systems-coordinates-functional-relations}
\HEAD{{\ref{two-systems-coordinates-functional-relations}}.~{\sf Two 
Systems of Coordinates and Functional Relations}
}{
Joël {\sc Merker} (Orsay)}

A submanifold of solutions $\mathcal{M}$ of codimension $c$ 
in $\K^n \times \K^c \times \K^m \times \K^c$
naturally comes equipped with two charts represented
by the following two projections:
\[
\aligned
\tau_\parameterssmall
\colon
\ \ \ \ \ \ \ \ \ \ \ \
\mathcal{M}
&
\,\longrightarrow\,
\K^n\times\K^m\times\K^c
\\
(x,y,a,b)
&
\,\longmapsto\,
(x,a,b)
\endaligned
\ \ \ \ \ \ \ \ \ \ \ \ \
\text{and}
\ \ \ \ \ \ \ \ \ \ \ \ \
\aligned
\tau_\variablessmall
\colon
\ \ \ \ \ \ \ \ \ \ \ \
\mathcal{M}
&
\,\longrightarrow\,
\K^m\times\K^n\times\K^c
\\
(x,y,a,b)
&
\,\longmapsto\,
(a,x,y),
\endaligned
\]
which correspond to the arguments $(x,a,b)$ in the graphing function
$Q(x,a,b)$ and to the arguments $(a,x,y)$ in the graphing function
$P(a,x,y)$. The theory will be systematic by playing simultaneously
with these {\em two} systems of coordinates.

Once the implicit function theorem has been applied, it is natural to
set, in the first system of coordinates $(x,a,b)$:
\[
R(x,y,a,b)
\,:=\,
-\,y
+
Q(x,a,b),
\] 
and in the second system of coordinates $(a,x,y)$:
\[
R(x,y,a,b)
\,:=\,
-\,b
+
P(a,x,y).
\]

\begin{FunctionalRelations}
\label{functional-relations-Q-P}
The two graphing functions $Q(x,a,b)$ and $P(a,x,y)$ satisfy
identically:
\reqnomode\usetagform{EngelLie}
\begin{align}
y
&
\,\equiv\,
Q\big(x,a,P(a,x,y)\big),
\notag
\\
b
&
\,\equiv\,
P\big(a,x,Q(x,a,b)\big).
\tag{\qed}
\end{align}
\end{FunctionalRelations}

\Section{\bf Two Fundamental Collections of Vector Fields 
$\mathcal{K}_{x_i}$ and $\mathcal{L}_{a_j}$}
\label{two-collections-vector-fields-K-L}
\HEAD{{\ref{two-collections-vector-fields-K-L}}.~{\sf Two Fundamental 
Collections of Vector Fields $\mathcal{K}_{x_i}$ and 
$\mathcal{L}_{a_j}$}
}{
Joël {\sc Merker} (Orsay)}

From now on, assume $c = 1$.
In the ambient, extrinsic coordinates $(x,y,a,b)$, some natural 
basic vector fields generating the two vector bundles 
$T^\parameterssmall \mathcal{M}$ of rank $m$ and
$T^\variablessmall \mathcal{M}$ of rank $n$ are:
\reqnomode\usetagform{EngelLie}
\begin{align}
\mathcal{L}_{a_j}
&
\,:=\,
\frac{\partial}{\partial a_j}
-
\frac{Q_{a_j}}{Q_b}
(x,a,b)\,
\frac{\partial}{\partial b}
\tag{(1\,\leqslant\,j\,\leqslant\,m),}
\\
\mathcal{K}_{x_i}
&
\,:=\,
\frac{\partial}{\partial x_i}
-
\frac{P_{x_i}}{P_y}
(a,x,y)\,
\frac{\partial}{\partial y}\,
\tag{(1\,\leqslant\,i\,\leqslant\,n).}
\end{align}

One easily convinces oneself that,
in the two systems of coordinates $(x,a,b)$ and $(a,x,y)$,
the pushes-forward
of these vector fields consist in just dropping the
extrinsinc coordinate fields $\frac{\partial}{\partial b}$ and
$\frac{\partial}{\partial y}$:
\[
\aligned
\tau_{\variablessmall\ast}
\bigg(
\frac{\partial}{\partial a_j}
-
\frac{Q_{a_j}}{Q_b}\,
\frac{\partial}{\partial b}
\bigg)
&
\,=\,
\frac{\partial}{\partial a_j},
\\
\tau_{\parameterssmall\ast}
\bigg(
\frac{\partial}{\partial x_i}
-
\frac{P_{x_i}}{P_y}\,
\frac{\partial}{\partial y}
\bigg)
&
\,=\,
\frac{\partial}{\partial x_i}.
\endaligned
\]
So intrinsically on $\mathcal{M}$, generators of $T^\variablessmall
\mathcal{M}$ and of $T^\parameterssmall \mathcal{M}$ write out in two
ways as:
\[
\boxed{\,
\aligned
\mathcal{L}_{a_j}
&
\,=\,
\frac{\partial}{\partial a_j}
-
\frac{Q_{a_j}}{Q_b}\,
\frac{\partial}{\partial b},
\\
\mathcal{K}_{x_i}
&
\,=\,
\frac{\partial}{\partial x_i},
\endaligned
\ \ \ \ \ \ \ \ \ \ \ \ \ \ \ \ \ \ \ \ \ \ \ \ \ \
\text{or}
\ \ \ \ \ \ \ \ \ \ \ \ \ \ \ \ \ \ \ \ \ \ \ \ \ \
\aligned
\mathcal{L}_{a_j}
&
\,=\,
\frac{\partial}{\partial a_j},
\\
\mathcal{K}_{x_i}
&
\,=\,
\frac{\partial}{\partial x_i}
-
\frac{P_{x_i}}{P_y}\,
\frac{\partial}{\partial y},
\endaligned\,}
\]
depending on the choice of coordinates, $(x,a,b)$ or
$(a,x,y)$.

As seen above, the two integrable
distributions generated by the
$\mathcal{L}_{a_j}$ and by
the $\mathcal{K}_{x_i}$ are invariant. 
It is useful to also introduce two {\sl transversal}
vector fields:
\[
\mathcal{T}_b
\,:=\,
\frac{\partial}{\partial b}
\ \ \ \ \ \ \ \ \ \ \ \ \
\text{and}
\ \ \ \ \ \ \ \ \ \ \ \ \
\mathcal{U}_y
\,:=\,
\frac{\partial}{\partial y},
\]
whose directions are {\em not} invariant under allowed
transformations $(F, \Phi) \in 
\Diff_\variablessmall \times \Diff_\parameterssmall$. 

In the coordinates $(x,a,b)$ on $\mathcal{M}$, a frame for the tangent
bundle $T\mathcal{M}$ is:
\[
\big\{
\mathcal{K}_{x_1},\dots,\mathcal{K}_{x_n},\,
\mathcal{L}_{a_1},\dots,\mathcal{L}_{a_m},\,
\mathcal{T}_b
\big\}.
\]
In the other coordinates $(a,x,y)$ on $\mathcal{M}$, a frame
for $T\mathcal{M}$ is:
\[
\big\{
\mathcal{L}_{a_1},\dots,\mathcal{L}_{a_m},\,
\mathcal{K}_{x_1},\dots,\mathcal{K}_{x_n},\,
\mathcal{U}_y
\big\}.
\]

Notice that, by construction:
\[
\big\{
dR=0
\big\}
\,=\,
\Span\,
\big\{
\mathcal{K}_{x_i},\,
\mathcal{L}_{a_j}
\big\}.
\]
a property that is independent of the choice of a defining
function $R$ for the submanifold of solutions, 
{\em cf.}~{\cite{Merker-Pocchiola-Sabzevari-2013-5-CR-II}};
this is clear with either $R := -y + Q(x,a,b)$ or 
$R := -b + P(a,x,y)$, since in this case:
\[
0
\,=\,
dR\big(\mathcal{K}_{x_i}\big)
\,=\,
dR\big(\mathcal{L}_{a_j}\big).
\]
It is also clear that:
\[
dR\big(\mathcal{T}_b\big)
\,\neq\,
0
\ \ \ \ \ \ \ \ \ \ \ \ \
\text{and}
\ \ \ \ \ \ \ \ \ \ \ \ \
dR\big(\mathcal{U}_y\big)
\,\neq\,0,
\]
which can also be seen while making in advance the linear 
normalization~({\ref{linear-normalization-R}}) which implies:
\[
dR
\,=\,
dy+db
+
{\rm O}(1).
\]
Consequently, the analytic differential $1$-forms on $\mathcal{M}$:
\[
\rho
\,:=\,
d\,
\Big(
\frac{R}{dR(\mathcal{T}_b)}
\Big)
\bigg\vert_{\mathcal{M}}
\ \ \ \ \ \ \ \ \ \ \ \ \
\text{and}
\ \ \ \ \ \ \ \ \ \ \ \ \
\sigma
\,:=\,
d
\Big(
\frac{R}{dR(\mathcal{U}_y)}
\Big)
\bigg\vert_{\mathcal{M}}
\]
satisfy:
\[
\rho\big(\mathcal{T}_b\big)
\,\equiv\,
1
\ \ \ \ \ \ \ \ \ \ \ \ \
\text{and}
\ \ \ \ \ \ \ \ \ \ \ \ \
\sigma\big(\mathcal{U}_y\big)
\,\equiv\,
1.
\]

\Section{\bf Nonholonomy Levi Bilinear Forms
$\Levi_\variablessmall$ and $\Levi_\parameterssmall$}
\label{nonholonomy-form-and-dual}
\HEAD{{\ref{nonholonomy-form-and-dual}}.~{\sf 
Nonholonomy Levi Bilinear Forms
$\Levi_\variablessmall$ and $\Levi_\parameterssmall$}
}{
Joël {\sc Merker} (Orsay)}

On any CR manifold, the so-called {\sl Levi form} is a CR-invariant
Hermitian form on the CR bundle. Locally, it is represented by a
square matrix, equal to its transposed conjugate, which means that the
other Hermitian form on the conjugate CR bundle {\em coincides} with
it.

On a submanifold of solutions, by analogy, there is not just one
(generalized) {\sl Levi form}, but {\em two}. 
Similarly, their matrices,
not necessarily square, are transposed of 
one another (up to a nowhere
vanishing factor), so that their ranks and kernels are anyway
essentially the same.

However, we will see later that some higher order jet invariants of
submanifolds of solutions can differ strongly, when seen either from
the variables space, or from the parameters space.

General sections of $T^\variablessmall \mathcal{M}$ and of
$T^\parameterssmall \mathcal{M}$ decompose along the chosen frames
as:
\[
\aligned
\mathcal{K}
&
\,=\,
\mu_1\mathcal{K}_{x_1}
+\cdots+
\mathcal{\mu}_n\mathcal{K}_{x_n},
\\
\mathcal{L}
&
\,=\,
\nu_1\mathcal{L}_{a_1}
+\cdots+
\nu_m\mathcal{L}_{a_m},
\endaligned
\]
in terms of certain arbitrary analytic functions $\mu_i$ and $\nu_j$.

The {\sl Levi form} of $\mathcal{M}$ can be intrinsically
defined as taking Lie brackets modulo the sum of the two
invariant distributions:
\[
\aligned
\Gamma\big(\mathcal{T}^\variablessmall\mathcal{M}\big)
\times
\Gamma\big(\mathcal{T}^\parameterssmall\mathcal{M}\big)
&
\,\,\longrightarrow\,\,
\Gamma\big(T\mathcal{M}\big)
\ \ \
\mod\,
\Big(
\Gamma\big(\mathcal{T}^\variablessmall\mathcal{M}\big)
\oplus
\Gamma\big(\mathcal{T}^\parameterssmall\mathcal{M}\big)
\Big)
\\
\big(
\mathcal{K},\,
\mathcal{L}
\big)
&
\,\,\longmapsto\,\,
\big[
\mathcal{K},\,
\mathcal{L}
\big]
\ \ \
\mod\,
\Big(
\mathcal{T}^\variablessmall\mathcal{M}
\oplus
\mathcal{T}^\parameterssmall\mathcal{M}
\Big),
\endaligned
\] 
but it is preferable, in order to fix ideas about such a moding out,
to employ some differential $1$-forms like $\rho$ and $\sigma$ above
whose kernel distributions represent exactly this direct sum
$\mathcal{T}^\variablessmall\mathcal{M} \oplus
\mathcal{T}^\parameterssmall\mathcal{M}$, the result being
independent, up to a nowhere vanishing factor, of such a choice.

Then {\em two} versions $\Levi_\variablessmall$ and 
$\Levi_\parameterssmall$ of this Levi form exist,
both of which will be useful later.

Firstly, work in coordinates $(x,a,b)$. For
$1 \leqslant i \leqslant n$ and for $1 \leqslant j \leqslant m$,
compute the basic Lie brackets:
\[
\big[
\mathcal{K}_{x_i},
\mathcal{L}_{a_j}
\big]
\,=\,
\bigg[
\frac{\partial}{\partial x_i},\,
\frac{\partial}{\partial a_j}
-
\frac{Q_{a_j}}{Q_b}\,
\frac{\partial}{\partial b}
\bigg]
\,=\,
\Big(
\frac{-Q_b\,Q_{x_ia_j}+Q_{a_j}\,Q_{x_ib}}{
Q_b\,Q_b}
\Big)\,
\frac{\partial}{\partial b}.
\]
With the $1$-form $\rho$ introduced above, compute then:
\[
\aligned
\rho
\Big(
\big[
\mu_1\mathcal{K}_{x_1}+\cdots+\mu_n\mathcal{K}_{x_n},\,\,
\nu_1\mathcal{L}_{a_1}+\cdots+\nu_m\mathcal{L}_{a_m}
\big]
\Big)
&
\,=\,
\sum_{i=1}^n\,
\sum_{j=1}^m\,
\mu_i\,\nu_j\,
\rho
\big(
\big[
\mathcal{K}_{x_i},
\mathcal{L}_{a_j}
\big]
\big)
\\
&
\,=\,
\sum_{i=1}^n\,
\sum_{j=1}^m\,
\mu_i\,\nu_j\,
\frac{-\,Q_bQ_{x_ia_j}+Q_{a_j}Q_{x_ib}}{Q_b\,Q_b},
\endaligned
\]
so that by introducing the $n \times m$ matrix: 
\[
\Levi_\parameterssmall(Q)
\,:=\,
\left(\!
\begin{array}{ccc}
\frac{-Q_bQ_{x_1a_1}+Q_{a_1}Q_{x_1b}}{Q_b\,Q_b}
& \cdots &
\frac{-Q_bQ_{x_1a_m}+Q_{a_m}Q_{x_1b}}{Q_b\,Q_b}
\\
\vdots & \ddots & \vdots
\\
\frac{-Q_bQ_{x_na_1}+Q_{a_1}Q_{x_nb}}{Q_b\,Q_b}
& \cdots &
\frac{-Q_bQ_{x_na_m}+Q_{a_m}Q_{x_nb}}{Q_b\,Q_b}
\end{array}
\!\right),
\]
it comes in terms of the row vectors 
${}^\TT \nu = (\nu_1, \dots, \nu_m)$ and
${}^\TT \mu = (\mu_1, \dots, \mu_n)$: 
\[
\rho
\big(
\big[\mathcal{K},\mathcal{L}\big]
\big)
\,=\,
{}^\TT \mu
\cdot
\Levi_\parameterssmall(Q)
\cdot
\nu.
\]

Secondly, work in coordinates $(a,x,y)$. 
The basic Lie brackets become:
\[
\big[
\mathcal{L}_{a_j},
\mathcal{K}_{x_i}
\big]
\,\,=\,\,
\bigg[
\frac{\partial}{\partial a_j},\,\,
\frac{\partial}{\partial x_i}
-
\frac{P_{x_i}}{P_y}\,
\frac{\partial}{\partial y}
\bigg]
\,\,=\,\,
\Big(
\frac{-P_yP_{a_jx_i}+P_{x_i}P_{a_jy}}{P_y\,P_y}
\Big)\,
\frac{\partial}{\partial y},
\]
hence:
\[
\sigma
\Big(
\big[
\nu_1\mathcal{L}_{a_1}
+\cdots+
\nu_m\mathcal{L}_{a_m},\,\,
\mu_1\mathcal{K}_{x_1}
+\cdots+
\mu_n\mathcal{K}_{x_n}
\big]
\Big)
\,\,=\,\,
\sum_{j=1}^m\,
\sum_{i=1}^n\,
\nu_j\,\mu_i\,
\frac{-P_yP_{a_jx_i}+P_{x_i}P_{a_jy}}{P_y\,P_y},
\]
so that with the $m \times n$ matrix:
\[
\Levi_\variablessmall(P)
\,:=\,
\left(\!
\begin{array}{ccc}
\frac{-P_yP_{a_1x_1}+P_{x_1}P_{a_1y}}{P_y\,P_y}
& \cdots &
\frac{-P_yP_{a_1x_n}+P_{x_n}P_{a_1y}}{P_y\,P_y}
\\
\vdots & \ddots & \vdots
\\
\frac{-P_yP_{a_mx_1}+P_{x_1}P_{a_my}}{P_y\,P_y}
& \cdots &
\frac{-P_yP_{a_mx_n}+P_{x_n}P_{a_my}}{P_y\,P_y}
\end{array}
\!\right),
\]
it comes analogously:
\[
\sigma
\big(
\big[
\mathcal{L},\mathcal{K}
\big]
\big)
\,=\,
{}^\TT\nu
\cdot
\Levi_\variablessmall(P)
\cdot
\mu.
\]

We now need to compare these two matrices
of the two Levi forms. The result
is that they are transposed one to another, up to 
a nowhere vanishing factor.

\begin{Lemma}
\label{Lemma-transfer-from-Q-to-P}
For all indices $1 \leqslant i \leqslant n$ and $1 \leqslant j
\leqslant m$, it holds identically on $\mathcal{M}$:
\[
\frac{-Q_bQ_{x_ia_j}+Q_{a_j}Q_{x_ib}}{
Q_b\,Q_b}
\,\,\equiv\,\,
-\,P_y\,
\bigg(
\frac{-P_yP_{x_ia_j}+P_{x_i}P_{a_jy}}{
P_y\,P_y}
\bigg).
\] 
\end{Lemma}

\proof
Start from the Functional 
Relations~{\ref{functional-relations-Q-P}}, 
rewritten in length as:
\[
\aligned
y
&
\,\equiv\,
Q
\big(
x_1,\dots,x_n,\,a_1,\dots,a_m,\,
P(a_1,\dots,a_m,x_1,\dots,x_n,y)
\big),
\\
b
&
\,\equiv\,
P
\big(
a_1,\dots,a_m,\,x_1,\dots,x_n,\,
Q(x_1,\dots,x_n,a_1,\dots,a_m,b)
\big),
\endaligned
\]
differentiate the first line with respect to 
\green{$x_i$}, to \green{$a_j$}, to \green{$y$}:
\[
\aligned
0
&
\overset{\green{x_i}}{\,\,\equiv\,\,}
Q_{x_i}
+
P_{x_i}\,Q_b,
\\
0
&
\overset{\green{a_j}}{\,\,\equiv\,\,}
Q_{a_j}
+
P_{a_j}\,Q_b,
\\
1
&
\overset{\green{y}}{\,\,\equiv\,\,}
\ \ \ \ \ \ \ \ \ \ \ \
P_y\,Q_b,
\endaligned
\]
and solve:
\[
\aligned
P_{a_j}
&
\,=\,
-\,\frac{Q_{a_j}}{Q_b},
\\
Q_b
&
\,=\,
\frac{1}{P_y}.
\endaligned
\]

Next, differentiate up to order $2$ with respect to
\green{$x_ia_j$}, to \green{$ya_j$}:
\[
\aligned
0
&
\overset{\green{x_ia_j}}{\,\,\equiv\,\,}
Q_{x_ia_j}
+
\widehat{P_{a_j}}Q_{x_ib}
+
P_{x_ia_j}\widehat{Q_b}
+
P_{x_i}\,
\big(
\underline{Q_{a_jb}+P_{a_j}Q_{bb}}
\big),
\\
0
&
\overset{\green{ya_j}}{\,\,\equiv\,\,}
\ \ \ \ \ \ \ \ \ \ \ \ \ \ \
P_y\underline{Q_{a_jb}}
+
P_{a_jy}\widehat{Q_b}
+
\underline{P_{a_j}}P_y\,\underline{Q_{bb}},
\endaligned
\]
solve the underlined terms in the second equation,
replace the result in the first equation,
and replace the wide hats to conclude by multiplying
at the end by $\frac{1}{Q_b} = P_y$:
\begin{align}
0
&
\,\equiv\,
Q_{x_ia_j}
-
\frac{Q_{a_j}}{Q_b}\,
Q_{x_ib}
+
P_{x_ia_j}\,
\frac{1}{P_y}
+
P_{x_i}\,
\Big(
-\,\frac{P_{a_jy}}{P_y}\,
\frac{1}{P_y}
\Big)
\notag
\\
&
\,\equiv\,
\frac{Q_bQ_{x_ia_j}-Q_{a_j}Q_{x_ib}}{Q_b}
+
\frac{P_{x_ia_j}P_y-P_{x_i}P_{a_jy}}{P_y\,P_y}.
\qedhere
\end{align}
\endproof

With this uniform relation between the entries of the two
matrices of the two Levi forms, equal to each other 
up to the nowhere vanishing factors:
\[
-\,P_y
\,\neq\,
0
\ \ \ \ \ \ \ \ \ \ \ \ \
\text{and}
\ \ \ \ \ \ \ \ \ \ \ \ \
-\,Q_b
\,\neq\,
0,
\]
we obtain:
\[
\Levi_\parameterssmall(Q)
\,=\,
-\,P_y\,
{}^\TT
\Levi_\variablessmall(P)
\ \ \ \ \ \ \ \ \ \ \ \ \
\Longleftrightarrow
\ \ \ \ \ \ \ \ \ \ \ \ \
-\,Q_b\,
{}^\TT
\Levi_\parameterssmall(Q)
\,=\,
\Levi_\variablessmall(P).
\]

\Section{\bf Invariance of the Nonholonomy Bilinear Forms
$\Levi_\variablessmall$ and $\Levi_\parameterssmall$}
\label{invariance-nonholonomy-form-and-dual}
\HEAD{{\ref{invariance-nonholonomy-form-and-dual}}.~{\sf 
Invariance of the Nonholonomy Bilinear Forms
$\Levi_\variablessmall$ and $\Levi_\parameterssmall$}
}{
Joël {\sc Merker} (Orsay)}

Consider an equivalence $(F, \Phi) \colon \mathcal{M}
\longrightarrow \mathcal{M}'$ between two submanifolds of
solutions, as in 
Section~{\ref{two-fundamental-foliations}},
hence satisfying:
\[
\aligned
(F,\Phi)_\ast
\big(
T^\variablessmall\mathcal{M}
\big)
&
\,=\,
T^\variablessmall\mathcal{M}',
\\
(F,\Phi)_\ast
\big(
T^\parameterssmall\mathcal{M}
\big)
&
\,=\,
T^\parameterssmall\mathcal{M}'.
\endaligned
\]
As in Section~{\ref{two-systems-coordinates-functional-relations}},
consider two frames for $T \mathcal{M}$ and 
$T \mathcal{M}'$, expressed in terms of graphing functions
$P$, $Q$ in coordinates $(x,y,a,b)$ and in terms
of graphing functions $P'$, $Q'$ in coordinates 
$(x', y', a', b')$:
\[
\aligned
{}
&
\mathcal{K}_{x_1},\dots,\mathcal{K}_{x_1},
\ \ \ \ \ \ \ \ \ \ \ \ \ \ \ \ \ \ \ \ \ \ \ \ \ \
\mathcal{K}_{x_1'}',\dots,\mathcal{K}_{x_n'}',
\\
{}
&
\mathcal{L}_{a_1},\dots,\mathcal{L}_{a_m},
\ \ \ \ \ \ \ \ \ \ \ \ \ \ \ \ \ \ \ \ \ \ \ \ \ \
\mathcal{L}_{a_1'}',\dots,\mathcal{L}_{a_m'}',
\\
{}
&
\mathcal{T}_b,
\ \ \ \ \ \ \ \ \ \ \ \ \ \ \ \ \ \ \ \ \ \ \ \ \ \
\ \ \ \ \ \ \ \ \ \ \ \ \ \ \ \ \ \ \ 
\mathcal{T}_{b'}',
\\
{}
&
\mathcal{U}_y,
\ \ \ \ \ \ \ \ \ \ \ \ \ \ \ \ \ \ \ \ \ \ \ \ \ \
\ \ \ \ \ \ \ \ \ \ \ \ \ \ \ \ \ \ \ 
\mathcal{U}_{y'}'.
\endaligned
\]
Dropping the symbol $(F, \Phi)_\ast^{-1}$ for the
push-forward of vector fields under the inverse map
$\mathcal{M} \longleftarrow \mathcal{M}' \,\,\colon\!
(F,\Phi)^{-1}$, we can therefore write:
\[
\aligned
\mathcal{K}_{x_i}
&
\,=\,
\sum_{1\leqslant i'\leqslant n}\,
X_{i,i'}'\,
\mathcal{K}_{x_{i'}'}',
\\
\mathcal{L}_{a_j}
&
\,=\,
\sum_{1\leqslant j'\leqslant m}\,
A_{j,j'}'\,
\mathcal{L}_{a_{j'}'}',
\endaligned
\]
in terms of certain two invertible $n \times n$ and
$m \times m$ matrices of functions on $\mathcal{M}'$.
Furthermore, there exists a nowhere vanishing 
function $c' \colon \mathcal{M}' \longrightarrow \K$
such that, after (still unwritten) pull-back:
\[
\rho
\,=\,
c'\,\rho'.
\]
Then:
\[
\aligned
\frac{-\,Q_bQ_{x_ia_j}+Q_{a_j}Q_{x_ib}}{Q_b\,Q_b}
&
\,=\,
\rho
\big(
\big[
\mathcal{K}_{x_i},
\mathcal{L}_{a_j}
\big]
\big)
\\
&
\,=\,
c'\,\rho'\,
\bigg(
\Big[
\smallsum{i'}\,
X_{i,i'}'\mathcal{K}_{x_{i'}'}',\,\,
\smallsum{j'}\,
A_{j,j'}'
\mathcal{L}_{a_{j'}'}'
\Big]
\bigg)
\\
&
\,=\,
c'\,\,
\smallsum{i'}\,
\smallsum{j'}\,
X_{i,i'}'\,
\rho'
\big(
\big[
\mathcal{K}_{x_{i'}'}',\,
\mathcal{L}_{a_{j'}'}'
\big]
\big)\,
A_{j,j'}'
\\
&
\,=\,
c'\,\,
\smallsum{i'}\,
\smallsum{j'}\,
X_{i,i'}'\,
\bigg(
\frac{-\,Q_{b'}'Q_{x_{i'}'a_{j'}'}'+
Q_{a_{j'}'}'Q_{x_{i'}'b'}'}{
Q_{b'}'\,Q_{b'}'}
\bigg)\,
A_{j,j'}',
\endaligned
\]
and the result writes in matrix form as:
\[
\footnotesize
\aligned
{}
&
\left(\!
\begin{array}{ccc}
\frac{-Q_bQ_{x_1a_1}+Q_{a_1}Q_{x_1b}}{Q_b\,Q_b}
& \cdots &
\frac{-Q_bQ_{x_1a_m}+Q_{a_m}Q_{x_1b}}{Q_b\,Q_b}
\\
\vdots & \ddots & \vdots
\\
\frac{-Q_bQ_{x_na_1}+Q_{a_1}Q_{x_nb}}{Q_b\,Q_b}
& \cdots &
\frac{-Q_bQ_{x_na_m}+Q_{a_m}Q_{x_nb}}{Q_b\,Q_b}
\end{array}
\!\right)
\,\,=\,\,
\\
&
\,\,=\,\,
c'\,
\left(\!
\begin{array}{ccc}
X_{1,1}' & \cdots & X_{1,n}'
\\
\vdots & \ddots & \vdots
\\
X_{n,1}' & \cdots & X_{n,n}'
\end{array}
\!\right)
\cdot
\left(\!
\begin{array}{ccc}
\frac{-Q_{b'}'Q_{x_1'a_1'}'+Q_{a_1'}'Q_{x_1'b'}'}{Q_{b'}'\,Q_{b'}'}
& \cdots &
\frac{-Q_{b'}'Q_{x_1'a_m'}'+Q_{a_m'}'Q_{x_1'b'}'}{Q_{b'}'\,Q_{b'}'}
\\
\vdots & \ddots & \vdots
\\
\frac{-Q_{b'}'Q_{x_n'a_1'}'+Q_{a_1'}'Q_{x_n'b'}'}{Q_{b'}'\,Q_{b'}'}
& \cdots &
\frac{-Q_{b'}'Q_{x_n'a_m'}'+Q_{a_m'}'Q_{x_n'b'}'}{Q_{b'}'\,Q_{b'}'}
\end{array}
\!\right)
\cdot
\left(\!
\begin{array}{ccc}
A_{1,1}' & \cdots & A_{m,1}'
\\
\vdots & \ddots & \vdots
\\
A_{1,m}' & \cdots & A_{m,m}'
\end{array}
\!\right),
\endaligned
\]
that is to say:
\[
\Levi_\parameterssmall(Q)
\,=\,
c'\cdot X'\cdot
\Levi_\parameterssmall(Q')\cdot
{}^\TT\!A'.
\]

Similarly, with $\sigma = d'\, \sigma'$:
\[
\Levi_\variablessmall(P)
\,=\,
d'\cdot A'\cdot
\Levi_\variablessmall(P')\cdot
{}^\TT\!X'.
\]

\begin{Corollary}
The common rank $r$ with $0 \leqslant r \leqslant \min(n,m)$
of the two matrices $\Levi_\parameterssmall(Q)$ and 
$\Levi_\variablessmall(P)$ of the two Levi forms
of $\mathcal{M}$ is independent of coordinates.\qed
\end{Corollary}

\Section{\bf Normalizations of $\Levi_\parameterssmall$ and of
$\Levi_\variablessmall$ in Local Coordinates}
\label{normalizations-Levi-par-Levi-var}
\HEAD{{\ref{normalizations-Levi-par-Levi-var}}.~{\sf Normalizations 
of $\Levi^\parameterssmall$ and of
$\Levi^\variablessmall$ in Local Coordinates}
}{
Joël {\sc Merker} (Orsay)}

Once this basic invariant is at hand:
\[
r
\,:=\,
\rank\,
\Levi_\parameterssmall(Q)
\,=\,
\rank\,
\Levi_\variablessmall(P),
\]
it is interesting to perform normalizing changes of coordinates
in order to view it in the equations of $\mathcal{M}$.

Assuming from the beginning that coordinates $(x,y,a,b)$
are normal in the sense of 
Proposition~{\ref{Proposition-normal-coordinates}}, 
when expanding in power series:
\[
y
\,=\,
b
+
\sum_{k=0}^\infty\,
Q_k(x,a)\,
b^k
\,=\,
b
+
{\rm O}_{x,a,b}(2),
\] 
it then comes:
\[
0
\,\equiv\,
Q_k(0,a)
\,\equiv\,
Q_k(x,0)
\eqno
{\scriptstyle{(\forall\,k\,\geqslant\,0)}},
\]
whence:
\[
Q_k(x,a)
\,=\,
{\rm O}_{x,a}(2)
\eqno
{\scriptstyle{(\forall\,k\,\geqslant\,0)}}.
\]

By specifying homogeneous terms of degree $2$ in $Q_0(x,a)$,
observing that $b^k\, Q_k(x,a) = b\, {\rm O}_{x,a}(2)$ for all
$k \geqslant 1$, we can write the equation of $\mathcal{M}$ as:
\[
y
\,=\,
b
+
\Lambda(x,a)
+
{\rm O}_{x,a}(3)
+
b\,{\rm O}_{x,a,b}(2),
\]
with a certain bilinear form:
\[
\Lambda(x,a)
\,=\,
\sum_{i=1}^n\,
\sum_{j=1}^m\,
\lambda_{i,j}\,
x_i\,a_j,
\]
which, as is easily devised, {\em represents} 
$\Levi_\parameterssmall(Q)$ 
at the origin.

Indeed, a computation of the entries of this matrix:
\[
\aligned
\frac{-Q_bQ_{x_ia_j}+Q_{a_j}Q_{x_ib}}{Q_b\,Q_b}
&
\,=\,
\frac{-(1+{\rm O}_2)\,(\lambda_{i,j}+{\rm O}_1)
+{\rm O}_1\,{\rm O}_1}{
(1+{\rm O}_1)^2}
\\
&
\,=\,
-\,\lambda_{i,j}
+
{\rm O}_{x,a,b}(1),
\endaligned
\]
confirms that at the origin, up to an innocuous overall minus sign:
\[
\Levi_\parameterssmall(Q)
\,=\,
\big(
-\,
\lambda_{i,j}
\big)_{1\leqslant i\leqslant n}^{1\leqslant j\leqslant m}.
\]

\begin{Lemma}
Any $\K$-bilinear form on $\K_x^n \times \K_a^m$ of rank
$0 \leqslant r \leqslant \min(n,m)$:
\[
\Lambda(x,a)
\,=\,
\sum_{1\leqslant i\leqslant n}\,
\sum_{1\leqslant j\leqslant m}\,
\lambda_{i,j}\,
x_i\,a_j,
\]
can be brought, by means of some {\em linear} transformation
in $\Diff_x \times \Diff_a$, to the form:
\[
x_1a_1
+\cdots+
x_ra_r.
\]
\end{Lemma}

\proof
Rewrite it as:
\[
\Lambda(x,a)
\,=\,
\sum_{i=1}^n\,
x_i\,
\Lambda_i(a_1,\dots,a_r)
\ \ \ \ \ \ \ \ \ \ \ \ \
\text{with}
\ \ \ \ \ \ \ \ \ \ \ \ \
\Lambda_i(a)
\,:=\,
\smallsum{j}\,
\lambda_{i,j}\,a_j.
\]
Among these $n$ linear forms $\Lambda_1, \dots, \Lambda_n$,
a certain maximal
number $r$ are linearly independent, with of course $r$
being the rank in question. After a renumeration, $\Lambda_1, \dots,
\Lambda_r$ {\em are} 
independent, hence can be taken as new coordinates
$a_1 := \Lambda_1$, \dots, $a_r := \Lambda_r$, in terms
of which the remaining linear forms therefore express: 
\[
\Lambda_{r+1}
\big(a_1,\dots,a_r\big),
\,\,\dots\dots\dots,
\Lambda_n
\big(a_1,\dots,a_r\big),
\]
whence a clever (easy) reorganization:
\[
\aligned
\Lambda(x,a)
&
\,=\,
x_1a_1+\cdots+x_ra_r
+
x_{r+1}
\big(
\lambda_{r+1,1}\,a_1
+\cdots+
\lambda_{r+1,r}\,a_r
\big)
\,+
\\
&
\ \ \ \ \ \ \ \ \ \ \ \ \ \ \ \ \ \ \ \ \ \ \ \ \ \ 
\ \ \ \ \ \ \ \ \ \,
+
\cdots\cdots\cdots\cdots\cdots\cdots\cdots\cdots\cdots
+
\\
&
\ \ \ \ \ \ \ \ \ \ \ \ \ \ \ \ \ \ \ \ \ \ \ \ \ \ 
\ \ \ \ \ \ \ \ \ \,
+
x_n
\big(
\lambda_{n,1}\,a_1
+\cdots+
\lambda_{n,r}\,a_r
\big)
\\
&
\,=\,
\big(
x_1
+
\lambda_{r+1,1}x_{r+1}
+\cdots+
\lambda_{n,1}x_n
\big)\,a_1
+\cdots+
\big(
x_r
+
\lambda_{r+1,r}x_{r+1}
+\cdots+
\lambda_{n,r}x_n
\big)\,a_r
\\
&
\,=:\,
x_1a_1
+\cdots+
x_ra_r,
\endaligned
\]
yielding an obvious linear change of coordinates in the $x$-space,
concludes.
\endproof

\begin{Corollary}
\label{normalization-Levi-form-origin}
A submanifold of solutions $\mathcal{M} \subset \K^{n+1} \times
\K^{m+1}$ whose Levi forms at the origin have (always equal) rank 
$0 \leqslant r \leqslant \min(n,m)$ can always be represented,
in suitable normalized coordinates $(x,y,a,b)$, as:
\[
y
\,=\,
b
+
x_1a_1
+\cdots+
x_ra_r
+
{\rm O}_{x,a}(3)
+
b\,{\rm O}_{x,a,b}(2),
\]
and equivalently as:
\[
b
\,=\,
y
-
x_1a_1
-\cdots-
x_ra_r
+
{\rm O}_{a,x}(3)
+
y\,{\rm O}_{a,x,y}(2).
\]
\end{Corollary}

\proof
An application of the implicit
function theorem shows that $P(a,x,y)$ is indeed
of this form once $Q(x,a,b)$ has been normalized.
\endproof

\Section{\bf Levi Forms and Order $1$ Jets of Invariant Leaves}
\label{Levi-form-order-1-jets-Segre}
\HEAD{{\ref{Levi-form-order-1-jets-Segre}}.~{\sf Levi Forms 
and Order $1$ jets of Invariant Leaves}
}{
Joël {\sc Merker} (Orsay)}

In CR geometry, and in integrable {\sc pde} geometry as well, the Levi
form is a too elementary invariant, insufficient to even characterize
some initial data before launching the search for deeper invariants by
means of Cartan's method of equivalence.

A wealth of results in CR geometry, based on finer invariants that
enjoy natural transfer to {\sc pde} geometry, have not yet been set up
in the literature (just a few of them appear
in~{\cite{Merker-2008}}). Let us present a family of fundamental
higher order invariants, and link them with the Levi form(s) before
going further.

Knowing, from Section~{\ref{two-fundamental-foliations}}, that the two
foliations:
\[
\aligned
\mathcal{F}_\variablessmall
&
\,:=\,
\bigcup_{a,b}\,
\big\{
(x,y)
\colon\,
y
=
Q(x,a,b)
\big\},
\\
\mathcal{F}_\parameterssmall
&
\,=\,
\bigcup_{x,y}\,
\big\{
(a,b)
\colon\,
b
=
P(a,x,y)
\big\},
\endaligned
\]
are invariant under equivalences in
$\Diff_\variablessmall \times \Diff_\parameterssmall$,
the main idea is to look at their leaves:
\[
\aligned
\mathcal{Q}_{a,b}
&
\,=\,
\big\{
(x,y)
\colon\,
y
=
Q(x,a,b)
\big\},
\\
\mathcal{P}_{x,y}
&
\,=\,
\big\{
(a,b)
\colon\,
b
=
P(a,x,y)
\big\},
\endaligned
\]
considered as parametrized by $x$ and by $a$, and then, for every
integer $k \geqslant 0$ and every integer $l \geqslant 0$, to
introduce the two fundamental collections of jets-of-leaves maps:
\[
\aligned
J_x^k\mathcal{Q}_{a,b}
&
\,:=\,
\Big(
x,\,\,
\big(
\partial_x^\beta Q(x,a,b)
\big)_{\vert\beta\vert\leqslant k}
\Big),
\\
J_a^l\mathcal{P}_{x,y}
&
\,:=\,
\Big(
a,\,\,
\big(
\partial_a^\gamma
P(a,x,y)
\big)_{\vert\gamma\vert\leqslant l}
\Big),
\endaligned
\]
using multiindex notation:
\[
\aligned
\beta
&
\,=\,
\big(\beta_1,\dots,\beta_n\big)
\,\in\,
\N^n,
\ \ \ \ \ \ \ \ \ \ \ \ \
\vert\beta\vert
\,:=\,
\beta_1+\cdots+\beta_n,
\\
\gamma
&
\,=\,
\big(\gamma_1,\dots,\gamma_m\big)
\,\in\,
\N^m,
\ \ \ \ \ \ \ \ \ \ \ \ \
\vert\gamma\vert
\,:=\,
\gamma_1+\cdots+\gamma_m.
\endaligned
\]

In the next 
Section~{\ref{higher-order-jets-invariant-leaves}}, we
will show rigorously that these two maps {\em are} invariants of
submanifolds of solutions considered modulo 
transformations in the pseudogroup
$\Diff_\variablessmall \times \Diff_\parameterssmall$, but before
doing this, to make a natural transition with what precedes, 
let us look at the special case
of order $1$ jets:
\[
k
\,=\,
1
\,=\,
l.
\]

Since the first entry $x$ of the order $1$ jet map
$J_x^1 \mathcal{Q}_{a,b}$ does not depend on $(a,b)$, 
its rank amounts to the rank
of the `submap':
\[
J_\variablessmall^1Q
\colon\ \ \ 
(a,b)
\,\,\longmapsto\,\,
\big(
Q,\,\,
Q_{x_1},\dots,Q_{x_n}
\big)
(x,a,b)
\]
---\,\,in which $x$ is fixed\,\,---\,\,and 
whose Jacobian $(n+1) \times (m+1)$ matrix is just:
\[
\Jac\,
\big(
J_\variablessmall^1Q
\big)
\,=\,
\left(\!
\begin{array}{cccc}
Q_{a_1} & \cdots & Q_{a_m} & Q_b
\\
Q_{x_1a_1} & \cdots & Q_{x_1a_m} & Q_{x_1b}
\\
\vdots & \ddots & \vdots & \vdots
\\
Q_{x_na_1} & \cdots & Q_{x_na_m} & Q_{x_nb}
\end{array}
\!\right).
\]

\begin{Question}
{\sl Can this $\Jac\, (J_\variablessmall^1Q)$ 
be compared with the matrix:}
\[
\Levi_\parameterssmall(Q)
\,=\,
\left(\!
\begin{array}{ccc}
\frac{-Q_bQ_{x_1a_1}+Q_{a_1}Q_{x_1b}}{Q_b\,Q_b}
& \cdots &
\frac{-Q_bQ_{x_1a_m}+Q_{a_m}Q_{x_1b}}{Q_b\,Q_b}
\\
\vdots & \ddots & \vdots
\\
\frac{-Q_bQ_{x_na_1}+Q_{a_1}Q_{x_nb}}{Q_b\,Q_b}
& \cdots &
\frac{-Q_bQ_{x_na_m}+Q_{a_m}Q_{x_nb}}{Q_b\,Q_b}
\end{array}
\!\right)
\text{\bf ?}
\]
\end{Question}

Yes! It suffices, since $Q_b \neq 0$ vanishes nowhere
by assumption, to divide all entries of 
$\Jac\big( J_\variablessmall^1 Q \big)$ by $Q_b$ 
and then to perform obvious row operations 
which leave the rank unchanged:
\[
\aligned
\Jac\,\big(J_\variablessmall^1Q\big)
&
\,\,\,\longmapsto\,\,\,
Q_b\,
\left(\!
\begin{array}{cccc}
\frac{Q_{a_1}}{Q_b} & \cdots & \frac{Q_{a_m}}{Q_b} & 1
\\
\frac{Q_{x_1a_1}}{Q_b} & \cdots & \frac{Q_{x_1a_m}}{Q_b} 
& \frac{Q_{x_1b}}{Q_b}
\\
\vdots & \ddots & \vdots & \vdots
\\
\frac{Q_{x_na_1}}{Q_b} & \cdots & \frac{Q_{x_na_m}}{Q_b} 
& \frac{Q_{x_nb}}{Q_b}
\end{array}
\!\right)
\\
&
\,\,\,\longmapsto\,\,\,
Q_b\,
\left(\!
\begin{array}{cccc}
\frac{Q_{a_1}}{Q_b} & \cdots & \frac{Q_{a_m}}{Q_b} & 1
\\
\frac{Q_{x_1a_1}}{Q_b}
-
\frac{Q_{a_1}}{Q_b}\,\frac{Q_{x_1b}}{Q_b}
& \cdots & 
\frac{Q_{x_1a_m}}{Q_b}
-
\frac{Q_{a_m}}{Q_b}\,\frac{Q_{x_1b}}{Q_b}
& 0
\\
\vdots & \ddots & \vdots & \vdots
\\
\frac{Q_{x_na_1}}{Q_b}
-
\frac{Q_{a_1}}{Q_b}\,\frac{Q_{x_nb}}{Q_b}
& \cdots & 
\frac{Q_{x_na_m}}{Q_b}
-
\frac{Q_{a_m}}{Q_b}\,\frac{Q_{x_nb}}{Q_b}
& 0
\end{array}
\!\right)
\\
&
\ \ \ \ \ 
\,=\,
Q_b\,
\left(\!
\begin{array}{cc}
\ast & 1
\\
-\,\Levi_\variablessmall(Q) & 0
\end{array}
\!\right),
\endaligned
\]
to realize that the rank of this order $1$ jet map
coincides with the rank of the Levi form!

Of course, a similar procedure can also be applied to the map:
\[
J_\parameterssmall^1P
\colon\ \ \
(x,y)
\,\,\longmapsto\,\,
\big(
P,\,\,
P_{a_1},\,\,
\dots,\,\,
P_{a_m}
\big)
(a,x,y)
\]
having $(m+1) \times (n+1)$ Jacobian matrix:
\[
\aligned
\Jac\,
\big(
J_\parameterssmall^1P
\big)
&
\,=\,
\left(\!
\begin{array}{cccc}
P_{x_1} & \cdots & P_{x_n} & P_y
\\
P_{a_1x_1} & \cdots & P_{a_1x_n} & P_{a_1y}
\\
\vdots & \ddots & \vdots & \vdots
\\
P_{a_mx_1} & \cdots & P_{a_mx_n} & P_{a_my}
\end{array}
\!\right)
\\
&
\longmapsto\,
P_y\,
\left(\!
\begin{array}{cc}
\ast & 1
\\
-\,\Levi_\parameterssmall(P) & 0
\end{array}
\!\right).
\endaligned
\]

\begin{Corollary}
\label{Corollary-jets-1-identical-Levi-form}
At every point $(x,y,a,b) \in \mathcal{M}$:
\[
\rank\,
\big(
J_\variablessmall^1Q
\big)
\,=\,
1
+
\rank\,
\big(
\Levi_\variablessmall(Q)
\big)
\,=\,
1
+
\rank\,
\big(
\Levi_\parameterssmall(P)
\big)
\,=\,
\rank\,
\big(
J_\parameterssmall^1P
\big),
\]
or equivalently:
\reqnomode\usetagform{EngelLie}
\begin{align}
\rank\,\big(J_x^1\mathcal{Q}_{a,b}\big)
&
\,=\,
n+1
+
\rank\,
\big(
\Levi_\variablessmall(Q)
\big),
\notag
\\
\rank\,\big(J_a^1\mathcal{P}_{x,y}\big)
&
\,=\,
m+1
+
\rank\,
\big(
\Levi_\parameterssmall(P)
\big).
\tag{\qed}
\end{align}
\end{Corollary}

\Section{\bf Higher 
Order Jets of the Pair of Invariant Leaves}
\label{higher-order-jets-invariant-leaves}
\HEAD{{\ref{higher-order-jets-invariant-leaves}}.~{\sf Higher Order 
Jets of the Pair of Invariant Leaves}
}{
Joël {\sc Merker} (Orsay)}

\footnotetext[5]{\,
In a CR context, 
partial aspects of the topics presented in
this Section~{\ref{higher-order-jets-invariant-leaves}}
and in the next 
Sections~{\ref{invariancies-jet-ranks}},
{\ref{two-finite-nondegeneracies}},
{\ref{generic-degeneracies}},
{\ref{local-Lie-group-structure}},
{\ref{generic-2-2-1-para-CR-structures}},
have already been studied
in~{\cite{Merker-2002, Merker-Porten-2006,
Merker-Pocchiola-Sabzevari-2013-5-CR-II}}.}

Having justified the interest of jets by linking them with the
(too) classical concept(s) of Levi form(s), let us now study our two 
favorite jet maps\footnotemark[5]
of arbitrary orders:
\[
\aligned
J_x^k\mathcal{Q}_{a,b}
&
\,:=\,
\Big(
x,\,\,
\big(
\partial_x^\beta Q(x,a,b)
\big)_{\vert\beta\vert\leqslant k}
\Big),
\\
J_a^l\mathcal{P}_{x,y}
&
\,:=\,
\Big(
a,\,\,
\big(
\partial_a^\gamma
P(a,x,y)
\big)_{\vert\gamma\vert\leqslant l}
\Big),
\endaligned
\]
which, for $k = 0 = l$, reduce plainly to:
\[
\aligned
(x,a,b)
&
\,\longmapsto\,
\big(
x,\,Q(x,a,b)
\big),
\\
(a,x,y)
&
\,\longmapsto\,
\big(
a,P(a,x,y)
\big).
\endaligned
\]
The goal is to understand in a precise
manner their invariancy under the pseudogroup
of allowed equivalences:
\[
\aligned
(x,y,a,b)
&
\,\,\longmapsto\,\,
\big(
f(x,y),\,g(x,y),\,\varphi(a,b),\,\psi(a,b)
\big)
\\
&
\ \ \ \ \
=:\,
\big(x',y',a',b'\big),
\endaligned
\]
which send $\mathcal{M}$ into another $\mathcal{M}'$
represented by similar equivalent equations:
\[
\aligned
y'
&
\,=\,
Q'(x',a',b'),
\\
b'
&
\,=\,
P'(a',x',y').
\endaligned
\] 

For jets of order $k = 0 = l$, the task is easy.

\begin{Observation}
The two maps:
\[
\aligned
{}
&
\xymatrix{
(x,y)
\ar[rr]^{R_{f,g}^0\,\,\,\,\,\,\,\,\,\,\,\,\,\,\,\,\,\,\,\,\,\,}
&
&
\big(f(x,y),g(x,y)\big),
}
\\
{}
&
\xymatrix{
(a,b)
\ar[rr]^{S_{\varphi,\psi}^0\,\,\,\,\,\,\,\,\,\,\,\,\,\,\,\,\,\,\,\,\,\,}
&
&
\big(\varphi(a,b),\psi(a,b)\big),
}
\endaligned
\]
make commutative the two diagrams:
\[
\xymatrix{
\mathcal{M}
\ar[rr]^{(f,g,\varphi,\psi)}
\ar[d]_{J_x^0\mathcal{Q}_{a,b}}
&
&
\mathcal{M}'
\ar[d]^{J_{x'}^0\mathcal{Q}_{a',b'}'}
\\
\K^{n+1}
\ar[rr]_{R_{f,g}^0}
&
&
{\K'}^{n+1},
}
\ \ \ \ \ \ \ \ \ \ \ \ \ \ \ \ \ \ \ \ \ \ \ \ \ \
\xymatrix{
\mathcal{M}
\ar[rr]^{(\varphi,\psi,f,g)}
\ar[d]_{J_a^0\mathcal{P}_{x,y}}
&
&
\mathcal{M}'
\ar[d]^{J_{a'}^0\mathcal{P}_{x',y'}'}
\\
\K^{m+1}
\ar[rr]_{S_{\varphi,\psi}^0}
&
&
{\K'}^{m+1}.
}
\]
\end{Observation}

\proof
By symmetry, we can treat only the first diagram. 
The up$\,\circ\,$right composition is:
\[
\xymatrix{
(x,a,b)
\ar[rr]
&
&
\big(
f(x,Q(x,a,b)),\varphi(a,b),\psi(a,b)
\big)
\ar[d]
\\
&
&
\Big(
f(x,Q(x,a,b)),\,
Q'\big(f(x,Q(x,a,b)),
\varphi(a,b),\psi(a,b)
\big)
\Big),
}
\]
while the down\,$\circ$\,right one is:
\[
\xymatrix{
(x,a,b)
\ar[d]
&
&
\\
\big(x,Q(x,a,b)\big)
\ar[rr]
&
&
\big(
f(x,Q(x,a,b)),\,
g(x,Q(x,a,b))
\big),
}
\]
and the coincidence of the
two bottom-right results comes instantly from the assumption
$(f,\varphi,\psi)(\mathcal{M}) \subset \mathcal{M}'$
which writes out as the following identity in $\C\{x, a, b\}$:
\[
g\big(x,Q(x,a,b)\big)
\,\equiv\,
Q'
\Big(
f\big(x,Q(x,a,b)\big),\,
\varphi(a,b),\psi(a,b)
\Big).
\qedhere
\]
\endproof

Similar pairs of commutative diagrams exist for jets $k$ and $l$
of arbitrary orders. To present them, some preliminaries are needed.

\begin{Lemma}
\label{Lemma-invertibility-CR-horizontal-map}
For every equivalence $(f,g,\varphi,\psi) \in 
\Diff_\variablessmall \times \Diff_\parameterssmall$
between two submanifolds of solutions $\mathcal{M} 
\longrightarrow \mathcal{M}'$, one has the two
everywhere nonvanishings:
\[
\aligned
0
&
\,\neq\,
\Delta(f,Q)
\,:=\,
\det\,
\Big(
\big(
f_{i',x_i}
+
Q_{x_i}\,f_{i',y}
\big)_{1\leqslant i'\leqslant n}^{
1\leqslant i\leqslant n}
\Big),
\\
0
&
\,\neq\,
\Delta(\varphi,P)
\,:=\,
\det\,
\Big(
\big(
\varphi_{j',a_j}
+
P_{a_j}\,\varphi_{j',b}
\big)_{1\leqslant j'\leqslant n}^{
1\leqslant j\leqslant n}
\Big).
\endaligned
\]
\end{Lemma}

\proof
By symmetry, we can focus on
just $\Delta(f,Q)$. 

We know that the restriction of $(F, \Phi)$ to the $n$-dimensional
leaf in $\mathcal{M}$:
\[
\mathcal{Q}_{a,b}
\,=\,
\big\{
(x,y)
\colon\,
y
=
Q(x,a,b)
\big\},
\]
maps it into the corresponding $n$-dimensional leaf in $\mathcal{M}'$:
\[
(F,\Phi)
\big(\mathcal{Q}_{a,b}\big)
\,\,\subset\,\,
\mathcal{Q}_{\varphi(a,b),\psi(a,b)}'.
\]
It is well known that the restriction of any (local) diffeomorphism
to a submanifold becomes a diffeomorphism onto its (submanifold)
image, hence we receive an analytic diffeomorphism:
\[
(F,\Phi)
\Big\vert_{\mathcal{Q}_{a,b}}
\colon\ \ \
\mathcal{Q}_{a,b}
\overset{\sim}{\,\,\longrightarrow\,\,}
\mathcal{Q}_{\varphi(a,b),\psi(a,b)}'.
\]

But since $x'$ is a natural coordinate on all the
$\mathcal{Q}_{a',b'}'$, which are graphed as $y' = Q'(x',a',b')$,
and since $x$ parametrizes all the $\mathcal{Q}_{a,b}$, 
the restricted map in question identifies with the
bottom map below:
\[
\xymatrix{
\big(x,Q(x,a,b),a,b\big)
\ar[rr]
&
&
\big(
f(x,Q(x,a,b)),\,g(x,Q(x,a,b)),\,\varphi(a,b),\,\psi(a,b)
\big)
\ar[d]
\\
x
\ar[rr]
\ar[u]
&
&
f(x,Q(x,a,b)),
}
\]
which is therefore a (local) diffeomorphism $\K^n \longrightarrow
\K^n$, hence has nowhere vanishing Jacobian determinant,
easily seen to be precisely $\Delta(f,Q)$, because:
\[
{\textstyle{\frac{\partial}{\partial x_i}}}\,
\big[
f_{i'}
\big(x,Q(x,a,b)\big)
\big]
\,=\,
f_{i',x_i}
+
Q_{x_i}\,f_{i',y}.
\qedhere
\]
\endproof

The fact that $(f,g,\varphi,\psi) \colon \mathcal{M}
\overset{\sim}{\longrightarrow} \mathcal{M}'$ is an equivalence
expresses as two pairs of equivalent identities:
\reqnomode\usetagform{EngelLie}
\begin{align}
g\big(x,Q(x,a,b)\big)
&
\,\equiv\,
Q'
\Big(
f\big(x,Q(x,a,b)\big),\,
\varphi(a,b),\psi(a,b)
\Big)
\tag{(\text{\rm in}\,\C\{x,a,b\}),}
\\
g(x,y)
&
\,\equiv\,
Q'\Big(
f(x,y),\,
\varphi\big(a,P(a,x,y)\big),\,
\psi\big(a,P(a,x,y)\big)
\Big)
\tag{(\text{\rm in}\,\C\{a,x,y\}),}
\\
\psi\big(
a,P(a,x,y)\big)
&
\,\equiv\,
P'\Big(
\varphi\big(
a,P(a,x,y)\big),\,
f(x,y),g(x,y)
\Big)
\Big)
\tag{(\text{\rm in}\,\C\{a,x,y\}),}
\\
\psi(a,b)
&
\,=\,
P'\Big(
\varphi(a,b),\,
f\big(x,Q(x,a,b)
\big),\,
g\big(x,Q(x,a,b)
\big)
\Big)
\tag{(\text{\rm in}\,\C\{x,a,b\}).}
\end{align}

By repeatedly differentiating these identities and then reorganizing
the outcomes, we will soon be in a position
to realize the invariancy of the two jet maps under
the form of an appropriate statement. For that purpose, let us
introduce some independent jet-coordinates:
\[
\big(
y_{x^\beta}
\big)_{\vert\beta\vert\leqslant k}
\ \ \ \ \ \ \ \ \ \ \ \ \
\text{and}
\ \ \ \ \ \ \ \ \ \ \ \ \
\big(
b_{a^\gamma}
\big)_{\vert\gamma\vert\leqslant l},
\]
corresponding to $y$ and $b$ being considered as functions of $x$
and $a$, a point of view which will also be useful later on
when we will introduce partial differential equations associated to 
submanifolds of solutions. 
The two jet maps 
take therefore their values into the two jet spaces:
\[
\aligned
J_x^k\mathcal{Q}_{a,b}
&
\,\,\in\,\,
\K^{n+\frac{(n+k)!}{n!\,k!}}
\,\,=\,\,
\big\{
\big(
x,\,(y_{x^\beta})_{\vert\beta\vert\leqslant k}
\big)
\big\},
\\
J_a^l\mathcal{P}_{x,y}
&
\,\,\in\,\,
\K^{m+\frac{(m+l)!}{m!\,l!}}
\,\,=\,\,
\big\{
\big(
a,\,
(b_{a^\gamma})_{\vert\gamma\vert\leqslant l}
\big)
\big\}.
\endaligned
\]
Lastly, let us employ abbreviated notations for jets:
\[
\aligned
j_x^kQ
&
\,:=\,
\big(
\partial_x^\beta Q
\big)_{\vert\beta\vert\leqslant k},
\ \ \ \ \ \ \ \ \ \ \ \ \ \ \ \ \ \ \ \ \ \ \ \ \ \
\ \ \ \ \ \ \ \ \ \ \ \ \
j_a^lP
\,:=\,
\big(
\partial_a^\gamma P
\big)_{\vert\gamma\vert\leqslant l},
\\
j_{a,b}^k\varphi
&
\,:=\,
\big(
\partial_a^\gamma\partial_b^l\varphi
\big)_{\vert\gamma\vert+l\leqslant k},
\ \ \ \ \ \ \ \ \ \ \ \ \ \ \ \ \ \ \ \ \ \ \ \ \ \ 
\ \ \ \ \ \
j_{x,y}^lf
\,:=\,
\big(
\partial_x^\beta\partial_y^kf
\big)_{\vert\beta\vert+k\leqslant l}
\\
j_{a,b}^k\psi
&
\,:=\,
\big(
\partial_a^\gamma\partial_b^l\psi
\big)_{\vert\gamma\vert+l\leqslant k},
\ \ \ \ \ \ \ \ \ \ \ \ \ \ \ \ \ \ \ \ \ \ \ \ \ \ 
\ \ \ \ \ \
j_{x,y}^lg
\,:=\,
\big(
\partial_x^\beta\partial_y^kg
\big)_{\vert\beta\vert+k\leqslant l}.
\endaligned
\]

\begin{Theorem}
\label{Thm-transfer-jets}
For each multiindex $\beta' \in \N^n$ with $1 \leqslant \vert
\beta'\vert$, and for each multiindex $\gamma' \in \N^m$ with $1
\leqslant \vert \gamma' \vert$, there exist two rational expressions
satisfying identically in $\C\{x, a, b\}$ and in $\C \{a, x, y\}$:
\[
\aligned
\frac{\Polynomial_{\beta'}
\big(
\big(
\partial_x^\beta Q\big)_{
1\leqslant\vert\beta\vert\leqslant\vert\beta'\vert},\,
\big(
\partial_x^\beta\partial_y^lf
\big)_{1\leqslant\vert\beta\vert+l\leqslant\vert\beta'\vert},\,
\big(
\partial_x^\beta\partial_y^lg
\big)_{1\leqslant\vert\beta\vert+l\leqslant\vert\beta'\vert}
\big)}{
\big[\det\,(f_{i',x_i}+Q_{x_i}f_{i',y})]^{2\vert\beta'\vert+1}}
\,\,\equiv\,\,
Q_{{x'}^{\beta'}}'
\big(
f,\,\varphi,\psi
\big),
\\
\frac{\Polynomial_{\gamma'}
\big(
\big(
\partial_a^\beta P\big)_{
1\leqslant\vert\gamma\vert\leqslant\vert\gamma'\vert},\,
\big(
\partial_a^\gamma\partial_b^k\varphi
\big)_{1\leqslant\vert\gamma\vert+k\leqslant\vert\gamma'\vert},\,
\big(
\partial_a^\gamma\partial_b^k\psi
\big)_{1\leqslant\vert\gamma\vert+k\leqslant\vert\gamma'\vert}
\big)}{
\big[\det\,(\varphi_{j',a_j}+P_{a_j}\varphi_{j',b})]^{
2\vert\gamma'\vert+1}}
\,\,\equiv\,\,
P_{{a'}^{\gamma'}}'
\big(
\varphi,f,g
\big).
\endaligned
\]

Moreover, for any two integers $k \geqslant 0$ and
$l \geqslant 0$, the two maps:
\[
\aligned
R_{f,g}^k
\colon\ \ \ \ \
\K^{n+\frac{(n+k)!}{n!\,k!}}
&
\xrightarrow[{\rule[0pt]{40pt}{0pt}}]{}
\K^{n+\frac{(n+k)!}{n!\,k!}},
\\
S_{\varphi,\psi}^l
\colon\ \ \ \ \
\K^{m+\frac{(m+l)!}{m!\,l!}}
&
\xrightarrow[{\rule[0pt]{40pt}{0pt}}]{}
\K^{m+\frac{(m+l)!}{m!\,l!}},
\endaligned
\]
defined by:
\[
\!\!\!\!\!\!\!\!\!\!\!\!\!\!\!
\footnotesize
\aligned
{}
&
R_{f,g}^k\,
\Big(
x,y,\,
\big(
y_{x^\beta}
\big)_{1\leqslant\vert\beta\vert\leqslant k}
\Big)
\,:=\,
\\
&
\,:=\,
\bigg(
f(x,y),g(x,y),\,\,
\bigg\{
\frac{\Polynomial_{\beta'}
\big(
\big(
y_{x^\beta}
\big)_{
1\leqslant\vert\beta\vert\leqslant\vert\beta'\vert},\,
\big(
\partial_x^\beta\partial_y^lf(x,y)
\big)_{1\leqslant\vert\beta\vert+l\leqslant\vert\beta'\vert},\,
\big(
\partial_x^\beta\partial_y^lg(x,y)
\big)_{1\leqslant\vert\beta\vert+l\leqslant\vert\beta'\vert}
\big)}{
\big[\det\,(f_{i',x_i}(x,y)+y_{x_i}f_{i',y}(x,y))]^{2\vert\beta'\vert+1}}
\bigg\}_{1\leqslant\vert\beta'\vert\leqslant k}
\Bigg),
\\
{}
&
S_{\varphi,\psi}^l\,
\Big(
a,b,\,
\big(
b_{a^\gamma}
\big)_{1\leqslant\vert\gamma\vert\leqslant l}
\Big)
\,:=\,
\\
&
\,:=\,
\bigg(
\varphi(a,b),\psi(a,b),\,\,
\bigg\{
\frac{\Polynomial_{\gamma'}
\big(
\big(
b_{a^\gamma}
\big)_{
1\leqslant\vert\beta\vert\leqslant\vert\beta'\vert},\,
\big(
\partial_a^\gamma\partial_b^k\varphi(a,b)
\big)_{1\leqslant\vert\gamma\vert+k\leqslant\vert\gamma'\vert},\,
\big(
\partial_a^\gamma\partial_b^k\psi(a,b)
\big)_{1\leqslant\vert\gamma\vert+k\leqslant\vert\gamma'\vert}
\big)}{
\big[\det\,(\varphi_{j',a_j}(a,b)+b_{a_j}\varphi_{j',b}(a,b))]^{
2\vert\gamma'\vert+1}}
\bigg\}_{1\leqslant\vert\gamma'\vert\leqslant l}
\Bigg),
\endaligned
\]
make commutative the two diagrams:
\[
\aligned
\xymatrix{
\mathcal{M}
\ar[rr]^{(f,g,\varphi,\psi)}
\ar[d]_{J_x^k\mathcal{Q}_{a,b}}
&
&
\mathcal{M}'
\ar[d]^{J_{x'}^k\mathcal{Q}_{a',b'}'}
\\
\K^{n+\frac{(n+k)!}{n!\,k!}}
\ar[rr]_{R_{f,g}^k}
&
&
{\K'}^{n+\frac{(n+k)!}{n!\,k!}},
}
\ \ \ \ \ \ \ \ \ \ \ \ \
\ \ \ \ \ \ \ \ \ \ \ \ \
\xymatrix{
\mathcal{M}
\ar[rr]^{(\varphi,\psi,f,g)}
\ar[d]_{J_a^l\mathcal{P}_{x,y}}
&
&
\mathcal{M}'
\ar[d]^{J_{a'}^l\mathcal{P}_{x',y'}'}
\\
\K^{m+\frac{(m+l)!}{m!\,l!}}
\ar[rr]_{S_{\varphi,\psi}^l}
&
&
{\K'}^{m+\frac{(m+l)!}{m!\,l!}}.
}
\endaligned
\]
\end{Theorem}

\proof
By symmetry, we can focus on the first statement, which 
has already been explained above for $k = 0$. 

For simplicity, let us treat the case $k = 1$. We start from the
identity:
\[
g\big(x,Q(x,a,b)\big)
\,\equiv\,
Q'\Big(
f\big(x,Q(x,a,b)\big),\,\varphi(a,b),\psi(a,b)
\Big),
\]
which we differentiate with respect to the $x_i$ for all $1 \leqslant i
\leqslant n$:
\leqnomode\usetagform{default}
\begin{align}
\label{jet-1-g-Q-Q-prime}
g_{x_i}
+
Q_{x_i}\,g_y
\,\equiv\,
\sum_{1\leqslant i'\leqslant n}\,
\big(
f_{i',x_i}+Q_{x_i}\,f_{i',y}
\big)\,
Q_{x_{i'}'}'.
\end{align}
This amounts to apply the already seen first-order differentiation
operators:
\[
\mathcal{K}_{x_i}
\,:=\,
\frac{\partial}{\partial x_i}
+
Q_{x_i}\,
\frac{\partial}{\partial y},
\]
namely in matrix form:
\[
\left(\!
\begin{array}{c}
\mathcal{K}_{x_1}(g)
\\
\vdots
\\
\mathcal{K}_{x_n}(g)
\end{array}
\!\right)
\,\,=\,\,
\left(\!
\begin{array}{ccc}
\mathcal{K}_{x_1}(f_1) & \cdots & \mathcal{K}_{x_1}(f_n)
\\
\vdots & \ddots & \vdots
\\
\mathcal{K}_{x_n}(f_1) & \cdots & \mathcal{K}_{x_n}(f_n)
\end{array}
\!\right)
\left(\!
\begin{array}{c}
Q_{x_1'}'
\\
\vdots
\\
Q_{x_n'}'
\end{array}
\!\right).
\]
Since Lemma~{\ref{Lemma-invertibility-CR-horizontal-map}} guarantees
that this $n \times n$ matrix is invertible, we deduce:
\[
\left(\!
\begin{array}{ccc}
\mathcal{K}_{x_1}(f_1) & \cdots & \mathcal{K}_{x_1}(f_n)
\\
\vdots & \ddots & \vdots
\\
\mathcal{K}_{x_n}(f_1) & \cdots & \mathcal{K}_{x_n}(f_n)
\end{array}
\!\right)^{\!-1}
\,
\left(\!
\begin{array}{c}
\mathcal{K}_{x_1}(g)
\\
\vdots
\\
\mathcal{K}_{x_n}(g)
\end{array}
\!\right)
\,\,=\,\,
\left(\!
\begin{array}{c}
Q_{x_1'}'
\\
\vdots
\\
Q_{x_n'}'
\end{array}
\!\right).
\]
For instance, when $n = 1$, this writes simply as:
\[
\frac{g_x+Q_x\,g_y}{f_x+Q_x\,f_y}
\,=\,
Q_{x'}'.
\]

If we therefore define:
\[
R_{f,g}^1
\big(x,y,\,
y_{x_1},\dots,y_{x_n}
\big)
\,=:\,
\big(
x',y',\,
y_{x_1'}',\dots,y_{x_n'}'
\big),
\]
with:
\[
\aligned
x'
\,:=\,
f(x,y),
\\
y'
\,:=\,
g(x,y),
\endaligned
\]
and with:
\[
\left(\!
\begin{array}{c}
y_{x_1'}'
\\
\vdots
\\
y_{x_n'}'
\end{array}
\!\right)
\,\,:=\,\,
\left(\!
\begin{array}{ccc}
f_{1,x_1}+y_{x_1}f_y & \cdots & f_{n,x_1}+y_{x_1}f_{n,y}
\\
\vdots\,\, & \ddots & \vdots\,\,\,\,\,
\\
f_{1,x_n}+y_{x_n}f_y & \cdots & f_{n,x_n}+y_{x_n}f_{n,y}
\end{array}
\!\right)^{\!-1}\,
\left(\!
\begin{array}{c}
g_{x_1}+y_{x_1}g_y
\\
\vdots\,\,\,\,\,
\\
g_{x_n}+y_{x_n}g_y
\end{array}
\!\right),
\]
a check left shows that commutativity
of the $1$-jet diagram:
\[
\xymatrix{
\mathcal{M}
\ar[rr]^{(f,g,\varphi,\psi)}
\ar[d]_{J_x^1\mathcal{Q}_{a,b}}
&
&
\mathcal{M}'
\ar[d]^{J_{x'}^1\mathcal{Q}_{a',b'}'}
\\
\K^{n+1+n}
\ar[rr]_{R_{f,g}^1}
&
&
{\K'}^{n+1+n},
}
\]
holds true, by virtue of the identity~({\ref{jet-1-g-Q-Q-prime}}).

The construction of the maps $R_{f,g}^k$ for higher
jet orders $k$ proceeds by induction.
\endproof

\Section{\bf Invariancies of Jet Mappings Ranks}
\label{invariancies-jet-ranks}
\HEAD{{\ref{invariancies-jet-ranks}}.~{\sf Invariancies 
of Jet Mappings Ranks}
}{
Joël {\sc Merker} (Orsay)}

Now, suppose that we have a third submanifold of solutions
$\mathcal{M}''$ represented by:
\[
y''
\,=\,
Q''(x'',a'',b'')
\ \ \ \ \ \ \ \ \ \ \ \ \
\text{or}
\ \ \ \ \ \ \ \ \ \ \ \ \
b''
\,=\,
P''(a'',x'',y''),
\]
and that we have a sequence of $2$ composable equivalences:
\[
\mathcal{M}
\xrightarrow[{\rule[0pt]{40pt}{0pt}}]{}
\mathcal{M}'
\xrightarrow[{\rule[0pt]{40pt}{0pt}}]{}
\mathcal{M}'',
\]
and let us write the second map as:
\[
f'(x',y')
\,=:\,
x'',
\ \ \ \ \ \ \ \ \ \ \ \
g'(x',y')
\,=:\,
y'', 
\ \ \ \ \ \ \ \ \ \ \ \
\varphi'(a',b')
\,=:\,
a'', 
\ \ \ \ \ \ \ \ \ \ \ \
\psi'(a',b')
\,=:\,
b''.
\]
Denote the composition of the two maps as:
\[
\big(f'',g'',\varphi'',\psi''\big)
\,:=\,
\big(f',g',\varphi',\psi'\big)
\circ
\big(f,g,\varphi,\psi\big),
\]
that is to say:
\[
\aligned
f''(x,y)
&
\,=\,
f'
\big(
f(x,y),\,
g(x,y)
\big),
\ \ \ \ \ \ \ \ \ \ \ \ \ \ \ \ \ \ \ \
g''(x,y)
\,=\,
g'
\big(
f(x,y),\,
g(x,y)
\big),
\\
\varphi''(a,b)
&
\,=\,
\varphi'
\big(\varphi(a,b),\,\psi(a,b)\big),
\ \ \ \ \ \ \ \ \ \ \ \ \ \ \ \ \ \ \ \
\psi''(a,b)
\,=\,
\psi'
\big(\varphi(a,b),\,\psi(a,b)\big).
\endaligned
\]

The naturality of the above construction yields:

\begin{Lemma}
In the composition diagram:
\[
\xymatrix{
\mathcal{M}
\ar[rr]^{(f,g,\varphi,\psi)}
\ar[d]_{J_x^k\mathcal{Q}_{a,b}}
\ar@/^3pc/[rrrr]^{(f'',g'',\varphi'',\psi'')}
&
&
\mathcal{M}'
\ar[rr]^{(f',g',\varphi',\psi')}
\ar[d]^{J_{x'}^k\mathcal{Q}_{a',b'}'}
&
&
\mathcal{M}''
\ar[d]^{J_{x''}^k\mathcal{Q}_{a'',b''}''}
\\
\K^{n+\frac{(n+k)!}{n!\,k!}}
\ar[rr]_{R_{f,g}^k}
\ar@/_3pc/[rrrr]_{R_{f'',g''}^k}
&
&
{\K'}^{n+\frac{(n+k)!}{n!\,k!}}
\ar[rr]_{R_{f',g'}^k}
&
&
{\K''}^{n+\frac{(n+k)!}{n!\,k!}},
}
\]
one has:
\[
R_{f'',g''}^k
\,=\,
R_{f',g'}^k
\circ
R_{f,g}^k.
\eqno\qed
\]
\end{Lemma}

As an application, take the inverses:
\[
(f',g')
\,:=\,
(f,g)^{-1}
\ \ \ \ \ \ \ \ \ \ \ \ \
\text{and}
\ \ \ \ \ \ \ \ \ \ \ \ \
(\varphi',\psi')
\,:=\,
(\varphi,\psi)^{-1}.
\]
Since an elementary check convinces that:
\[
\Id
\,=\,
R_\Idsmall^k
\,=\,
R_{(f,g)^{-1}}^k
\circ
R_{f,g}^k,
\]
we deduce that $R_{f,g}^k$ is {\em also} an equivalence
(an invertible map):
\[
\K^{n+\frac{(n+k)!}{n!\,k!}}
\overset{\sim}{\xrightarrow[{\rule[0pt]{40pt}{0pt}}]{}}
\K^{n+\frac{(n+k)!}{n!\,k!}}.
\]

\begin{Corollary}
\label{invariances-ranks-jet-maps}
For every equivalence
$\mathcal{M} \overset{\sim}{\longrightarrow}
\mathcal{M}'$, and for every integers $k \geqslant 0$ and
$l \geqslant 0$, 
the ranks at all points of the two pairs of jet maps:
\[
\aligned
\xymatrix{
\mathcal{M}
\ar[rr]
\ar[d]_{J_x^k\mathcal{Q}_{a,b}}
&
&
\mathcal{M}'
\ar[d]^{J_{x'}^k\mathcal{Q}_{a',b'}'}
\\
\K^{n+\frac{(n+k)!}{n!\,k!}}
&
&
{\K'}^{n+\frac{(n+k)!}{n!\,k!}},
}
\ \ \ \ \ \ \ \ \ \ \ \ \
\ \ \ \ \ \ \ \ \ \ \ \ \
\xymatrix{
\mathcal{M}
\ar[rr]
\ar[d]_{J_a^l\mathcal{P}_{x,y}}
&
&
\mathcal{M}'
\ar[d]^{J_{a'}^l\mathcal{P}_{x',y'}'}
\\
\K^{m+\frac{(m+l)!}{m!\,l!}}
&
&
{\K'}^{m+\frac{(m+l)!}{m!\,l!}},
}
\endaligned
\]
are the same:
\[
\aligned
\rank\,
\big(
J_x^k\mathcal{Q}_{a,b}
\big)
\,=\,
\rank\,
\Big(
J_{f(x,y)}^k
Q_{\varphi(a,b),\psi(a,b)}'
\Big),
\\
\rank\,
\big(
J_a^l\mathcal{P}_{x,y}
\big)
\,=\,
\rank\,
\Big(
J_{\varphi(a,b)}^l
P_{f(x,y),g(x,y)}'
\Big),
\endaligned
\]
and their generic ranks:
\[
\genrank\,
\big(
J_\smallbullet^k\mathcal{Q}_{\smallbullet,\smallbullet}
\big)
\,=\,
\genrank\,
\big(
J_\smallbullet^k\mathcal{Q}_{\smallbullet,\smallbullet}'
\big)
\ \ \ \ \ \ \ \ \ \ \ \ \
\text{and}
\ \ \ \ \ \ \ \ \ \ \ \ \
\genrank\,
\big(
J_\smallbullet^l\mathcal{P}_{\smallbullet,\smallbullet}
\big)
\,=\,
\genrank\,
\big(
J_\smallbullet^l\mathcal{P}_{\smallbullet,\smallbullet}'
\big),
\]
are also identical.\qed
\end{Corollary}

These ranks are therefore {\sl invariants} of the para-CR
structures under study. 

\smallskip

Another basic fact is a stabilization property.

\begin{Proposition}
\label{Proposition-stabilization-jet-ranks}
If, for some integer $k \geqslant 1$:
\[
\genrank\,
\big(
J_x^k\mathcal{Q}_{a,b}
\big)
\,=\,
\genrank\,
\big(
J_x^{k+1}\mathcal{Q}_{a,b}
\big),
\]
then for all integers $l \geqslant 1$ as well:
\[
\genrank\,
\big(
J_x^k\mathcal{Q}_{a,b}
\big)
\,=\,
\genrank\,
\big(
J_x^{k+l}\mathcal{Q}_{a,b}
\big).
\eqno\qed
\]
\end{Proposition}

A particular, most degenerate case, is when the generic rank is 
smallest possible. 

\begin{Proposition}
If:
\[
\genrank\,
\big(
J_x^1\mathcal{Q}_{a,b}
\big)
\,=\,
n+1,
\]
then $\mathcal{M}$ is equivalent to the flat hyperplane:
\[
\big\{
y
=
b
+
0
\big\}.\eqno\qed
\]
\end{Proposition}

This statement uses the fact that already:
\[
\genrank\,
\big(
J_x^0\mathcal{Q}_{a,b}
\big)
\,=\,
n+1
\eqno
{\scriptstyle{(\text{\rm at every point})}}.
\]

\Section{\bf Finite Nondegeneracy with Respect to Parameters
and to Variables}
\label{two-finite-nondegeneracies}
\HEAD{{\ref{two-finite-nondegeneracies}}.~{\sf Finite 
Nondegeneracy with Respect to Parameters and to Variables}
}{
Joël {\sc Merker} (Orsay)}

Contrary to what happens for CR manifolds,
in the context of submanifolds of solutions,
there are two distinct and {\em non-equivalent} conditions of
finite nondegeneracy.

\begin{Definition}
A submanifold of solutions $\mathcal{M} \subset 
\K_\variablessmall^{n+1} \times \K_\parameterssmall^{m+1}$ 
is called {\sl finitely nondegenerate with respect to 
parameters} at a point $(x_0, a_0, b_0)$ if there exists
an integer $k \geqslant 0$ such that the invariant jet map
$J_x^k \mathcal{Q}_{a,b}$:
\[
(x,a,b)
\,\,\longmapsto\,\,
\Big(
x,\,\,
\big(
\partial_x^\beta Q(x,a,b)
\big)_{\vert\beta\vert\leqslant k}
\Big)
\]
is of maximal possible rank $n + 1 + m$ at $(x_0, a_0, b_0)$.
\end{Definition}

Equivalently, the map:
\[
(a,b)
\,\,\longmapsto\,\,
\big(
\partial_x^\beta Q(x_0,a,b)
\big)_{\vert\beta\vert\leqslant k}
\]
is of maximal possible rank $1 + m$ at $(a_0, b_0)$.

\begin{Terminology}
Say that $\mathcal{M}$ is {\sl $k$-nondegenerate with 
respect to parameters}
when $k < \infty$ is the smallest such integer.
\end{Terminology}

Thanks to the fundamental
Corollary~{\ref{invariances-ranks-jet-maps}}, this
first condition
is independent of coordinates, as is the next, second,
symmetric condition.

\begin{Definition}
A submanifold of solutions $\mathcal{M} \subset 
\K_\variablessmall^{n+1} \times \K_\parameterssmall^{m+1}$ 
is called {\sl finitely nondegenerate with respect to 
variables} at a point $(a_0, x_0, y_0)$ if there exists
an integer $l \geqslant 0$ such that the invariant jet map
$J_a^l \mathcal{P}_{x,y}$:
\[
(a,x,y)
\,\,\longmapsto\,\,
\Big(
a,\,\,
\big(
\partial_a^\gamma P(a,x,y)
\big)_{\vert\gamma\vert\leqslant l}
\Big)
\]
is of maximal possible rank $m + 1 + n$ at $(a_0, x_0, y_0)$.
\end{Definition}

Equivalently, the map:
\[
(x,y)
\,\,\longmapsto\,\,
\big(
\partial_a^\gamma P(a_0,x,y)
\big)_{\vert\gamma\vert\leqslant l}
\]
is of maximal possible rank $1 + n$ at $(x_0, y_0)$.

\begin{Terminology}
Say that $\mathcal{M}$ is {\sl $l$-nondegenerate with 
respect to variables}
when $l < \infty$ is the smallest such integer.
\end{Terminology}

\begin{Example}
For $n = m = 1$, the submanifold of solutions
$\{ y = b + xa\}$ is simultaneously finitely nondegenerate
with respect to parameters and to variables, but with
$n = m = 2$, the submanifold:
\[
\big\{
y
=
b+xa+xxb
\big\},
\]
is $2$-nondegenerate with respect to parameters,
and {\em not} $l$-nondegenerate with respect to variables
for any $l \in \N$.\qed
\end{Example}

Our goal in the subsequent paragraphs is to explore more what
happens when $m = n = 2$. But before studying specifically 
$5$-dimensional submanifolds $\mathcal{M} \subset
\K^{2+1} \times \K^{2+1}$, we need to yet expose some further
aspects of the general theory.

\Section{\bf Generic Ranks, Degeneracies, and Local Product Structures}
\label{generic-degeneracies}
\HEAD{{\ref{generic-degeneracies}}.~{\sf Generic Ranks, Degeneracies, 
and Local Product Structures}
}{
Joël {\sc Merker} (Orsay)}

Another view on the two invariant jet maps, in which the jet order
is now pushed to infinity:
\[
\aligned
J_x^\infty\mathcal{Q}_{a,b}
\colon\ \ \ \ \
\big(
x,a,b
\big)
&
\,\,\longmapsto\,\,
\Big(
x,\,\,
\big(
\partial_x^\beta Q(x,a,b)
\big)_{\beta\in\N^n}
\Big)
\,\,\in\,\,
\K^\infty,
\\
J_a^\infty\mathcal{P}_{x,y}
\colon\ \ \ \ \
\big(
a,x,y
\big)
&
\,\,\longmapsto\,\,
\Big(
a,\,\,
\big(
\partial_a^\gamma P(a,x,y)
\big)_{\gamma\in\N^m}
\Big)
\,\,\in\,\,
\K^\infty,
\endaligned
\]
consists in partly expanding the two equations of $\mathcal{M}$ as:
\[
\aligned
y
&
\,=\,
Q(x,a,b)
\,=\,
\sum_{\beta\in\N^n}\,
x^\beta\,Q_\beta(a,b),
\\
b
&
\,=\,
P(a,x,y)
\,=\,
\sum_{\gamma\in\N^m}\,
a^\gamma\,P_\gamma(x,y).
\endaligned
\]
Then by introducing some two infinite coefficients maps:
\[
\aligned
\mathcal{Q}_\parameterssmall^\infty
&
\,:=\,
(a,b)
\,\,\longmapsto\,\,
\big(
Q_\beta(a,b)
\big)_{\beta\in\N^n},
\\
\mathcal{P}_\variablessmall^\infty
&
\,:=\,
(x,y)
\,\,\longmapsto\,\,
\big(
P_\gamma(x,y)
\big)_{\gamma\in\N^n},
\endaligned
\]
a detailed analysis of the power series expansions: 
\[
\aligned
\partial_x^\beta Q(x,a,b)
&
\,=\,
\sum_{\beta_1\in\N^n}\,
x^{\beta_1}\,
\frac{(\beta+\beta_1)!}{\beta!\,\beta_1!}\,
Q_{\beta+\beta_1}(a,b),
\\
\partial_a^\gamma P(a,x,y)
&
\,=\,
\sum_{\gamma_1\in\N^m}\,
a^{\gamma_1}\,
\frac{(\gamma+\gamma_1)!}{\gamma!\,\gamma_1!}\,
P_{\gamma+\gamma_1}(x,y),
\endaligned
\]
conducts to a proof of

\begin{Proposition}
The generic ranks of these $4$ infinite maps are related by:
\reqnomode\usetagform{EngelLie}
\begin{align}
\genrank\,
\big(
J_\smallbullet^\infty\mathcal{Q}_{\smallbullet,\smallbullet}
\big)
\,\,=\,\,
n
+
\genrank\,
\big(
\mathcal{Q}_\parameterssmall^\infty
\big),
\notag
\\
\genrank\,
\big(
J_\smallbullet^\infty\mathcal{P}_{\smallbullet,\smallbullet}
\big)
\,\,=\,\,
m
+
\genrank\,
\big(
\mathcal{P}_\variablessmall^\infty
\big).
\tag{\qed}
\end{align}
\end{Proposition}

Recall that from the beginning, we assume
that the implicit defining function of
$\mathcal{M} = \{ R(x,y,a,b) = 0\}$ satisfies:
\[
R_y(0,0,0,0)
\,\neq\,
0
\,\neq\,
R_b(0,0,0,0),
\]
an assumption equivalent to:
\[
Q_b(0,0,0)
\,\neq\,
0
\,\neq\,
P_y(0,0,0),
\]
so that it comes that the ranks of the $0$-jet maps:
\[
\aligned
(x,a,b)
&
\,\,\longmapsto\,\,
Q(x,a,b),
\ \ \ \ \ \ \ \ \ \ \ \ \ \ \ \ \ \ \ \ \ \ \ \ \ \
(a,b)
\,\,\longmapsto\,\,
Q_0(a,b),
\\
(a,x,y)
&
\,\,\longmapsto\,\,
P(a,x,y),
\ \ \ \ \ \ \ \ \ \ \ \ \ \ \ \ \ \ \ \ \ \ \ \ \
(x,y)
\,\,\longmapsto\,\,
P_0(x,y),
\endaligned
\]
are all already equal to $1$. Consequently, the above
generic ranks are all $\geqslant 1$.

Let us therefore introduce a notation for them:
\[
\aligned
1
&
\,\,\leqslant\,\,
\genrank\,
\big(
\mathcal{Q}_\parameterssmall^\infty
\big)
\,\,=:\,\,
1
+
m_{\mathcal{M}}
\,\,\leqslant\,\,
1+m,
\\
1
&
\,\,\leqslant\,\,
\genrank\,
\big(
\mathcal{P}_\variablessmall^\infty
\big)
\,\,=:\,\,
1
+
n_{\mathcal{M}}\,
\,\,\leqslant\,\,
1+n.
\endaligned
\]
These two invariant integers $m_{\mathcal{M}}$ and $n_{\mathcal{M}}$
are intrinsice to $\mathcal{M}$, and they represent the
true-dimensional space $\K^{n_{\mathcal{M}}+1} \times
\K^{m_{\mathcal{M}}+1}$ in which $\mathcal{M}$ really lives (at least
generically), as expresses the next

\begin{Theorem}
\label{Theorem-two-degeneracies-jets-var-par}
{\bf (1)}\,
Locally in a neighborhood of a generic point of $\mathcal{M}$,
there exists a change of coordinates:
\[
(x,y,a,b)
\,\,\longmapsto\,\,
\big(
f(x,y),g(x,y),\varphi(a,b),\psi(a,b)
\big),
\]
which transforms $\mathcal{M}$ into a new $\mathcal{M}' 
\subset \K^{n+1} \times \K^{m+1}$ having equations:
\[
y'
\,=\,
Q'(x',a',b')
\ \ \ \ \ \ \ \ \ \ \ \ \
\Longleftrightarrow
\ \ \ \ \ \ \ \ \ \ \ \ \
b'
\,=\,
P'(a',x',y'),
\]
in which both $Q'$ and $P'$ are independent of the coordinates:
\[
\big(
x_{n_{\mathcal{M}}+1}',\dots,x_n'
\big)
\ \ \ \ \ \ \ \ \ \ \ \ \
\text{and}
\ \ \ \ \ \ \ \ \ \ \ \ \
\big(
a_{m_{\mathcal{M}}+1}',\dots,a_m'
\big).
\]

\smallskip\noindent{\bf (2)}\,
Furthermore, the reduced submanifold:
\[
\mathcal{M}_\reducedsmall'
\,\,\subset\,\,
\K^{n_{\mathcal{M}}}
\times
\K
\times
\K^{m_{\mathcal{M}}}
\times
\K
\]
having these equations is finitely nondegenerate 
at every point both with respect to parameters and to variables.

\smallskip\noindent{\bf (3)}\,
Lastly, $\mathcal{M}_\reducedsmall'$ is {\em not} locally 
equivalent to any product with either $\K_\variablessmall^1$
or $\K_\parameterssmall^1$.\qed
\end{Theorem}

A particular case is when $n_{\mathcal{M}} = m_{\mathcal{M}} = 0$, 
in which $\mathcal{M}$ is {\sl maximally flat:}
\[
\mathcal{M}
\,\,\cong\,\,
\big\{
y=b
\big\}.
\]

\Section{\bf Local Lie Group Structure}
\label{local-Lie-group-structure}
\HEAD{{\ref{local-Lie-group-structure}}.~{\sf Local Lie Group Structure}
}{
Joël {\sc Merker} (Orsay)}

These theoretical preliminaries justify to work only
with submanifolds of solutions that are finitely
nondegenerate both with respect to parameters and
to variables.

As is well known in the context of CR geometry, the group of local
automorphisms of an arbitrary CR submanifold of $\C^\NN$ may be {\em
not} finite-dimensional, and hence, it is necessary to impose
nondegeneracy conditions. A similar phenomenon occurs in our setting.
As in the preceding sections, two independent nondegeneracy conditions
are needed, instead of the single finite nondegeneracy condition for
CR manifolds.

Now, let us state
and present a {\em new general result}, 
not contained in~{\cite{Merker-2008}}.

\begin{Theorem}
\label{Thm-GM-para-CR}
Let $\mathcal{M} \subset \K^{n+1} \times \K^{m+1}$ be a hypersurface
which is $k$-nondegenerate with respect to parameters and
$l$-nondegenerate with respect to variables. Then its (pseudo) group
of self-transformations close to the identity:
\[
\Aut(\mathcal{M})
\,:=\,
\big\{
(F,\Phi)
\in
\Diff_\variablessmall\times\Diff_\parameterssmall
\colon\,
(F,\Phi)(\mathcal{M})
\subset
\mathcal{M}
\big\},
\]
is a local Lie group of dimension:
\[
\dim\,
\Aut
(\mathcal{M})
\,\,\leqslant\,\,
(n+1)\,
\binom{n+1+2k+2l}{n+1}
+
(m+1)\,
\binom{m+1+2k+2l}{m+1}.
\]
\end{Theorem}

The two nondegeneracy conditions are necessary
to apply
the Implicit Function Theorem at many steps in the proof of jet
parametrizations of $F$ and $\Phi$, 
{\em cf.}~\cite[Section 11]{Merker-2008}.
In fact, it follows from 
Theorem~{\ref{Theorem-two-degeneracies-jets-var-par}} above
that 
if $\mathcal{M}$ (connected) is either {\em not} finitely nondegenerate
with respect to variables, or
{\em not} finitely nondegenerate with respect to parameters
(at a Zariski-generic point), then 
$\Aut (\mathcal{M})$ is {\em infinite}-dimensional (at
such points).

Note that in general, a further condition on submanifolds, that is the
condition of {\sl minimality}~\cite{Gaussier-Merker-2004,
  Merker-2008}, is needed for finiteness of the dimension of $\Aut
(\mathcal{M})$. In our case of a hypersurface $\mathcal{M}$, however,
minimality follows from either one of the two finite nondegeneracy
conditions.

\proof
Now, by taking inspiration from~{\cite{Gaussier-Merker-2004}}, let us
present ideas of the proof.
Recall that a {\sl local Lie transformation group} consists of
a family of local analytic diffeomorphisms $x' = \tau(x; g)$, where
$x = (x_1, \dots, x_n)$ and $x' = (x_1', \dots, x_n')$ are
source and target space coordinates, and where 
$g = (g_1, \dots, g_r)$ are group parameters,
with $0' = \tau(0; 0)$, satisfying:

\smallskip\noindent{\bf (1)}\,
$\tau\big( \tau(x;g); g'\big) = \tau\big(x; \mu(g,g')\big)$,
for some local analytic (multiplication) map $\mu$ with $\mu(0,0) = 0$,
with $\mu(g,0) = g = \mu(0,g)$, which enjoys associativity:
\[
\mu\big(g,\,\mu(g',g'')\big)
\,=\,
\mu\big(\mu(g,g'),\,g''\big);
\]

\smallskip\noindent{\bf (2)}\,
$\tau(x; 0) = x$ and $\tau\big( \tau(x;g); \iota(g) \big) = x$,
for some local analytic (inversion) map $g \longmapsto \iota(g)$
with $\iota(0) = 0$ satisfying:
\[
\mu\big(g,\iota(g)\big)
\,=\,
g
\,=\,
\mu\big(\iota(g),g).
\]

These properties are supposed to hold near $0$, up to some
neighborhood shinking. The identity transformation $x \longmapsto x =
x'$ corresponds to $g = 0$, and considerations are localized near $x =
0$ and $x' = 0'$.

Next, abbreviate:
\[
G_{{\sf v},{\sf p}}
\,:=\,
\Diff_\variablessmall\times\Diff_\parameterssmall,
\]
and, for small $\varepsilon > 0$, with 
the nondegeneracy invariants $k$, $l$ of
Theorem~{\ref{Thm-GM-para-CR}}, introduce
a space of self-maps of $\mathcal{M}$
which are close to the identity map $\Id$ in the jet sense:
\[
\aligned
\Aut_\varepsilon(\mathcal{M})
\,:=\,
\Big\{
(F,\Phi)
\in
G_{{\sf v},{\sf p}}
\colon\,\,
&
(F,\Phi)(\mathcal{M})
\subset
\mathcal{M},
\\
&
\bignorm
j_{x,y}^{2k+2l}F
-
j_{x,y}^{2k+2l}\Id
\bignorm
<
\varepsilon,
\ \ \
\bignorm
j_{a,b}^{2k+2l}\Phi
-
j_{a,b}^{2k+2l}\Id
\bignorm
<
\varepsilon
\Big\},
\endaligned
\]
where order $\kappa$ jets are:
\[
j_{x,y}^\kappa F
\,:=\,
\big(
\partial_x^\beta\partial_y^iF
\big)_{\vert\beta\vert+i\leqslant\kappa},
\ \ \ \ \ \ \ \ \ \ \ \ \ \ \ \ \ \ \ \ \ \ \ \ \ \
j_{a,b}^\kappa\Phi
\,:=\,
\big(
\partial_a^\gamma\partial_b^j\Phi
\big)_{\vert\gamma\vert+j\leqslant\kappa}.
\]

In our particular case of a hypersurface $M$ in $\mathbb{C}^{n+1} \times
\mathbb{C}^{n+1}$, the type $(\mu,\mu^*)$ of $M$ as defined in
\cite[Definition 10.21]{Merker-2008} is always $(2,2)$ (see
\cite[Section 10]{Merker-2008}), and the jet parametrization reads as
follows.

\begin{Theorem} 
\cite[Theorem 11.6]{Merker-2008} 
\label{Thm-jet-parametrization}
Suppose that the hypersurface $\mathcal{M}$ 
is $k$-nondegenerate with respect to
parameters and $l$-nondegenerate with respect to variables. Then there
exist two local analytic maps:
\[
\aligned
H\colon\ \ \ \ \ 
\mathbb{K}^{n+1+(n+1)\binom{n+1+2k+2l}{2k+2}} 
&
\,\,\,\longrightarrow\,\,\,
\mathbb{K}^{n+1},
\\
\Pi\colon\ \ \ 
\mathbb{K}^{m+1+(m+1)\binom{m+1+2k+2l}{2k+2l}} 
&
\,\,\,\longrightarrow\,\,\,
\mathbb{K}^{n+1},
\endaligned
\]
which are constructed only from the defining
functions $Q,P$ of $\mathcal{M}$ 
so that the maps $F,\Phi$ satisfy:
\[
\aligned
F(x,y) 
&\,\equiv\, 
H
\big(x,y,
j^{2k+2l}_{x,y}F(0)
\big), 
\\
\Phi(a,b) 
&\,\equiv\, 
\Pi
\big(a,b,
j^{2k+2l}_{a,b}\Phi(0)
\big),
\endaligned
\]
as power series. 
\end{Theorem}

This theorem is in fact the para-CR analog of the jet
parametrization theorem for CR manifolds 
in~\cite[Theorem 6.4]{Gaussier-Merker-2004}. 
It is also clear
from the jet parametrizations of $F$ and $\Phi$ that the jet order
$\kappa$ in the definition of ${\sf Aut}_{\epsilon}
(\mathcal{M})$ is $\kappa =
2k+2l$.

Since Theorem~{\ref{Thm-jet-parametrization}}
concerns {\em any} self-map of $\mathcal{M}$ sufficiently
close to the identity, it applies to compositions
of such maps (provided both are sufficiently close
to the identity), and then this theorem shows
that the composition is again of the same form,
with the same parametrizing functions $H$, $\Pi$.
Also, the inverse of any self-map of $\mathcal{M}$
is again of this parametrized form. 

So Theorem~{\ref{Thm-jet-parametrization}}, provides exactly Lie's
original definition of a local (Lie) transformation
group~{\cite[pp.~3--4]{Engel-Lie-1888}},
{\cite[Section~3.1]{Lie-Merker-2015}}, hence Lie's general theory
applies. 
This would be enough to conclude the proof, 
but it yet remains to bound $\dim\, 
\Aut\, \mathcal{M}$, and this will show 
other aspects.

Taking either $R := Y - Q(x,a,b)$ or $R := b - P(a,x,y)$, or even
any other definining equation for the hypersurface $\mathcal{M}$:
\[
\mathcal{M}
\colon\ \ \ \ \
0
\,=\,
R\big(x,y,\,a,b\big),
\]
the assumption that $(F, \Phi)$ is a local self-map of 
$\mathcal{M}$ (close to the identity, up to shrinking neighborhoods),
reads as:
\[
(x,y,\,a,b)
\,\in\,
\mathcal{M}
\ \ \ \ \ \ \ \
\Longrightarrow
\ \ \ \ \ \ \ \ 
\big(
F(x,y),\,
\Phi(a,b)
\big)
\,\in\,
\mathcal{M},
\]
that is, it expresses as the identical vanishing of the following
convergent power series in the variables $(x, a, b)$:
\[
0
\,\equiv\,
R\Big(
F\big(x,\,Q(x,a,b)\big),\,\,
\Phi(a,b)
\Big).
\]

The key act is then to plug the two jet parametrization formulas
of Theorem~{\ref{Thm-jet-parametrization}} in this identity:
\[
0
\,\equiv\,
R
\bigg(
H\Big(
x,\,Q\big(x,a,b)\big),\,
j_{x,y}^{2k+2l}F(0)
\Big),\,\,\,
\Pi\Big(
a,b,\,
j_{a,b}^{2k+2l}\Phi(0)
\Big)
\bigg),
\]
and then to power-expand:
\[
0
\,\equiv\,
\sum_{\beta\in\N^n}\,
\sum_{\gamma\in\N^m}\,
\sum_{j\in\N}\,
x^\beta\,a^\gamma\,
b^j\,
C_{\beta,\gamma,j}
\Big(
j_{x,y}^{2k+2l}F(0),\,
j_{a,b}^{2k+2l}\Phi(0)
\Big),
\]
getting an infinite collection of equations in the finite-order
jets of $F$ and $\Phi$:
\[
0
\,=\,
C_{\beta,\gamma,j}
\Big(
j_{x,y}^{2k+2l}F(0),\,
j_{a,b}^{2k+2l}\Phi(0)
\Big)
\eqno
{\scriptstyle{(\forall\,\beta,\,\forall\,\gamma,\,\forall\,j)}},
\]
which define a certain analytic subvariety of $\K^{\NN_{n,m,k,l}}$,
with the following integer counting all such jets:
\[
\NN
\,:=\,
\NN_{n,m,k,l}
\,:=\,
\big(
n+1
\big)\,
\binom{n+1+2k+2l}{n+1}
+
\big(
m+1
\big)\,
\binom{m+1+2k+2l}{m+1}.
\]
By Noetherianity of the local ring of $\K$-analytic functions,
only a finite number of equations suffices. 

By arguing as in~{\cite[Lemma~6.5]{Gaussier-Merker-2004}},
one can show that this analytic
subvariety is in fact geometrically smooth
at the origin, as a consequence 
that compositions and inversions
of elements of $\Aut_\varepsilon (\mathcal{M})$
are still of the form 
given by Theorem~{\ref{Thm-jet-parametrization}}.
But this was already done by Engel and Lie 
in~{\cite{Engel-Lie-1888}} in their general theory, what the authors 
of~{\cite{Gaussier-Merker-2004}} were not aware of.

Thus, we certainly have the announced dimension bound:
\[
\dim\,\mathcal{M}
\,\leqslant\,
\big(
n+1
\big)\,
\binom{n+1+2k+2l}{n+1}
+
\big(
m+1
\big)\,
\binom{m+1+2k+2l}{m+1}.
\qedhere
\]
\endproof

However, let us close this section by mentioning that analyzing these
equations $C_{\beta, \gamma, j} = 0$ is the `{\sl hard core}' of the
equivalence\big/classification problem, which is 
very complicated and very
ramified even in small dimensions.
By mastering Cartan's method of equivalence, 
Hill and Nurowski~{\cite{Hill-Nurowski-2010}}
were able to understand and to analyze the `{\sl hard core}'
of quite a number of para-CR structures of specific
dimensions and codimensions.

Anyway, in the para-CR context, this Theorem~{\ref{Thm-GM-para-CR}}
seems to appear nowhere in the literature.

Also, Theorem~{\ref{Thm-GM-para-CR}} generalizes easily to para-CR
structures of any codimension, with the supplementary assumption of
{\em minimality}, or of {\sl covering property}, {\em
see}~{\cite[Section~9]{Merker-2008}}. The proof is the same as in
the hypersurface case (to which this paper restricts itself), with
the same proof, since the jet parametrization
Theorem~{\ref{Thm-jet-parametrization}} was
already stated and proved in any
codimension as Theorem~11.6 in~{\cite{Merker-2008}}.

\Section{\bf Generic $(2,2,1)$ Para-CR Structures}
\label{generic-2-2-1-para-CR-structures}
\HEAD{{\ref{generic-2-2-1-para-CR-structures}}.~{\sf Generic 
$(2,2,1)$ Para-CR Structures}
}{
Joël {\sc Merker} (Orsay)}

Let us therefore examine the special case $n = m = 2$,
still in codimension $c = 1$. The rest of the paper will
in fact not touch anymore the general theory. Hence,
we allow ourselves to change notation from now 
on, and to write out the two
graphed equations of $\mathcal{M}$ as:
\[
z
\,=\,
Q(x,y,a,b,c)
\ \ \ \ \ \ \ \ \ \ \ \ \
\Longleftrightarrow
\ \ \ \ \ \ \ \ \ \ \ \ \
c
\,=\,
P(a,b,x,y,z).
\]

We already saw that when $\Levi_\variablessmall \equiv 0 \equiv
\Levi_\parameterssmall$, there is an equivalence
to the (very) flat model $\{y = b\}$.

At the other extreme lies the case where the (identical) rank(s) of
the (two) Levi form(s) is (are) equal to $2$. This case 
has already been treated by Hachtroudi~{\cite{Hachtroudi-1937}} 
in 1937.

What about intermediate cases?

Since we are working at a generic point, we have therefore
to examine the remaining branch:
\[
\genrank\,
\big(
\Levi_\variablessmall(Q)
\big)
\,=\,
1
\,=\,
\genrank\,
\big(
\Levi_\parameterssmall(P)
\big).
\]
But remind that in the previous paragraphs, the two Levi forms have
been inserted in the wider, more adequate and more general concept of
jets of leaves of the two invariant foliations on $\mathcal{M}$. More
precisely, recall that according to 
Corollary~{\ref{Corollary-jets-1-identical-Levi-form}}:
\[
\rank\,
\big(
J_x^1\mathcal{Q}_{a,b}
\big)
\,=\,
3
+
\rank\,
\big(
\Levi_\variablessmall(Q)
\big)
\,=\,
3
+
\rank\,
\big(
\Levi_\parameterssmall(P)
\big)
\,=\,
\rank\,
\big(
J_a^1\mathcal{P}_{x,y}
\big),
\]
whence we are looking at the branch:
\[
\genrank\,
\big(
J_x^1\mathcal{Q}_{a,b}
\big)
\,=\,
4
\,=\,
\genrank\,
\big(
J_a^1\mathcal{P}_{x,y}
\big).
\]

Now, remind that 
Proposition~{\ref{Proposition-stabilization-jet-ranks}}
has shown that second order (generic) jet ranks 
can either increase, or stabilize, 
hence it may happen that they jump to
the maximal value $5$, or stay at $4$:
\[
\genrank\,
\big(
J_x^2\mathcal{Q}_{a,b}
\big)
\,=\,
\left\{
\aligned
{}
&
4,
\\
{}
&
5,
\endaligned\right.
\ \ \ \ \ \ \ \ \ \ \ \ \
\text{and}
\ \ \ \ \ \ \ \ \ \ \ \ \
\genrank\,
\big(
J_a^2\mathcal{P}_{x,y}
\big)
\,=\,
\left\{
\aligned
{}
&
4,
\\
{}
&
5.
\endaligned\right.
\] 
Consequently, it is obvious that 
we have to examine three degenerate cases:
\[
\aligned
\text{\bf Case~I}
\ \ \ \ \
&
\left[
\aligned
{\bf 4}
&
\,=\,
\genrank\,
\big(
J_\smallbullet^{\bf 1}\mathcal{Q}_{\smallbullet,\smallbullet}
\big)
\ \ \ \ \ \ \ \ \ \ \ \ \
\text{and}
\ \ \ \ \ \ \ \ \ \ \ \ \ \ \ \ 
{\bf 4}
\,=\,
\genrank\,
\big(
J_\smallbullet^{\bf 2}\mathcal{Q}_{\smallbullet,\smallbullet}
\big),
\\
{\bf 4}
&
\,=\,
\genrank\,
\big(
J_\smallbullet^{\bf 1}\mathcal{P}_{\smallbullet,\smallbullet}
\big)
\ \ \ \ \ \ \ \ \ \ \ \ \,
\text{and}
\ \ \ \ \ \ \ \ \ \ \ \ \ \ \ \
{\bf 4}
\,=\,
\genrank\,
\big(
J_\smallbullet^{\bf 2}\mathcal{P}_{\smallbullet,\smallbullet}
\big),
\endaligned\right.
\\
\text{\bf Case~II}
\ \ \ \ \
&
\left[
\aligned
{\bf 4}
&
\,=\,
\genrank\,
\big(
J_\smallbullet^{\bf 1}\mathcal{Q}_{\smallbullet,\smallbullet}
\big)
\ \ \ \ \ \ \ \ \ \ \ \ \
\text{and}
\ \ \ \ \ \ \ \ \ \ \ \ \ \ \ \ \,
{\bf 4}
\,=\,
\genrank\,
\big(
J_\smallbullet^{\bf 2}\mathcal{Q}_{\smallbullet,\smallbullet}
\big),
\\
{\bf 4}
&
\,=\,
\genrank\,
\big(
J_\smallbullet^{\bf 1}\mathcal{P}_{\smallbullet,\smallbullet}
\big)
\ \ \ \ \ \ \ \ \ \ \ \ \
\text{while}
\ \ \ \ \ \ \ \ \ \ \ \ \
{\bf 5}
\,=\,
\genrank\,
\big(
J_\smallbullet^{\bf 2}\mathcal{P}_{\smallbullet,\smallbullet}
\big),
\endaligned\right.
\\
\text{\bf Case~III}
\ \ \ \ \
&
\left[
\aligned
{\bf 4}
&
\,=\,
\genrank\,
\big(
J_\smallbullet^{\bf 1}\mathcal{Q}_{\smallbullet,\smallbullet}
\big)
\ \ \ \ \ \ \ \ \ \ \ \ \ \,
\text{while}
\ \ \ \ \ \ \ \ \ \ \ \ \
{\bf 5}
\,=\,
\genrank\,
\big(
J_\smallbullet^{\bf 2}\mathcal{Q}_{\smallbullet,\smallbullet}
\big),
\\
{\bf 4}
&
\,=\,
\genrank\,
\big(
J_\smallbullet^{\bf 1}\mathcal{P}_{\smallbullet,\smallbullet}
\big)
\ \ \ \ \ \ \ \ \ \ \ \ \
\text{and}
\ \ \ \ \ \ \ \ \ \ \ \ \ \ \ \
{\bf 4}
\,=\,
\genrank\,
\big(
J_\smallbullet^{\bf 2}\mathcal{P}_{\smallbullet,\smallbullet}
\big),
\endaligned\right.
\endaligned
\]
plus one nondegenerate case:
\[
\text{\bf Case~IV}
\ \ \ \ \
\left[
\aligned
{\bf 4}
&
\,=\,
\genrank\,
\big(
J_\smallbullet^{\bf 1}\mathcal{Q}_{\smallbullet,\smallbullet}
\big)
\ \ \ \ \ \ \ \ \ \ \ \ \ \,
\text{while}
\ \ \ \ \ \ \ \ \ \ \ \ \
{\bf 5}
\,=\,
\genrank\,
\big(
J_\smallbullet^{\bf 2}\mathcal{Q}_{\smallbullet,\smallbullet}
\big),
\\
{\bf 4}
&
\,=\,
\genrank\,
\big(
J_\smallbullet^{\bf 1}\mathcal{P}_{\smallbullet,\smallbullet}
\big)
\ \ \ \ \ \ \ \ \ \ \ \ \
\text{while}
\ \ \ \ \ \ \ \ \ \ \ \ \ 
{\bf 5}
\,=\,
\genrank\,
\big(
J_\smallbullet^{\bf 2}\mathcal{P}_{\smallbullet,\smallbullet}
\big).
\endaligned\right.
\]

Serendipitously, Theorem~{\ref{Theorem-two-degeneracies-jets-var-par}}
offers without work 

\begin{Proposition}
At a generic point:
\reqnomode\usetagform{EngelLie}
\begin{align}
\text{\bf Case~I}
\ \ \ \ \
\mathcal{M}
&
\,\,\cong\,\,
\big\{
z
=
Q(x,\emptyset,a,\emptyset,c)
\big\}
\ \ \ \ \
\text{\rm or}
\ \ \ \ \
\big\{
c
\,=\,
P(a,\emptyset,x,\emptyset,z)
\big\},
\notag
\\
\text{\bf Case~II}
\ \ \ \ \
\mathcal{M}
&
\,\,\cong\,\,
\big\{
z
=
Q(x,y,a,\emptyset,c)
\big\}
\ \ \ \ \
\text{\rm or}
\ \ \ \ \
\big\{
c
\,=\,
P(a,\emptyset,x,y,z)
\big\},
\notag
\\
\text{\bf Case~III}
\ \ \ \ \
\mathcal{M}
&
\,\,\cong\,\,
\big\{
z
=
Q(x,\emptyset,a,b,c)
\big\}
\ \ \ \ \
\text{\rm or}
\ \ \ \ \
\big\{
c
\,=\,
P(a,b,x,\emptyset,z)
\big\}.
\tag{\qed}
\end{align}
\end{Proposition}

These three degenerate cases could be treated separately,
but we prefer to study only Case IV.
In this case, also according
to Theorem~{\ref{Theorem-two-degeneracies-jets-var-par}},
we can assume that $\mathcal{M}$ is both
$2$-nondegenerate with respect to parameters 
and $2$-nondegenerate with respect to variables.

\Section{\bf Vanishing of the Two $3 \times 3$ Levi Determinants}
\label{vanishing-two-3-3-Levi-determinants}
\HEAD{{\ref{vanishing-two-3-3-Levi-determinants}}.~{\sf Vanishing 
of the Two $3 \times 3$ Levi Determinants}
}{
Joël {\sc Merker} (Orsay)}

Thus, recall that we assume that the two Levi
forms are of rank $1$ at every point, which expresses as:
\[
\aligned
0
&
\,\equiv\,
\det\,
\Levi_\parameterssmall(Q)
\,=\,
\left\vert\!
\begin{array}{ccc}
\frac{-Q_cQ_{xa}+Q_aQ_{xc}}{Q_c\,Q_c}
&
\frac{-Q_cQ_{xb}+Q_bQ_{xc}}{Q_c\,Q_c}
\\
\frac{-Q_cQ_{ya}+Q_aQ_{yc}}{Q_c\,Q_c}
&
\frac{-Q_cQ_{yb}+Q_bQ_{yc}}{Q_c\,Q_c}
\end{array}
\!\right\vert,
\\
0
&
\,\equiv\,
\det\,
\Levi_\variablessmall(P)
\,=\,
\left\vert\!
\begin{array}{ccc}
\frac{-P_zP_{ax}+P_xP_{az}}{P_z\,P_z}
&
\frac{-P_zP_{ay}+P_xP_{ay}}{P_z\,P_z}
\\
\frac{-P_zP_{bx}+P_xP_{bz}}{P_z\,P_z}
&
\frac{-P_zP_{by}+P_yP_{bz}}{P_z\,P_z}
\end{array}
\!\right\vert.
\endaligned
\]
Moreover, the general relation (known from the end of
Section~{\ref{nonholonomy-form-and-dual}}):
\[
\Levi_\parameterssmall(Q)
\,=\,
-\,P_z\,
{}^\TT
\Levi_\variablessmall(P),
\]
gives after taking determinants, using $Q_c\, P_z = 1$:
\[
\det\,
\Levi_\parameterssmall(Q)
\,=\,
P_zP_z\,
\det\,
\Levi_\variablessmall(P)
\ \ \ \ \ \ \ \ \ \ \ \ \
\Longleftrightarrow
\ \ \ \ \ \ \ \ \ \ \ \ \
Q_cQ_c\,
\det\,
\Levi_\parameterssmall(Q)
\,=\,
\det\,
\Levi_\variablessmall(P).
\]

By what precedes, the two degeneracies of the two Levi forms
are linked with the degeneracies of the two 
ranks of the two (reduced) first jet maps:
\[
(a,b,c)
\,\,\longmapsto\,\,
\big(
Q,\,Q_x,\,Q_y
\big)
\ \ \ \ \ \ \ \ \ \ \ \ \
\text{and}
\ \ \ \ \ \ \ \ \ \ \ \ \
(x,y,z)
\,\,\longmapsto\,\,
\big(
P,\,P_a,\,P_b
\big),
\]
whose respective Jacobian determinants are:
\[
{\sf L}_{3\times 3}(Q)
\,:=\,
\left\vert\!
\begin{array}{ccc}
Q_a & Q_b & Q_c
\\
Q_{xa} & Q_{xb} & Q_{xc}
\\
Q_{ya} & Q_{yb} & Q_{yc}
\end{array}
\!\right\vert
\ \ \ \ \ \ \ \ \ \ \ \ \
\text{\rm and}
\ \ \ \ \ \ \ \ \ \ \ \ \
\left\vert\!
\begin{array}{ccc}
P_x & P_y & P_z
\\
P_{xa} & P_{ya} & P_{za}
\\
P_{xb} & P_{yb} & P_{zb}
\end{array}
\!\right\vert
\,=:\,
{\sf L}_{3\times 3}(P),
\]
and indeed here, a direct check 
(just develop determinants) convinces that:
\[
\det\,\Levi_\parameterssmall(Q)
\,=\,
\frac{1}{Q_cQ_cQ_c}\,
{\sf L}_{3\times 3}(Q)
\ \ \ \ \ \ \ \ \ \ \ \ \
\text{and}
\ \ \ \ \ \ \ \ \ \ \ \ \
\det\,\Levi_\variablessmall(P)
\,=\,
\frac{1}{P_zP_zP_z}\,
{\sf L}_{3\times 3}(P).
\]
As a byproduct, it comes:
\[
{\sf L}_{3\times 3}(Q)
\,=\,
Q_cQ_cQ_c\,
\det\,\Levi_\parameterssmall(Q)
\,=\,
Q_c\,
\det\,\Levi_\variablessmall(P)
\,=\,
\frac{1}{P_zP_zP_zP_z}\,
{\sf L}_{3\times 3}(P),
\]
which shows (again) that the vanishings of these two $3 \times 3$
determinants are equivalent.
In summary, we are assuming in Case IV that:
\[
0
\,\equiv\,
\left\vert\!
\begin{array}{ccc}
Q_a & Q_b & Q_c
\\
Q_{xa} & Q_{xb} & Q_{xc}
\\
Q_{ya} & Q_{yb} & Q_{yc}
\end{array}
\!\right\vert
\ \ \ \ \ \ \ \ \ \ \ \ \
\Longleftrightarrow
\ \ \ \ \ \ \ \ \ \ \ \ \
0
\,\equiv\,
\left\vert\!
\begin{array}{ccc}
P_x & P_y & P_z
\\
P_{xa} & P_{ya} & P_{za}
\\
P_{xb} & P_{yb} & P_{zb}
\end{array}
\!\right\vert.
\]

\Section{\bf Local Graphs for Doubly $2$-Nondegenerate Submanifolds}
\label{local-graphs-twin-2-nondegenerate-M}
\HEAD{{\ref{local-graphs-twin-2-nondegenerate-M}}.~{\sf Local 
Graphs for Doubly $2$-Nondegenerate Submanifolds}
}{
Joël {\sc Merker} (Orsay)}

Assume therefore (Case IV)
that $\mathcal{M}$, represented in coordinates
$(x,y,z,a,b,c)$ by two graphed equations:
\[
z
\,=\,
Q(x,y,a,b,c)
\ \ \ \ \ \ \ \ \ \ \ \ \
\text{and}
\ \ \ \ \ \ \ \ \ \ \ \ \
c
\,=\,
P(a,b,x,y,z),
\]
is simultaneously $2$-nondegenerate with respect to parameters
and with respect to variables. If either $(x,y,a,b,c)$
or $(a,b,x,y,z)$ are taken as (horizontal) 
coordinates on $\mathcal{M}$,
we have two choices of fundamental vector fields:
\[
\left\{
\aligned
\mathcal{K}_x
&
\,:=\,
\frac{\partial}{\partial x},
\\
\mathcal{K}_y
&
\,:=\,
\frac{\partial}{\partial y},
\\
\mathcal{L}_a
&
\,:=\,
\frac{\partial}{\partial a}
-
\frac{Q_a}{Q_c}\,
\frac{\partial}{\partial c},
\\
\mathcal{L}_b
&
\,:=\,
\frac{\partial}{\partial b}
-
\frac{Q_b}{Q_c}\,
\frac{\partial}{\partial c},
\endaligned\right.
\ \ \ \ \ \ \ \ \ \ \ \ \
\text{or}
\ \ \ \ \ \ \ \ \ \ \ \ \
\left\{
\aligned
\mathcal{K}_x
&
\,:=\,
\frac{\partial}{\partial x}
-
\frac{P_x}{P_z}\,
\frac{\partial}{\partial z},
\\
\mathcal{K}_y
&
\,:=\,
\frac{\partial}{\partial y}
-
\frac{P_y}{P_z}\,
\frac{\partial}{\partial z},
\\
\mathcal{L}_a
&
\,:=\,
\frac{\partial}{\partial a},
\\
\mathcal{L}_b
&
\,:=\,
\frac{\partial}{\partial b}.
\endaligned\right.
\]

Corollary~{\ref{normalization-Levi-form-origin}}
has already shown that there exist coordinates
with:
\[
c
\,\equiv\,
Q(0,0,a,b,c)
\,\equiv\,
Q(x,y,0,0,c)
\ \ \ \ \ \ \ \ \ \ \ \ \
\text{and}
\ \ \ \ \ \ \ \ \ \ \ \ \
z
\,\equiv\,
P(0,0,x,y,z)
\,\equiv\,
P(a,b,0,0,z),
\]
in which the two graphing functions $P$ and $Q$ have
normalized second-order terms:
\[
\aligned
z
&
\,=\,
c
+
xa
+
{\rm O}_{x,y,a,b}(3)
+
c\,{\rm O}_{x,y,a,b,c}(2),
\\
c
&
\,=\,
z
-
ax
+
{\rm O}_{a,b,x,y}(3)
+
z\,{\rm O}_{a,b,x,y,z}(2).
\endaligned
\]

Next, specifying the homogeneous order $3$ terms, we can write:
\[
z
\,=\,
c
+
xa
+
Q_3(x,y,a,b)
+
c\,Q_2^\sim(x,y,a,b)
+
{\rm O}_{x,y,a,b,c}(4),
\]
with $Q_3$ and $Q_2^\sim$ homogeneous of degrees $3$ and $2$,
independent of $c$. In the $2 \times 2$ matrix
$\Levi_\parameterssmall (Q)$, it is easy to see that
the $(1,2)$ and the $(2,1)$ entries are ${\rm O}_1$,
while the $(2,2)$-entry is:
\[
\aligned
\frac{-Q_cQ_{yb}+Q_bQ_{yc}}{Q_c\,Q_c}
&
\,=\,
\frac{(-1+{\rm O}_1)\,\big(Q_{3,yb}+c\,Q_{2,yb}^\sim)+{\rm O}_2}{
(-1+{\rm O}_1)^2}
\\
&
\,=\,
\underbrace{
-\,Q_{3,yb}
-
c\,Q_{2,yb}^\sim}_{=\,\,{\rm O}_1}
+
{\rm O}_2,
\endaligned
\]
and therefore, the identical vanishing:
\[
0
\,\equiv\,
\left\vert\!
\begin{array}{cc}
-1+{\rm O}_1 & {\rm O}_1
\\
{\rm O}_1 & -Q_{3,yb}-c\,Q_{2,yb}^\sim+{\rm O}_2
\end{array}
\!\right\vert,
\]
yields an interesting annihilation [the second one will not be used]:
\[
0
\,\equiv\,
Q_{3,yb}
\ \ \ \ \ \ \ \ \ \ \ \ \
\Big[\,
\text{and}
\ \ \ \ \ \ \ \ \ \ \ \ \
0
\,\equiv\,
Q_{2,yb}^\sim
\Big].
\]

Then, expanding the cubic homogeneous polynomial $Q_3$
without pure $(x,y)$-terms or $(a,b)$-terms, all of which
monomials divisible by $yb$ should be absent, 
it remains:
\[
\aligned
z
\,=\,
c
+
xa
&
+
\alpha\,xxa
+
\underline{\alpha}\,xaa
+
\beta\,xxb
+
\underline{\beta}\,yaa
\,+
\\
&
+
\gamma\,xya
+
\underline{\gamma}\,xab
+
\delta\,yya
+
\underline{\delta}\,xbb
\,+
\\
&
+
c\,{\rm O}_{x,y,a,b}(2)
+
{\rm O}_{x,y,a,b,c}(4).
\endaligned
\]
But by replacing:
\[
\aligned
x
+
\alpha\,xx
+
\gamma\,xy
+
\delta\,yy
&
\,=:\,
x',
\\
a
+
\underline{\alpha}\,aa
+
\underline{\gamma}\,ab
+
\underline{\delta}\,bb
&
\,=:\,
a',
\endaligned
\]
and dropping primes, we come to:
\[
z
\,=\,
c
+
xa
+
\beta\,xxb
+
\underline{\beta}\,yaa
+
c\,{\rm O}_{x,y,a,b}(2)
+
{\rm O}_{x,y,a,b,c}(4).
\]

\begin{Assertion}
$\mathcal{M}$ is $2$-nondegenerate with respect to parameters
if and only if:
\[
\beta
\,\neq\,
0.
\]
\end{Assertion}

\proof
We have to guarantee that the rank at $0$ of the second jet map:
\[
(a,b,c)
\,\,\longmapsto\,\,
\Big(
Q,\,\,
Q_x,\,\,
Q_y,\,\,
Q_{xx},\,\,
Q_{xy},\,\,
Q_{yy}
\Big),
\]
is maximal equal to $3$. But after all these normalizations,
it is the map:
\[
(a,b,c)
\,\,\longmapsto\,\,
\Big(
c+{\rm O}_2,\,\,
a+{\rm O}_2,\,\,
\zero{{\rm O}_2},\,\,
2\beta\,b+{\rm O}_2,\,\,
\zero{{\rm O}_2},\,\,
\zero{{\rm O}_2}
\Big)
\]
whose Jacobian determinant at the origin (after dropping
its useless components $3$, $5$, $6$) is:
\[
\det\,
\left(\!
\begin{array}{ccc}
0 & 0 & 1
\\
1 & 0 & 0
\\
0 & 2\,\beta & 0
\end{array}
\!\right)
\,=\,
2\,\beta.
\qedhere
\]
\endproof

It is important to realize (repeat) that {\em another}, 
{\em different}, condition of $2$-nondegeneracy 
exists in the context of submanifolds of solutions
of {\sc pde} systems, contrary to the CR case in which
both $2$-nondegeneracy conditions are complex-conjugate
one to another, hence reduce to a {\em single} condition.

\begin{Assertion}
$\mathcal{M}$ is $2$-nondegenerate with respect to variables
if and only if:
\[
\underline{\beta}
\,\neq\,
0.
\]
\end{Assertion}

\proof
As we know, the implicit function theorem enables to solve $c$
as $c = P (a,b,x,y,z)$, and the result writes out as:
\[
c
\,=\,
z
-
ax
-
bxx
-
\beta\,xxb
-
\underline{\beta}\,aay
+
z\,{\rm O}_{a,b,x,y}(2)
+
{\rm O}_{a,b,x,y,z}(4).
\]
Similarly, we have to guarantee that the rank at $0$ of the
second jet map:
\[
\big(x,y,z\big)
\,\,\longmapsto\,\,
\Big(
P,\,\,
P_a,\,\,
P_b,\,\,
P_{aa},\,\,
P_{ab},\,\,
P_{bb}
\Big),
\]
is maximal equal to $3$, and this is also clear from its
expression:
\[
\big(x,y,z\big)
\,\,\longmapsto\,\,
\Big(
-\,z+{\rm O}_2,\,\,
-\,x+{\rm O}_2,\,\,
\zero{{\rm O}_2},\,\,
-\,2\underline{\beta}\,y+{\rm O}_2,\,\,
\zero{{\rm O}_2},\,\,
\zero{{\rm O}_2}
\Big).
\qedhere
\]
\endproof

\begin{Proposition}
If $\mathcal{M} \subset \K^{2+1} \times \K^{2+1}$ 
has Levi form(s) of constant rank $1$ and is 
$2$-nondegenerate both with respect to
parameters and to variables, then 
there exist normalized coordinates $(x,y,z,a,b,c)$ in which
its two equations read as:
\[
\aligned
z
&
\,=\,
c
+
xa
+
xxb
+
yaa
+
c\,{\rm O}_{x,y,a,b}(2)
+
{\rm O}_{x,y,a,b,c}(4),
\\
c
&
\,=\,
z
-
ax
-
aay
-
bxx
+
z\,{\rm O}_{a,b,x,y}(2)
+
{\rm O}_{a,b,x,y,z}(4).
\endaligned
\]
\end{Proposition}

\proof
Indeed, plain dilations make $\beta = 1 = \underline{\beta}$.
\endproof

\Section{\bf Determinantal Expressions of the Two 
$2$-Nondegeneracy Assumptions}
\label{determinantal-two-2-nondegeneracies}
\HEAD{{\ref{determinantal-two-2-nondegeneracies}}.~{\sf Determinantal 
Expressions of the Two $2$-Nondegeneracy Assumptions}
}{
Joël {\sc Merker} (Orsay)}

Now, and because it will be regularly used in what follows, it is
advisable to express in a concrete manner the two $2$-nondegeneracy
assumptions:
\[
\left\vert\!
\begin{array}{ccc}
Q_a & Q_b & Q_c
\\
Q_{xa} & Q_{xb} & Q_{xc}
\\
Q_{xxa} & Q_{xxb} & Q_{xxc}
\end{array}
\!\right\vert
\,\,\neq\,\,
0
\ \ \ \ \ \ \ \ \ \ \ \ \
\text{and}
\ \ \ \ \ \ \ \ \ \ \ \ \
0
\,\,\neq\,\,
\left\vert\!
\begin{array}{ccc}
P_x & P_y & P_z
\\
P_{ax} & P_{ay} & P_{az}
\\
P_{aax} & P_{aay} & P_{aaz}
\end{array}
\!\right\vert.
\]
Let us therefore abbreviate:
\[
\Delta(Q)
\,:=\,
\left\vert\!
\begin{array}{ccc}
Q_a & Q_b & Q_c
\\
Q_{xa} & Q_{xb} & Q_{xc}
\\
Q_{xxa} & Q_{xxb} & Q_{xxc}
\end{array}
\!\right\vert
\ \ \ \ \ \ \ \ \ \ \ \ \
\text{and}
\ \ \ \ \ \ \ \ \ \ \ \ \
\Box(P)
\,:=\,
\left\vert\!
\begin{array}{ccc}
P_x & P_y & P_z
\\
P_{ax} & P_{ay} & P_{az}
\\
P_{aax} & P_{aay} & P_{aaz}
\end{array}
\!\right\vert.
\]
We repeat here that the two hypotheses $\Delta(Q) \neq 0$ and
$\Box(P) \neq 0$ are independent, as show the two examples:
\[
\aligned
z
&
\,=\,
c
+
xa
+
xxb,
\\
\Delta(Q)
&
\,=\,
2,
\\
\Box(P)
&
\,=\,
0,
\endaligned
\ \ \ \ \ \ \ \ \ \ \ \ \ \ \ \ \ \ \ \ \ \ \ \ \ \
\aligned
z
&
\,=\,
c
+
ax
+
aay,
\\
\Delta(Q)
&
\,=\,
0,
\\
\Box(P)
&
\,=\,
2.
\endaligned
\]

\Section{\bf {\sc pde} System and Dual {\sc pde} System 
\\
Associated to Doubly $2$-Nondegenerate $(2,2,1)$ Para-CR Structures}
\label{pde-dual-pde-2-2-1}
\HEAD{{\ref{pde-dual-pde-2-2-1}}.~{\sf {\sc pde} System and Dual 
{\sc pde} System Associated
to Doubly $2$-Nondegenerate $(2,2,1)$ Para-CR Structures}
}{
Joël {\sc Merker} (Orsay)}

Now, in the equation $z = Q(x,y,a,b,c)$ of $\mathcal{M}$, 
we may view $z = z(x,y)$ as a function of $(x,y)$, and differentiate
it with respect to $x$, $xx$:
\[
\aligned
z
&
\,=\,
Q(x,y,a,b,c),
\\
z_x
&
\,=\,
Q_x(x,y,a,b,c),
\\
z_{xx}
&
\,=\,
Q_{xx}(x,y,a,b,c).
\endaligned
\]
The Jacobian of the map:
\[
\big(a,b,c\big)
\,\,\longmapsto\,\,
\Big(
Q(x,y,a,b,c),\,
Q_x(x,y,a,b,c),\,
Q_{xx}(x,y,a,b,c)
\Big),
\]
being precisely:
\[
\Delta(Q)
\,\neq\,
0,
\]
we can solve by means of the implicit function
theorem these three equations for $(a,b,c)$:
\[
\left[
\aligned
z
&
\,=\,
Q(x,y,a,b,c)
\\
z_x
&
\,=\,
Q_x(x,y,a,b,c)
\\
z_{xx}
&
\,=\,
Q_{xx}(x,y,a,b,c)
\endaligned\right.
\ \ \ \ \ \ \ \ \ \ \ \ \
\Longleftrightarrow
\ \ \ \ \ \ \ \ \ \ \ \ \
\left[
\aligned
a
&
\,=\,
A\big(x,y,z,z_x,z_{xx}\big)
\\
b
&
\,=\,
B\big(x,y,z,z_x,z_{xx}\big)
\\
c
&
\,=\,
C\big(x,y,z,z_x,z_{xx}\big),
\endaligned\right.
\]
and then replace these solutions in the two equations obtained
firstly by differentiating with respect to $y$:
\[
\aligned
z_y
&
\,=\,
Q_y(x,y,a,b,c)
\\
&
\,=\,
Q_y
\Big(
x,y,\,
A\big(x,y,z,z_x,z_{xx}\big),
B\big(x,y,z,z_x,z_{xx}\big),
C\big(x,y,z,z_x,z_{xx}\big)
\Big)
\\
&
\,=:\,
F\big(x,y,z,z_x,z_{xx}\big),
\endaligned
\] 
and secondly with respect to $xxx$:
\[
\aligned
z_{xxx}
&
\,=\,
Q_{xxx}
(x,y,a,b,c)
\\
&
\,=\,
Q_{xxx}
\Big(
x,y,\,
A\big(x,y,z,z_x,z_{xx}\big),
B\big(x,y,z,z_x,z_{xx}\big),
C\big(x,y,z,z_x,z_{xx}\big)
\Big)
\\
&
\,=:\,
H\big(x,y,z,z_x,z_{xx}\big).
\endaligned
\] 

\begin{center}
\scalebox{0.75}{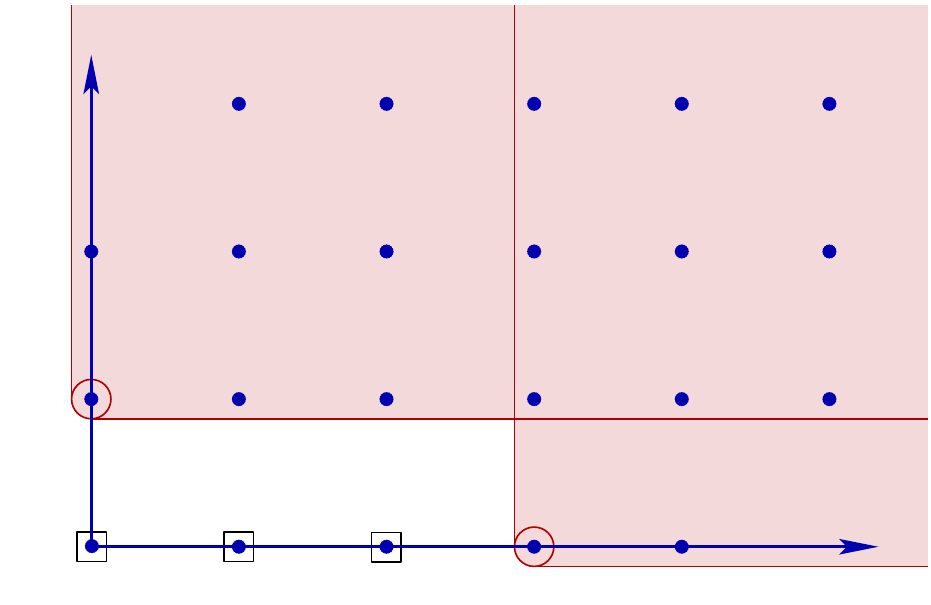}
\end{center}

Diagrammatically, we may
represent these two derivatives $z_y \longleftrightarrow
(0,1)$ and $z_{xxx} \longleftrightarrow (3,0)$,
and generally, all derivatives $z_{x^k y^l}$ by pairs of
integers $(k,l) \in \N \times \N$. 
And then beyond, all other derivatives
$z_{x^k y^l}$ with either $k \geqslant 3$ or $l \geqslant 1$
express in terms of the horizontal 
$(x,y,z, z_x, z_{xx})$, for instance:
\[
\aligned
z_{xy}
&
\,=\,
F_x
+
F_y
+
z_x\,F_z
+
z_{xx}\,F_{z_x}
+
z_{xxx}\,F_{z_{xx}}
\\
&
\,=\,
F_x
+
F_y
+
z_x\,F_z
+
z_{xx}\,F_{z_x}
+
H\,F_{z_{xx}}
\\
&
\,=:\,
\mathcal{F}_{1,1}
\big(x,y,z,z_x,z_{xx}\big),
\endaligned
\]
that is to say generally:
\[
z_{x^ky^l}
\,=\,
\mathcal{F}_{k,l}
\big(x,y,z,z_x,z_{xx}\big)
\eqno
{\scriptstyle{(k\,\geqslant\,3\,\,\text{\rm or}\,\,
l\,\geqslant\,1)}}.
\]
Later, we will see how to transfer equivalences of
submanifolds of solutions to equivalences of associated
{\sc pde} systems. 

For now, let us point out that we
have used only {\em one} assumption of $2$-nondegeneracy,
$\Delta(Q) \neq 0$. Unfortunately, we do not see what
is happening with the other $\Box(P) \neq 0$.

At least, this $\Box(P) \neq 0$ is useful to set
up a certain {\em dual}
{\sc pde} system. Indeed, in $c = P(a,b,x,y,z)$,
we can similarly view $c$ as a function of $(a,b)$, then solve
by means of the implicit function theorem:
\[
\left[
\aligned
c
&
\,=\,
P(a,b,x,y,z)
\\
c_a
&
\,=\,
P_a(a,b,x,y,z)
\\
c_{aa}
&
\,=\,
P_{aa}(a,b,x,y,z)
\endaligned\right.
\ \ \ \ \ \ \ \ \ \ \ \ \
\Longleftrightarrow
\ \ \ \ \ \ \ \ \ \ \ \ \
\left[
\aligned
x
&
\,=\,
X
\big(a,b,c,c_a,c_{aa}\big)
\\
y
&
\,=\,
Y
\big(a,b,c,c_a,c_{aa}\big)
\\
z
&
\,=\,
Z
\big(a,b,c,c_a,c_{aa}\big),
\endaligned\right.
\]
thanks to the fact that the Jacobian determinant of the map:
\[
(x,y,z)
\,\,\longmapsto\,\,
\Big(
P(a,b,x,y,z),\,
P_a(a,b,x,y,z),\,
P_{aa}(a,b,x,y,z)
\Big)
\]
is precisely:
\[
\Box(P)
\,\neq\,
0,
\]
and then replace\,\,---\,\,very similarly!\,\,---\,\,in 
the other two relevant partial derivatives:
\[
\aligned
c_b
&
\,=\,
P_b(a,b,x,y,z)
\\
&
\,=\,
P_b
\Big(
a,b,\,
X\big(a,b,c,c_a,c_{aa}\big),\,
Y\big(a,b,c,c_a,c_{aa}\big),\,
Z\big(a,b,c,c_a,c_{aa}\big)
\Big)
\\
&
\,=:\,
E\big(a,b,c,c_a,c_{aa}\big),
\\
c_{aaa}
&
\,=\,
P_{aaa}(a,b,x,y,z)
\\
&
\,=\,
P_{aaa}
\Big(
a,b,\,
X\big(a,b,c,c_a,c_{aa}\big),\,
Y\big(a,b,c,c_a,c_{aa}\big),\,
Z\big(a,b,c,c_a,c_{aa}\big)
\Big)
\\
&
\,=:\,
G\big(a,b,c,c_a,c_{aa}\big).
\endaligned
\]

However, the symmetric question arises here: how to view the second
$2$-nondegeneracy condition $\Delta(Q) \neq 0$ in this
dual {\sc pde} system?

The (symmetric) answer to both questions is accessible, but it
needs little preliminaries.

\Section{\bf Transfer of Derivations 
$\mathcal{M} \longleftrightarrow 
\text{\sc pde}_\variablessmall(\mathcal{M})$
and 
$\mathcal{M} \longleftrightarrow 
\text{\sc pde}_\parameterssmall(\mathcal{M})$
}
\label{transfer-vector-fields-M-PDE}
\HEAD{{\ref{transfer-vector-fields-M-PDE}}.~{\sf Transfer of 
Derivations 
$\mathcal{M} \longleftrightarrow 
\text{\sc pde}_\variablessmall(\mathcal{M})$
and 
$\mathcal{M} \longleftrightarrow 
\text{\sc pde}_\parameterssmall(\mathcal{M})$}
}{
Joël {\sc Merker} (Orsay)}

We shall treat only the transfer of derivations (of vector fields,
of frames) from the solution space $\mathcal{M}$
equipped with coordinates $(x,y,a,b,c)$ to 
the jet space equipped with coordinates $(x,y,z,z_x,z_{xx})$:
\[
\big(x,y,a,b,c\big)
\xleftrightarrow[{\rule[0pt]{50pt}{0pt}}]{}
\big(x,y,z,z_x,z_{xx}\big).
\]
Given a function $\Gaux = \Gaux (x,y, a, b, c)$ on the left
space and a function $\Faux = \Faux (x,y,z,z_x,z_{xx})$ on
the right space,
the transfer of basic coordinate vector fields is obtained
by differentiating the composition identity:
\[
\Faux
\big(x,y,z,z_x,z_{xx}\big)
\,\equiv\,
\Gaux
\Big(
x,y,\,
A\big(x,y,z,z_x,z_{xx}\big),\,
B\big(x,y,z,z_x,z_{xx}\big),\,
C\big(x,y,z,z_x,z_{xx}\big)
\Big),
\]
with respect to all five variables:
\[
\aligned
\frac{\partial\Faux}{\partial x}
&
\,=\,
\frac{\partial\Gaux}{\partial x}
\ \ \ \ \ \ \ \ \ \ \ \
+
A_x\,
\frac{\partial\Gaux}{\partial a}
+
B_x\,
\frac{\partial\Gaux}{\partial b}
+
C_x\,
\frac{\partial\Gaux}{\partial c},
\\
\frac{\partial\Faux}{\partial y}
&
\,=\,
\ \ \ \ \ \ \ \ \ \ \ \
\frac{\partial\Gaux}{\partial y}
+
A_y\,
\frac{\partial\Gaux}{\partial a}
+
B_y\,
\frac{\partial\Gaux}{\partial b}
+
C_y\,
\frac{\partial\Gaux}{\partial c},
\\
\frac{\partial\Faux}{\partial z}
&
\,=\,
\ \ \ \ \ \ \ \ \ \ \ \
\ \ \ \ \ \ \ \ \ \ \ \
A_z\,
\frac{\partial\Gaux}{\partial a}
+
B_z\,
\frac{\partial\Gaux}{\partial b}
+
C_z\,
\frac{\partial\Gaux}{\partial c},
\\
\frac{\partial\Faux}{\partial z_x}
&
\,=\,
\ \ \ \ \ \ \ \ \ \ \ \
\ \ \ \ \ \ \ \ \ \ \ \
A_{z_x}\,
\frac{\partial\Gaux}{\partial a}
+
B_{z_x}\,
\frac{\partial\Gaux}{\partial b}
+
C_{z_x}\,
\frac{\partial\Gaux}{\partial c},
\\
\frac{\partial\Faux}{\partial z_{xx}}
&
\,=\,
\ \ \ \ \ \ \ \ \ \ \ \
\ \ \ \ \ \ \ \ \ \ \ \
A_{z_{xx}}\,
\frac{\partial\Gaux}{\partial a}
+
B_{z_{xx}}\,
\frac{\partial\Gaux}{\partial b}
+
C_{z_{xx}}\,
\frac{\partial\Gaux}{\partial c}.
\endaligned
\]
Yet, the appearing coefficients
$A_\smallbullet$, $B_\smallbullet$, $C_\smallbullet$ 
are not expressed in terms of the left coordinates
$(x,y,a,b,c)$. 

To re-express them as required, 
from the three identically satisfied equations:
\[
\aligned
a
&
\,\equiv\,
A
\Big(
x,y,\,
Q(x,y,a,b,c),\,
Q_x(x,y,a,b,c),\,
Q_{xx}(x,y,a,b,c)
\Big),
\\
b
&
\,\equiv\,
B
\Big(
x,y,\,
Q(x,y,a,b,c),\,
Q_x(x,y,a,b,c),\,
Q_{xx}(x,y,a,b,c)
\Big),
\\
c
&
\,\equiv\,
C
\Big(
x,y,\,
Q(x,y,a,b,c),\,
Q_x(x,y,a,b,c),\,
Q_{xx}(x,y,a,b,c)
\Big),
\endaligned
\]
differentiations with respect to $x$, $y$, $a$, $b$, $c$ provide
firstly:
\[
\aligned
0
&
\,=\,
A_x
\ \ \ \ \ \ \ \
+
Q_x\,A_z
+
Q_{xx}\,A_{z_x}
+
Q_{xxx}\,A_{z_{xx}},
\\
0
&
\,=\,
\ \ \ \ \ \ \ \
A_y
+
Q_y\,A_z
+
Q_{xy}\,A_{z_x}
+
Q_{xxy}\,A_{z_{xx}},
\\
1
&
\,=\,
\ \ \ \ \ \ \ \ \ \ \ \ \ \ \ \ \ \
Q_a\,A_z
+
Q_{xa}\,A_{z_x}
+
Q_{xxa}\,A_{z_{xx}},
\\
0
&
\,=\,
\ \ \ \ \ \ \ \ \ \ \ \ \ \ \ \ \ \
Q_b\,A_z
+
Q_{xb}\,A_{z_x}
+
Q_{xxb}\,A_{z_{xx}},
\\
0
&
\,=\,
\ \ \ \ \ \ \ \ \ \ \ \ \ \ \ \ \ \
Q_c\,A_z
+
Q_{xc}\,A_{z_x}
+
Q_{xxc}\,A_{z_{xx}},
\endaligned
\]
secondly:
\[
\aligned
0
&
\,=\,
B_x
\ \ \ \ \ \ \ \
+
Q_x\,B_z
+
Q_{xx}\,B_{z_x}
+
Q_{xxx}\,B_{z_{xx}},
\\
0
&
\,=\,
\ \ \ \ \ \ \ \
B_y
+
Q_y\,B_z
+
Q_{xy}\,B_{z_x}
+
Q_{xxy}\,B_{z_{xx}},
\\
0
&
\,=\,
\ \ \ \ \ \ \ \ \ \ \ \ \ \ \ \ \ \
Q_a\,B_z
+
Q_{xa}\,B_{z_x}
+
Q_{xxa}\,B_{z_{xx}},
\\
1
&
\,=\,
\ \ \ \ \ \ \ \ \ \ \ \ \ \ \ \ \ \
Q_b\,B_z
+
Q_{xb}\,B_{z_x}
+
Q_{xxb}\,B_{z_{xx}},
\\
0
&
\,=\,
\ \ \ \ \ \ \ \ \ \ \ \ \ \ \ \ \ \
Q_c\,B_z
+
Q_{xc}\,B_{z_x}
+
Q_{xxc}\,B_{z_{xx}},
\endaligned
\]
and thirdly:
\[
\aligned
0
&
\,=\,
C_x
\ \ \ \ \ \ \ \
+
Q_x\,C_z
+
Q_{xx}\,C_{z_x}
+
Q_{xxx}\,C_{z_{xx}},
\\
0
&
\,=\,
\ \ \ \ \ \ \ \
C_y
+
Q_y\,C_z
+
Q_{xy}\,C_{z_x}
+
Q_{xxy}\,C_{z_{xx}},
\\
0
&
\,=\,
\ \ \ \ \ \ \ \ \ \ \ \ \ \ \ \ \ \
Q_a\,C_z
+
Q_{xa}\,C_{z_x}
+
Q_{xxa}\,C_{z_{xx}},
\\
0
&
\,=\,
\ \ \ \ \ \ \ \ \ \ \ \ \ \ \ \ \ \
Q_b\,C_z
+
Q_{xb}\,C_{z_x}
+
Q_{xxb}\,C_{z_{xx}},
\\
1
&
\,=\,
\ \ \ \ \ \ \ \ \ \ \ \ \ \ \ \ \ \
Q_c\,C_z
+
Q_{xc}\,C_{z_x}
+
Q_{xxc}\,C_{z_{xx}}.
\endaligned
\]
Resolutions of the three $3 \times 3$ linear systems 
consisting each time of the last three lines provide:
\[
\footnotesize
\aligned
A_z
&
\,=\,
\frac{1}{\Delta(Q)}\,
\left\vert\!
\begin{array}{cc}
Q_{xb} & Q_{xxb}
\\
Q_{xc} & Q_{xxc}
\end{array}
\!\right\vert,
\ \ \ \ \ \ \ \ \ \ \ \ \ \ \ \ \ \
A_{z_x}
\,=\,
-\,
\frac{1}{\Delta(Q)}\,
\left\vert\!
\begin{array}{cc}
Q_b & Q_{xxb}
\\
Q_c & Q_{xxc}
\end{array}
\!\right\vert,
\ \ \ \ \ \ \ \ \ \ \ \ \ \ \ \ \ \
A_{z_{xx}}
\,=\,
\frac{1}{\Delta(Q)}\,
\left\vert\!
\begin{array}{cc}
Q_b & Q_{xb}
\\
Q_c & Q_{xc}
\end{array}
\!\right\vert,
\\
B_z
&
\,=\,
-\,
\frac{1}{\Delta(Q)}\,
\left\vert\!
\begin{array}{cc}
Q_{xa} & Q_{xxa}
\\
Q_{xc} & Q_{xxc}
\end{array}
\!\right\vert,
\ \ \ \ \ \ \ \ \ \ \ \ \ \
B_{z_x}
\,=\,
\frac{1}{\Delta(Q)}\,
\left\vert\!
\begin{array}{cc}
Q_a & Q_{xxa}
\\
Q_c & Q_{xxc}
\end{array}
\!\right\vert,
\ \ \ \ \ \ \ \ \ \ \ \ \ \ \ \ \ \ \ \ \ \
B_{z_{xx}}
\,=\,
-\,
\frac{1}{\Delta(Q)}\,
\left\vert\!
\begin{array}{cc}
Q_a & Q_{xa}
\\
Q_c & Q_{xc}
\end{array}
\!\right\vert,
\\
C_z
&
\,=\,
\frac{1}{\Delta(Q)}\,
\left\vert\!
\begin{array}{cc}
Q_{xa} & Q_{xxa}
\\
Q_{xb} & Q_{xxb}
\end{array}
\!\right\vert,
\ \ \ \ \ \ \ \ \ \ \ \ \ \ \ \ \ \ 
C_{z_x}
\,=\,
-\,
\frac{1}{\Delta(Q)}\,
\left\vert\!
\begin{array}{cc}
Q_a & Q_{xxa}
\\
Q_b & Q_{xxb}
\end{array}
\!\right\vert,
\ \ \ \ \ \ \ \ \ \ \ \ \ \ \ \ \ \ \ 
C_{z_{xx}}
\,=\,
\frac{1}{\Delta(Q)}\,
\left\vert\!
\begin{array}{cc}
Q_a & Q_{xa}
\\
Q_b & Q_{xb}
\end{array}
\!\right\vert,
\endaligned
\]
and consequently, the transfer of the last three (among five)
vector fields is:
\[
\footnotesize
\aligned
\frac{\partial}{\partial z}
&
\,=\,
\frac{
\left\vert\!
\begin{array}{cc}
Q_{xb} & Q_{xxb}
\\
Q_{xc} & Q_{xxc}
\end{array}
\!\right\vert}{\Delta(Q)}\,
\frac{\partial}{\partial a}
-
\frac{
\left\vert\!
\begin{array}{cc}
Q_{xa} & Q_{xxa}
\\
Q_{xc} & Q_{xxc}
\end{array}
\!\right\vert
}{\Delta(Q)}\,
\frac{\partial}{\partial b}
+
\frac{
\left\vert\!
\begin{array}{cc}
Q_{xa} & Q_{xxa}
\\
Q_{xb} & Q_{xxb}
\end{array}
\!\right\vert
}{\Delta(Q)}\,
\frac{\partial}{\partial c},
\\
\frac{\partial}{\partial z_x}
&
\,=\,
-\,
\frac{
\left\vert\!
\begin{array}{cc}
Q_b & Q_{xxb}
\\
Q_c & Q_{xxc}
\end{array}
\!\right\vert
}{\Delta(Q)}\,
\frac{\partial}{\partial a}
+
\frac{
\left\vert\!
\begin{array}{cc}
Q_a & Q_{xxa}
\\
Q_c & Q_{xxc}
\end{array}
\!\right\vert
}{\Delta(Q)}\,
\frac{\partial}{\partial b}
-
\frac{
\left\vert\!
\begin{array}{cc}
Q_a & Q_{xxa}
\\
Q_b & Q_{xxb}
\end{array}
\!\right\vert
}{\Delta(Q)}\,
\frac{\partial}{\partial c},
\\
\frac{\partial}{\partial z_{xx}}
&
\,=\,
\frac{
\left\vert\!
\begin{array}{cc}
Q_b & Q_{xb}
\\
Q_c & Q_{xc}
\end{array}
\!\right\vert
}{\Delta(Q)}\,
\frac{\partial}{\partial a}
-
\frac{
\left\vert\!
\begin{array}{cc}
Q_a & Q_{xa}
\\
Q_c & Q_{xc}
\end{array}
\!\right\vert
}{\Delta(Q)}\,
\frac{\partial}{\partial b}
+
\frac{
\left\vert\!
\begin{array}{cc}
Q_a & Q_{xa}
\\
Q_b & Q_{xb}
\end{array}
\!\right\vert
}{\Delta(Q)}\,
\frac{\partial}{\partial c}.
\endaligned
\]

\begin{Proposition}
\label{Prp-F-zxx-zero}
If the submanifold of solutions $\mathcal{M} \subset
\K^2 \times \K^1 \times \K^2 \times \K^1$ has (degenerate)
Levi form of constant rank $1$ and if it is $2$-nondegenerate
with respect to parameters, then in
its associated system $\pde_\variablessmall(\mathcal{M})$:
\[
\aligned
z_y
&
\,=\,
F
\big(
x,y,z,z_x,z_{xx}
\big),
\\
z_{xxx}
&
\,=\,
H
\big(
x,y,z,z_x,z_{xx}
\big),
\endaligned
\]
the function $F$ is independent of $z_{xx}$:
\[
0
\,\equiv\,
F_{z_{xx}}.
\]
\end{Proposition}

\proof
By construction:
\[
F
\big(
x,y,z,z_x,z_{xx}
\big)
\,\,=:\,\,
Q_y
\Big(
x,y,\,
A\big(x,y,z,z_x,z_{xx}\big),\,
B\big(x,y,z,z_x,z_{xx}\big),\,
C\big(x,y,z,z_x,z_{xx}\big)
\Big),
\]
whence a differentiation with respect to $z_{xx}$ makes re-appear
the Levi determinant ${\sf L}_{3\times 3}(Q) \equiv 0$:
\begin{align}
F_{z_{xx}}
&
\,=\,
A_{z_{xx}}\,Q_{ya}
+
B_{z_{xx}}\,Q_{yb}
+
C_{z_{xx}}\,Q_{yc}
\notag
\\
&
\,=\,
\frac{
\left\vert\!
\begin{array}{cc}
Q_b & Q_{xb}
\\
Q_c & Q_{xc}
\end{array}
\!\right\vert
}{\Delta(Q)}\,
Q_{ya}
-
\frac{
\left\vert\!
\begin{array}{cc}
Q_a & Q_{xa}
\\
Q_c & Q_{xc}
\end{array}
\!\right\vert
}{\Delta(Q)}\,
Q_{yb}
+
\frac{
\left\vert\!
\begin{array}{cc}
Q_a & Q_{xa}
\\
Q_b & Q_{xb}
\end{array}
\!\right\vert
}{\Delta(Q)}\,
Q_{yc}
\notag
\\
&
\,=\,
\frac{1}{\Delta(Q)}\,\,
\left\vert\!
\begin{array}{ccc}
Q_a & Q_b & Q_c
\\
Q_{xa} & Q_{xb} & Q_{xc}
\\
Q_{ya} & Q_{yb} & Q_{yc}
\end{array}
\!\right\vert
\notag
\\
&
\,\equiv\,
0.
\qedhere
\end{align}
\endproof

Of course, we have in a symmetric way for 
$\pde_\parameterssmall (\mathcal{M})$:
\[
0
\,\equiv\,
E_{c_{aa}}.
\]

One must observe that if, on the contrary, on would have:
\[
F_{z_{xx}}
\,\not\equiv\,
0,
\]
then in a neighborhood of a generic point 
near which $F_{z_{xx}} \neq 0$,
one would be able to solve for $z_{xx}$
in the first partial differential equation $z_y = F$:
\[
z_{xx}
\,=\,
\Lambda\,
\big(
x,y,z,z_x,z_y
\big),
\]
and further replacements and differentiations 
would show that one 
would come to a {\sc pde} system of the kind already
studied by Hachtroudi~{\cite{Hachtroudi-1937}}:
\[
\aligned
z_{xx}
&
\,=\,
\Faux_{2,0}
\big(
x,y,z,z_x,z_y
\big),
\\
z_{xy}
&
\,=\,
\Faux_{1,1}
\big(
x,y,z,z_x,z_y
\big),
\\
z_{yy}
&
\,=\,
\Faux_{0,2}
\big(
x,y,z,z_x,z_y
\big).
\endaligned
\]
This is equivalent to the observation that, on a hypersurface
$M^5 \subset \C^3$ whose Levi form is not identically degenerate,
at a generic point, the reduction to an $\{e\}$-structure
and to a Cartan connection was already studied by Chern-Moser.
So we will definitely assume $F_{z_{xx}} \equiv 0$ up
to the end of our considerations.

\Section{\bf The Two Kernels of the Two Levi Forms}
\label{two-kernels-two-Levi-forms}
\HEAD{{\ref{two-kernels-two-Levi-forms}}.~{\sf The two Kernels 
of the Two Levi Forms}
}{
Joël {\sc Merker} (Orsay)}

Recall that with the two transversal vector fields: 
\[
\mathcal{T}_b
\,=\,
\frac{\partial}{\partial b}
\ \ \ \ \ \ \ \ \ \ \ \ \
\text{and}
\ \ \ \ \ \ \ \ \ \ \ \ \
\mathcal{U}_y
\,=\,
\frac{\partial}{\partial y},
\]
we have introduced
two $1$-forms $\rho$ and $\sigma$ on $\mathcal{M}$ which
satisfy:
\[
\rho\big(\mathcal{T}_b\big)
\,\equiv\,
1
\ \ \ \ \ \ \ \ \ \ \ \ \
\text{and}
\ \ \ \ \ \ \ \ \ \ \ \ \
\sigma\big(\mathcal{U}_y\big)
\,\equiv\,
1,
\]
and that:
\[
\Levi_\variablessmall
\big(
\mathcal{K},\mathcal{L}
\big)
\,=\,
\rho
\big(
[\mathcal{K},\mathcal{L}]
\big)
\ \ \ \ \ \ \ \ \ \ \ \ \
\text{and}
\ \ \ \ \ \ \ \ \ \ \ \ \
\Levi_\parameterssmall
\big(
\mathcal{L},
\mathcal{K}
\big)
\,=\,
\sigma
\big(
[\mathcal{L},\mathcal{K}]
\big).
\]

By assumption, the two matrices of these two Levi forms:
\[
\footnotesize
\aligned
\def\arraystretch{1.5}
\left\vert\!
\begin{array}{ccc}
\frac{-Q_cQ_{xa}+Q_aQ_{xc}}{Q_c\,Q_c}
&
\frac{-Q_cQ_{xb}+Q_bQ_{xc}}{Q_c\,Q_c}
\\
\frac{-Q_cQ_{ya}+Q_aQ_{yc}}{Q_c\,Q_c}
&
\frac{-Q_cQ_{yb}+Q_bQ_{yc}}{Q_c\,Q_c}
\end{array}
\!\right\vert
\big(x,y,a,b,c\big)
\ \ \ 
\text{and}
\ \ \ 
\def\arraystretch{1.5}
\left\vert\!
\begin{array}{ccc}
\frac{-P_zP_{ax}+P_xP_{az}}{P_z\,P_z}
&
\frac{-P_zP_{ay}+P_xP_{ay}}{P_z\,P_z}
\\
\frac{-P_zP_{bx}+P_xP_{bz}}{P_z\,P_z}
&
\frac{-P_zP_{by}+P_yP_{bz}}{P_z\,P_z}
\end{array}
\!\right\vert
\big(a,b,x,y,z\big),
\endaligned
\]
have rank $1$ at every point, and their two upper-left
$(1,1)$-entries are nowhere vanishing, since we have already
normalized:
\[
Q
\,=\,
c
+
xa
+
{\rm O}_3
\ \ \ \ \ \ \ \ \ \ \ \ \
\Longleftrightarrow
\ \ \ \ \ \ \ \ \ \ \ \ \
P
\,=\,
-\,z
+
ax
+
{\rm O}_3,
\]
or because this can always be achieved after an allowed change of
coordinates.

\begin{Definition}
The two kernels of the two Levi forms are:
\[
\aligned
\Ker\,\Levi_\variablessmall
&
\,:=\,
\Big\{
\mathcal{K}
\in
\Gamma\big(T^\variablessmall\mathcal{M}\big)
\colon\,
0
=
\Levi_\variablessmall
\big(
\mathcal{K},
\mathcal{L}
\big),
\ \ 
\forall\,
\mathcal{L}
\in
\Gamma
\big(
T^\parameterssmall\mathcal{M}
\big)
\Big\}, 
\\
\Ker\,\Levi_\parameterssmall
&
\,:=\,
\Big\{
\mathcal{L}
\in
\Gamma\big(T^\parameterssmall\mathcal{M}\big)
\colon\,
0
=
\Levi_\parameterssmall
\big(
\mathcal{L},
\mathcal{K}
\big),
\ \ 
\forall\,
\mathcal{K}
\in
\Gamma
\big(
T^\variablessmall\mathcal{M}
\big)
\Big\}.
\endaligned
\]
\end{Definition}

Since the above two Levi matrices have constant rank $1$ all over
$\mathcal{M}$, their two kernels define two analytic (smooth) rank $1$
distributions. Two vector field generators can be explicitly written
in terms of two quotients of the entries of the first rows of the
two Levi matrices above, as follows.

\begin{Lemma}
Two natural generators of these two Levi forms
kernels are the two vector fields:
\[
\aligned
\mathcal{K}_\kersmall^\variablessmall
&
\,:=\,
-\,
\frac{-P_zP_{ay}+P_yP_{az}}{-P_zP_{ax}+P_xP_{az}}\,
\bigg(
\frac{\partial}{\partial x}
-
\frac{P_x}{P_z}\,
\frac{\partial}{\partial z}
\bigg)
+
\frac{\partial}{\partial y}
-
\frac{P_y}{P_z}\,
\frac{\partial}{\partial z},
\\
\mathcal{L}_\kersmall^\parameterssmall
&
\,:=\,
-\,
\frac{-Q_cQ_{xb}+Q_bQ_{xc}}{-Q_cQ_{xa}+Q_aQ_{xc}}\,
\bigg(
\frac{\partial}{\partial a}
-
\frac{Q_a}{Q_c}\,\frac{\partial}{\partial c}
\bigg)
+
\frac{\partial}{\partial b}
-
\frac{Q_b}{Q_c}\,\frac{\partial}{\partial c}.
\endaligned
\]
\end{Lemma}

\proof
It suffices to verify that:
\[
\big[
{\textstyle{\frac{\partial}{\partial a}}},\,
\mathcal{K}_\kersmall^\variablessmall
\big]
\,\,\equiv\,\,
0
\,\,\equiv\,\,
\big[
{\textstyle{\frac{\partial}{\partial b}}},\,
\mathcal{K}_\kersmall^\variablessmall
\big],
\]
which is done by a direct computation,
the first equality to zero being true for free, 
while the second one comes from the Levi degeneracy assumption
${\sf L}_{3\times 3}(P) \equiv 0$.

The other pair of annihilations is checked
in a quite symmetric manner:
\[
\big[
{\textstyle{\frac{\partial}{\partial x}}},\,
\mathcal{L}_\kersmall^\parameterssmall
\big]
\,\,\equiv\,\,
0
\,\,\equiv\,\,
\big[
{\textstyle{\frac{\partial}{\partial y}}},\,
\mathcal{L}_\kersmall^\parameterssmall
\big].
\qedhere
\]
\endproof

For simplicity, we will sometimes write:
\[
\mathcal{K}_\kersmall^\variablessmall
\,=\,
\kaux\,
\mathcal{K}_x
+
\mathcal{K}_y
\ \ \ \ \ \ \ \ \ \ \ \ \
\text{and}
\ \ \ \ \ \ \ \ \ \ \ \ \
\mathcal{L}_\kersmall^\parameterssmall
\,=\,
\laux\,
\mathcal{L}_a
+
\mathcal{L}_b,
\] 
after giving a name to the two Levi entries quotients in question
(mind the two minus signs):
\[
\kaux
\,:=\,
-\,
\frac{-P_zP_{ay}+P_yP_{az}}{-P_zP_{ax}+P_xP_{az}}\,
\ \ \ \ \ \ \ \ \ \ \ \ \
\text{and}
\ \ \ \ \ \ \ \ \ \ \ \ \
\laux
\,:=\,
-\,
\frac{-Q_cQ_{xb}+Q_bQ_{xc}}{-Q_cQ_{xa}+Q_aQ_{xc}}.
\]

\Section{\bf Reexpressions of $F_{z_x}$ in Terms of $Q$ and
in Terms of $P$}
\label{reexpressions-F-z-x-with-Q-P}
\HEAD{{\ref{reexpressions-F-z-x-with-Q-P}}.~{\sf Reexpressions 
of $F_{z_x}$ in Terms of $Q$ and in Terms of $P$}
}{
Joël {\sc Merker} (Orsay)}

By differentiating with respect to $z_x$ the function 
$F$ of $\pde_\variablessmall (\mathcal{M})$ from its definition:
\[
F
\big(x,y,z,z_x,z_{xx}\big)
\,:=\,
Q_y
\Big(
x,y,\,
A\big(x,y,z,z_x,z_{xx}\big),\,
B\big(x,y,z,z_x,z_{xx}\big),\,
C\big(x,y,z,z_x,z_{xx}\big)
\Big),
\]
we get:
\[
\aligned
F_{z_x}
&
\,=\,
A_{z_x}\,
Q_{ya}
+
B_{z_x}\,
Q_{yb}
+
C_{z_x}\,
Q_{yc}
\\
&
\,=\,
-\,
\frac{
\left\vert\!
\begin{array}{cc}
Q_b & Q_{xxb}
\\
Q_c & Q_{xxc}
\end{array}
\!\right\vert
}{
\Delta(Q)}\,
Q_{ya}
+
\frac{
\left\vert\!
\begin{array}{cc}
Q_a & Q_{xxa}
\\
Q_c & Q_{xxc}
\end{array}
\!\right\vert
}{
\Delta(Q)}\,
Q_{yb}
-
\frac{
\left\vert\!
\begin{array}{cc}
Q_a & Q_{xxa}
\\
Q_b & Q_{xxb}
\end{array}
\!\right\vert
}{
\Delta(Q)}\,
Q_{yc}
\\
&
\,=\,
\frac{
\left\vert\!
\begin{array}{ccc}
Q_a & Q_b & Q_c
\\
Q_{ya} & Q_{yb} & Q_{yc}
\\
Q_{xxa} & Q_{xxb} & Q_{xxc}
\end{array}
\!\right\vert
}{
\left\vert\!
\begin{array}{ccc}
Q_a & Q_b & Q_c
\\
Q_{xa} & Q_{xb} & Q_{xc}
\\
Q_{xxa} & Q_{xxb} & Q_{xxc}
\end{array}
\!\right\vert
}.
\endaligned
\]
So we have a quotient of two determinants, 
one in the denominator which is nowhere vanishing, 
and one in the numerator in which only the second row differs.
A surprising (and useful) simplification occurs.

\begin{Lemma}
\label{Lm-F-zx}
One has in fact:
\[
F_{z_x}
\,=\,
\frac{-Q_cQ_{ya}+Q_aQ_{yc}}{-Q_cQ_{xa}+Q_aQ_{xc}}.
\]
\end{Lemma}

\proof

An expansion of both determinants along their last (identical)
rows gives:
\[
\aligned
F_{z_x}
&
\,\,=\,\,
\frac{
Q_{xxa}\,
\left\vert\!
\begin{array}{cc}
Q_b & Q_c
\\
Q_{yb} & Q_{yc}
\end{array}
\!\right\vert
-
Q_{xxb}\,
\left\vert\!
\begin{array}{cc}
Q_a & Q_c
\\
Q_{ya} & Q_{yc}
\end{array}
\!\right\vert
+
Q_{xxc}\,
\left\vert\!
\begin{array}{cc}
Q_a & Q_b
\\
Q_{ya} & Q_{yb}
\end{array}
\!\right\vert
}{
Q_{xxa}\,
\left\vert\!
\begin{array}{cc}
Q_b & Q_c
\\
Q_{xb} & Q_{xc}
\end{array}
\!\right\vert
-
Q_{xxb}\,
\left\vert\!
\begin{array}{cc}
Q_a & Q_c
\\
Q_{xa} & Q_{xc}
\end{array}
\!\right\vert
+
Q_{xxc}\,
\left\vert\!
\begin{array}{cc}
Q_a & Q_b
\\
Q_{xa} & Q_{xb}
\end{array}
\!\right\vert
}
\\
&
\overset{\text{\bf ?}}{\,\,=\,\,}
\frac{-Q_cQ_{ya}+Q_aQ_{yc}}{-Q_cQ_{xa}+Q_aQ_{xc}}.
\endaligned
\]
Then by eliminating (cross-producting) the denominators in this 
equality $\overset{\text{\bf ?}}{\,=\,}$ under question,
one recovers a multiple of
${\sf L}_{3 \times 3}(Q) \equiv 0$.
\endproof

\begin{Proposition}
\label{Proposition-k-equal-F-zx}
One has:
\[
\aligned
F_{z_x}
&
\,=\,
\frac{-P_zP_{ay}+P_yP_{az}}{-P_zP_{ax}+P_xP_{az}}
\\
&
\,=\,
-\,
\kaux.
\endaligned
\]
\end{Proposition}

\proof
Indeed, the very useful 
Lemma~{\ref{Lemma-transfer-from-Q-to-P}} applies:
\begin{align}
F_{z_x}
&
\,=\,
\frac{
\frac{-Q_cQ_{ya}+Q_aQ_{yc}}{Q_c\,Q_c}}{
\frac{-Q_cQ_{xa}+Q_aQ_{xc}}{Q_c\,Q_c}}
\notag
\\
&
\,=\,
\frac{-P_z\,\frac{-P_zP_{ay}+P_yP_{az}}{P_zP_z}
}{
-\,P_z\,\frac{-P_zP_{ax}+P_xP_{az}}{P_zP_z}}.
\qedhere
\end{align}
\endproof

Lastly, two definitions of invariant higher order Levi forms,
analogous to the {\sl Freeman form} in CR Geometry
{\cite[Section~9]{Merker-Pocchiola-Sabzevari-2013-5-CR-II}}, 
exist, but
because the jet theory is more general, we will dispense ourselves of
introducing these two concepts. Indeed,
the nondegeneracy of one of these two
`Freeman forms' could then be expressed as the nonvanishing
$\frac{\partial}{ \partial a} (\kaux) \neq 0$, entirely analogous to
the nonvanishing condition $\overline{\mathcal{L}}_1 (k) \neq 0$ which
was central in~{\cite{Merker-Pocchiola-2018}}, 
but we recover one of our two favorite nonvanishing $3
\times 3$ determinants anyway.

\begin{Lemma}
\label{Lm-da-F-zx}
One has the nowhere vanishing invariant expression:
\[
\frac{\partial}{\partial a}\,
\big(
F_{z_x}
\big)
\,=\,
P_z\,
\frac{
\left\vert\!
\begin{array}{ccc}
P_x & P_y & P_z
\\
P_{ax} & P_{ay} & P_{az}
\\
P_{aax} & P_{aay} & P_{aaz}
\end{array}
\!\right\vert
}{
\big(
-\,P_z\,P_{ax}+P_x\,P_{az}
\big)^2
}
\,\,\neq\,\,
0.
\]
\end{Lemma}

\proof
This just amounts to differentiate and to reorganize
properly:
\[
\frac{\partial}{\partial a}\,
\bigg[
\frac{-P_zP_{ay}+P_yP_{az}}{-P_zP_{ax}+P_xP_{az}}
\bigg].
\qedhere
\]
\endproof

In summary, the system $\pde_\variablessmall (\mathcal{M})$
and in its dual system
$\pde_\parameterssmall (\mathcal{M})$:
\[
\left[
\aligned
z_y
&
\,=\,
F
\big(
x,y,z,z_x,z_{xx}
\big),
\\
z_{xxx}
&
\,=\,
H
\big(
x,y,z,z_x,z_{xx}
\big),
\endaligned\right.
\ \ \ \ \ \ \ \ \ \ \ \ \
\text{and}
\ \ \ \ \ \ \ \ \ \ \ \ \
\left[
\aligned
c_b
&
\,=\,
E
\big(
a,b,c,c_a,c_{aa}
\big),
\\
c_{aaa}
&
\,=\,
G
\big(
a,b,c,c_a,c_{aa}
\big),
\endaligned\right.
\]
which come from a submanifold of solutions $\mathcal{M}$
whose two Levi forms have constant rank $1$, and
which is $2$-nondegenerate both with respect to
parameters and to variables, must satisfy:
\[
F_{z_{xx}}
\,\equiv\,
0
\ \ \ \ \ \ \ \ \ \ \ \ \
\text{and}
\ \ \ \ \ \ \ \ \ \ \ \ \
0
\,\equiv\,
E_{c_{aa}}.
\]

\Section{\bf Reexpression of the Hypothesis of
$2$-Nondegeneracy 
\\
with Respect to Variables as $F_{z_xz_x} \neq 0$}
\label{reexpression-2-nondegeneracy-F-pp-neq-0}
\HEAD{{\ref{reexpression-2-nondegeneracy-F-pp-neq-0}}.~{\sf 
Reexpression of the Hypothesis of $2$-Nondegeneracy with 
Respect to Variables as $F_{z_xz_x} \neq 0$}
}{
Joël {\sc Merker} (Orsay)}

To conclude this preliminary trip, we must answer the
question raised in 
Section~{\ref{pde-dual-pde-2-2-1}}, 
which we can now formulate a bit more precisely
(of course, a symmetric question can also be fomulated). 

\begin{Question}
{\sl How to view the hypothesis of $2$-nondegeneracy 
with respect to {\em variables} in
the system $\pde_\variablessmall (\mathcal{M})${\bf ?}}
\end{Question}

The thing is that the hypothesis of $2$-nondegeneracy
with respect to {\em parameters} (and not variables!)
has already been used to set up this system
$\pde_\variablessmall (\mathcal{M})$.

\begin{Proposition}
The system $\pde_\variablessmall (\mathcal{M})$ coming from
an $\mathcal{M}$ which is $2$-nondegenerate with respect
to parameters is also $2$-nondegenerate with respect
to variables if and only if:
\[
F_{z_xz_x}
\,\neq\,
0.
\]
\end{Proposition}

\proof
Start by rewriting, in the variables $(x,y,a,b,c)$ instead
of in the variables 
$(x,y,z,z_x,z_{xx})$, the functional equality
which was used in Section~{\ref{transfer-vector-fields-M-PDE}}
to perform a transfer of derivations:
\[
\Faux
\Big(
x,y,\,
Q\big(x,y,a,b,c\big),\,
Q_x\big(x,y,a,b,c\big),\,
Q_{xx}\big(x,y,a,b,c\big)
\Big)
\,\,\equiv\,\,
\Gaux
\big(
x,y,a,b,c
\big),
\]
and differentiate it with respect to $a$, $b$, $c$ to get
alternative formulas for vector fields (derivations):
\[
\aligned
Q_a\,
\frac{\partial}{\partial z}
+
Q_{xa}\,
\frac{\partial}{\partial z_x}
+
Q_{xxa}\,
\frac{\partial}{\partial z_{xx}}
&
\,=\,
\frac{\partial}{\partial a},
\\
Q_b\,
\frac{\partial}{\partial z}
+
Q_{xb}\,
\frac{\partial}{\partial z_x}
+
Q_{xxb}\,
\frac{\partial}{\partial z_{xx}}
&
\,=\,
\frac{\partial}{\partial b},
\\
Q_c\,
\frac{\partial}{\partial z}
+
Q_{xc}\,
\frac{\partial}{\partial z_x}
+
Q_{xxc}\,
\frac{\partial}{\partial z_{xx}}
&
\,=\,
\frac{\partial}{\partial c}.
\endaligned
\]

The path of computations is to eliminate 
$\frac{\partial}{\partial z}$ from lines $1$ and $3$:
\[
0
+
\big(
Q_c\,Q_{xa}
-
Q_a\,Q_{xc}
\big)\,
\frac{\partial}{\partial z_x}
+
\big(
Q_c\,Q_{xxa}
-
Q_a\,Q_{xxc}
\big)\,
\frac{\partial}{\partial z_{xx}}
\,=\,
Q_c\,
\frac{\partial}{\partial a}
-
Q_a\,
\frac{\partial}{\partial c},
\] 
to apply this derivation to the identity of
Proposition~{\ref{Proposition-k-equal-F-zx}}:
\[
F_{z_x}
\big(
x,y,z,z_x,z_{xx}
\big)
\,=\,
\frac{-P_zP_{ay}+P_yP_{az}}{-P_zP_{ax}+P_xP_{az}}
\big(
a,b,x,y,z
\big),
\]
taking advantage of two facts: firstly that $F_{z_{xx}} \equiv 0$ 
by degeneracy of the Levi form; secondly that
$P$ is independent of $c$; this conducts to 
observe that two terms, one in each side, disappear:
\[
\big(
Q_cQ_{xa}-Q_aQ_{xc}
\big)\,
F_{z_xz_x}
+
\big(
Q_cQ_{xxa}-Q_aQ_{xxc}
\big)\,
\zero{F_{z_xz_{xx}}}
\,\,=\,\,
Q_c\,
\frac{\partial}{\partial a}\,
\bigg(
\frac{-P_zP_{ay}+P_yP_{az}}{-P_zP_{ax}+P_xP_{az}}
\bigg)
-
Q_a\,
\zero{
\frac{\partial}{\partial c}\,
\big(
\same
\big)};
\]
it then suffices to solve to conclude:
\begin{align}
F_{z_xz_x}
&
\,=\,
\frac{Q_c}{Q_cQ_{xa}-Q_aQ_{xc}}\,
\frac{\partial}{\partial a}\,
\bigg(
\frac{-P_zP_{ay}+P_yP_{az}}{-P_zP_{ax}+P_xP_{az}}
\bigg)
\notag
\\
&
\,=\,
-\,
\frac{Q_c}{-Q_cQ_{xa}-Q_aQ_{xc}}\,
P_z\,
\frac{
\left\vert\!
\begin{array}{ccc}
P_x & P_y & P_z
\\
P_{ax} & P_{ay} & P_{az}
\\
P_{aax} & P_{aay} & P_{aaz}
\end{array}
\!\right\vert
}{
\big(
-\,P_z\,P_{ax}+P_x\,P_{az}
\big)^2
}
\notag
\\
&
\,=\,
\frac{P_z\,P_z}{\big(-P_zP_{ax}+P_xP_{az}\big)^3}\,
\left\vert\!
\begin{array}{ccc}
P_x & P_y & P_z
\\
P_{ax} & P_{ay} & P_{az}
\\
P_{aax} & P_{aay} & P_{aaz}
\end{array}
\!\right\vert
\,\,\,\neq\,\,
0.
\qedhere
\end{align}
\endproof

\Section{\bf Complete Integrability}
\label{complete-integrability}
\HEAD{{\ref{complete-integrability}}.~{\sf Complete Integrability}
}{
Joël {\sc Merker} (Orsay)}

By construction, in the system $z_y = F$, $z_{xxx} = H$, we
have $F = Q_y$ and $H = Q_{xxx}$, hence from $\big( Q_y\big)_{xxx}
= \big( Q_{xxx} \big)_y$, we deduce a compatibility constraint
on $F$ and $H$:
\[
\Dop_x
\big(
\Dop_x
(\Dop_x(F))
\big)
\,=\,
\Dop_y
\big(H\big),
\] 
in terms of the two {\sl total differentiation operators:}
\[
\aligned
\Dop_x
&
\,:=\,
\frac{\partial}{\partial x}
+
z_x\,
\frac{\partial}{\partial z}
+
z_{xx}\,
\frac{\partial}{\partial z_x}
+
H\,
\frac{\partial}{\partial z_{xx}},
\\
\Dop_y
&
\,:=\,
\frac{\partial}{\partial y}
+
F\,
\frac{\partial}{\partial z}
+
\Dop_x(F)\,
\frac{\partial}{\partial z_x}
+
\Dop_x\big(\Dop_x(F)\big)\,
\frac{\partial}{\partial z_{xx}}.
\endaligned
\]
The converse holds true, and is an elementary consequence of the 
Frobenius theorem~{\cite{Merker-2008}}.

\begin{Theorem}
If two analytic functions $F$ and $H$ of 
$\big(x,y,z,z_x,z_{xx})$ satisfy:
\[
\Dop_x
\big(
\Dop_x(\Dop_x(F))
\big)
\,=\,
\Dop_y
\big(H\big),
\]
then the analytic $\pde$ system $z_y = F$, $z_{xxx} = H$ is 
completely integrable in the sense that there exists a 
$\mathcal{C}^\omega$ family of solutions:
\[
z
\,=\,
Q
\big(
x,y,z(0,0),z_x(0,0),z_{xx}(0,0)
\big),
\]
parametrized by initial conditions:
\reqnomode\usetagform{EngelLie}
\begin{align}
Q
\big(0,0,z(0,0),z_x(0,0),z_{xx}(0,0)\big)
&
\,\equiv\,
z(0,0),
\notag
\\
Q_x
\big(0,0,z(0,0),z_x(0,0),z_{xx}(0,0)\big)
&
\,\equiv\,
z_x(0,0),
\notag
\\
Q_{xx}
\big(0,0,z(0,0),z_x(0,0),z_{xx}(0,0)\big)
&
\,\equiv\,
z_{xx}(0,0).
\tag{\qed}
\end{align}
\end{Theorem}

\Section{\bf Initial $G$-Structure for Equivalences 
$\mathcal{M} \overset{\sim}{\longrightarrow} \mathcal{M}'$}
\label{initial-G-structure-M-M}
\HEAD{{\ref{initial-G-structure-M-M}}.~{\sf Initial 
$G$-Structure for Equivalences 
$\mathcal{M} \overset{\sim}{\longrightarrow} \mathcal{M}'$}
}{
Joël {\sc Merker} (Orsay)}

Suppose given an equivalence in $\Diff_\variablessmall \times
\Diff_\parameterssmall$: 
\[
(F,\Phi)
\,=\,
(f,g,\varphi,\psi) 
\colon\ \ \
\mathcal{M}
\,\longrightarrow\,
\mathcal{M}'.
\]
Work in coordinates $(x,y,a,b,c)$ on $\mathcal{M}$ and
$(x',y',a',b',c')$ on $\mathcal{M}'$, and take
the two frames of vector fields (now written in this precise
order):
\[
\aligned
\mathcal{T}_c
&
\,=\,
\frac{\partial}{\partial c},
\\
\mathcal{L}_a
&
\,=\,
\frac{\partial}{\partial a}
-
\frac{Q_a}{Q_c}\,
\frac{\partial}{\partial c},
\\
\mathcal{L}_b
&
\,=\,
\frac{\partial}{\partial b}
-
\frac{Q_b}{Q_c}\,
\frac{\partial}{\partial c},
\\
\mathcal{K}_x
&
\,=\,
\frac{\partial}{\partial x},
\\
\mathcal{K}_y
&
\,=\,
\frac{\partial}{\partial y},
\endaligned
\ \ \ \ \ \ \ \ \ \ \ \ \
\text{and}
\ \ \ \ \ \ \ \ \ \ \ \ \
\aligned
\mathcal{T}_{c'}'
&
\,=\,
\frac{\partial}{\partial c'},
\\
\mathcal{L}_{a'}'
&
\,=\,
\frac{\partial}{\partial a'}
-
\frac{Q_{a'}'}{Q_{c'}'}\,
\frac{\partial}{\partial c'},
\\
\mathcal{L}_{b'}'
&
\,=\,
\frac{\partial}{\partial b'}
-
\frac{Q_{b'}'}{Q_{c'}'}\,
\frac{\partial}{\partial c'},
\\
\mathcal{K}_{x'}'
&
\,=\,
\frac{\partial}{\partial x'},
\\
\mathcal{K}_{y'}'
&
\,=\,
\frac{\partial}{\partial y'}.
\endaligned
\]
On $\mathcal{M}$ and on $\mathcal{M}'$, 
the two pairs of Levi kernel direction fields are
generated by:
\[
\aligned
\mathcal{K}_\kersmall^\variablessmall
&
\,=\,
\kaux\,
\mathcal{K}_x
+
\mathcal{K}_y,
\\
\mathcal{L}_\kersmall^\parameterssmall
&
\,=\,
\laux\,
\mathcal{L}_a
+
\mathcal{L}_b,
\endaligned
\ \ \ \ \ \ \ \ \ \ \ \ \
\text{and}
\ \ \ \ \ \ \ \ \ \ \ \ \
\aligned
{\mathcal{K}'}_\kersmall^\variablessmall
&
\,=\,
\kaux'\,
\mathcal{K}_{x'}'
+
\mathcal{K}_{y'}',
\\
{\mathcal{L}'}_\kersmall^\parameterssmall
&
\,=\,
\laux'\,
\mathcal{L}_{a'}'
+
\mathcal{L}_{b'}',
\endaligned
\] 
and Lemma~{\ref{Lemma-transfer-from-Q-to-P}} enables to
reexpress $\kaux$ in terms of $Q$:
\[
\aligned
\kaux
&
\,=\,
-\,
\frac{-P_zP_{ay}+P_yP_{az}}{-P_zP_{ax}+P_xP_{az}},
\ \ \ \ \ \ \ \ \ \ \ \ \ \ \ \ \ \ \ \ \ \ \ \ \ \
\laux
\,=\,
-\,
\frac{-Q_cQ_{xb}+Q_bQ_{xc}}{-Q_cQ_{xa}+Q_aQ_{xc}},
\\
&
\,=\,
-\,
\frac{-Q_cQ_{ya}+Q_aQ_{yc}}{-Q_cQ_{xa}+Q_aQ_{xc}},
\endaligned
\]
with similar expressions for $\kaux'$ and $\laux'$.

\begin{Observation}
Through an equivalence $(F, \Phi) \colon \mathcal{M}
\longrightarrow \mathcal{M}'$, Levi-kernel directions
transfer one to another:
\[
\aligned
(F,\Phi)_*
\big(
T_\kersmall^\variablessmall\mathcal{M}
\big)
&
\,=\,
T_\kersmall^\variablessmall\mathcal{M}',
\\
(F,\Phi)_*
\big(
T_\kersmall^\parameterssmall\mathcal{M}
\big)
&
\,=\,
T_\kersmall^\parameterssmall\mathcal{M}'.
\endaligned
\]
\end{Observation}

\proof
This comes from 
$(F,\Phi)_\ast \big( T^{\variablessmall \,/\,
\parameterssmall} \mathcal{M} \big) = 
T^{\variablessmall \,/\,
\parameterssmall} \mathcal{M}'$, from 
the definitions of Levi kernels,
and from the invariancy of Levi forms
expressed by the matrix identities of
Section~{\ref{invariance-nonholonomy-form-and-dual}}.
\endproof

As an important consequence, there are certain functions 
${\sf f}_3$ and ${\sf h}_3$ such that, after unwritten
push-forward:
\[
\aligned
\laux'\,\mathcal{L}_{a'}'
+
\mathcal{L}_{b'}'
&
\,=\,
{\sf f}_3
\cdot
\big(
\laux\,\mathcal{L}_a
+
\mathcal{L}_b
\big),
\\
\kaux'\,\mathcal{K}_{x'}'
+
\mathcal{K}_{y'}'
&
\,=\,
{\sf h}_3
\cdot
\big(
\kaux\,\mathcal{K}_x
+
\mathcal{K}_y
\big).
\endaligned
\]
Consequently, there exist $11$ functions on 
$\mathcal{M}$ such that:
\[
\left(\!
\begin{array}{c}
\mathcal{T}_{c'}'
\\
\mathcal{L}_{a'}'
\\
\laux'\mathcal{L}_{a'}'+\mathcal{L}_{b'}'
\\
\mathcal{K}_{x'}'
\\
\kaux'\mathcal{K}_{x'}'+\mathcal{K}_{y'}'
\end{array}
\!\right)
\,=\,
\left(\!
\begin{array}{ccccc}
{\sf a} & {\sf b}_1 & {\sf b}_2 & {\sf c}_1 & {\sf c}_2
\\
0 & {\sf f}_1 & {\sf f}_2 & 0 & 0
\\
0 & 0 & {\sf f}_3 & 0 & 0
\\
0 & 0 & 0 & {\sf h}_1 & {\sf h}_2
\\
0 & 0 & 0 & 0 & {\sf h}_3
\end{array}
\!\right)\,
\left(\!
\begin{array}{c}
\mathcal{T}_{c}
\\
\mathcal{L}_{a}
\\
\laux\mathcal{L}_{a}+\mathcal{L}_{b}
\\
\mathcal{K}_{x}
\\
\kaux\mathcal{K}_{x}+\mathcal{K}_{y}
\end{array}
\!\right).
\]
The initial $G$-structure for equivalences of such $\mathcal{M}$ 
(in terms of vector fields)
is therefore represented by such (invertible) matrices.
Before starting Cartan's method,
it remains only to re-express this in terms of differential 
$1$-forms.

The coframe dual to the frame on $\mathcal{M}$ (in this order):
\[
\Big\{
\mathcal{T}_c,\ \ \ \ \
\mathcal{L}_a,\ \ \ \ \
\laux\mathcal{L}_a+\mathcal{L}_b,\ \ \ \ \
\mathcal{K}_x,\ \ \ \ \
\kaux\,\mathcal{K}_x+\mathcal{K}_y
\Big\}
\]
is:
\[
\Big\{
{\textstyle{\frac{Q_a}{Q_c}}}da+
{\textstyle{\frac{Q_b}{Q_c}}}db+dc,\ \ \ \ \
da-\laux db,\ \ \ \ \
db,\ \ \ \ \
dx-\kaux dy,\ \ \ \ \
dy
\Big\}.
\]
A plain transposition of the above $5 \times 5$ matrix then yields
with new functions:
\[
\left(\!
\begin{array}{c}
Q_{a'}'da'+Q_{b'}'db'+Q_{c'}'dc'
\\
da'-\laux'db'
\\
db'
\\
dx'-\kaux'dy'
\\
dy'
\end{array}
\!\right)
\,=\,
\left(\!
\begin{array}{ccccc}
{\sf a} & 0 & 0 & 0 & 0
\\
{\sf b}_1 & {\sf f}_1 & 0 & 0 & 0
\\
{\sf b}_2 & {\sf f}_2 & {\sf f}_3 & 0 & 0
\\
{\sf c}_1 & 0 & 0 & {\sf h}_1 & 0
\\
{\sf c}_2 & 0 & 0 & {\sf h}_2 & {\sf h}_3
\end{array}
\!\right)\,
\left(\!
\begin{array}{c}
Q_ada+Q_bdb+dc
\\
da-\laux db
\\
db
\\
dx-\kaux dy
\\
dy
\end{array}
\!\right).
\]

\Section{\bf Triangular Initial $G$-structure for Equivalences 
$\pde_\variablessmall(\mathcal{M})
\overset{\sim}{\longrightarrow} 
\pde_\variablessmall(\mathcal{M}')$}
\label{initial-G-structure-pde-var-M}
\HEAD{{\ref{initial-G-structure-pde-var-M}}.~{\sf Triangular Initial 
$G$-Structure for Equivalences 
$\pde_\variablessmall(\mathcal{M})
\overset{\sim}{\longrightarrow} 
\pde_\variablessmall(\mathcal{M}')$}
}{
Joël {\sc Merker} (Orsay)}

By what precedes, on $\mathcal{M}$ equipped with coordinates
$(x,y,a,b,c)$, there is a natural coframe:
\[
\Big\{
{\textstyle{\frac{Q_a}{Q_c}}}\,da+
{\textstyle{\frac{Q_b}{Q_c}}}\,db+dc,\ \ \ \ \
da-\laux\,db,\ \ \ \ \
db,\ \ \ \ \
dx-\kaux\,dy,\ \ \ \ \
dy
\Big\},
\]
on which the two invariant
plane fields $T^\variablessmall \mathcal{M}$ and
$T^\parameterssmall \mathcal{M}$ and the two invariant
Levi-kernel direction fields
$T_\kersmall^\variablessmall \mathcal{M}$ and
$T_\kersmall^\parameterssmall \mathcal{M}$
are visible. The next goal is to transmit this initial 
geometry to the associated system
$\pde_\variablessmall(\mathcal{M})$.

As is known, the transfer:
\[
\big(x,y,a,b,c\big)
\xleftrightarrow[{\rule[0pt]{40pt}{0pt}}]{}
\big(x,y,z,z_x,z_{xx}\big),
\]
is represented by a known map accompanied with its inverse:
\leqnomode\usetagform{default}
\begin{align}
\label{transfer-5-jet-M}
\aligned
x
&
\,=\,
x,
\\
y
&
\,=\,
y,
\\
Q(x,y,a,b,c)
&
\,=\,
z,
\\
Q_x(x,y,a,b,c)
&
\,=\,
z_x,
\\
Q_{xx}(x,y,a,b,c)
&
\,=\,
z_{xx}
\endaligned
\ \ \ \ \ \ \ \ \ \ \ \ \ \ \ \ \ \ \ \ \ \ \ \ \ \
\aligned
x
&
\,=\,
x,
\\
y
&
\,=\,
y,
\\
a
&
\,=\,
A\big(x,y,z,z_x,z_{xx}\big),
\\
b
&
\,=\,
B\big(x,y,z,z_x,z_{xx}\big),
\\
c
&
\,=\,
C\big(x,y,z,z_x,z_{xx}\big).
\endaligned
\end{align}

Next, on the system $\pde_\variablessmall(\mathcal{M})$:
\[
\aligned
z_y
&
\,=\,
F\big(x,y,z,z_x,z_{xx}\big)
\,=\,
Q_y
\Big(
x,y,\,
A\big(x,y,z,z_x,z_{xx}\big),
B\big(x,y,z,z_x,z_{xx}\big),
C\big(x,y,z,z_x,z_{xx}\big)
\Big),
\\
z_{xxx}
&
\,=\,
H\big(x,y,z,z_x,z_{xx}\big)
\,=\,
Q_{xxx}
\Big(
x,y,\,
A\big(x,y,z,z_x,z_{xx}\big),
B\big(x,y,z,z_x,z_{xx}\big),
C\big(x,y,z,z_x,z_{xx}\big)
\Big),
\endaligned
\]
there is a natural coframe consisting of $5$ differential
$1$-forms, the first $3$ being contact forms pulled-back to the 
$\pde$ system:
\[
\aligned
\lambda
&
\,:=\,
dz
-
z_x\,dx
-
F\,dy,
\\
\mu_1
&
\,:=\,
dz_x
-
z_{xx}\,dx
-
\Dop_x(F)\,dy,
\\
\mu_2
&
\,:=\,
dz_{xx}
-
H\,dx
-
\Dop_x(\Dop_x(F))\,dy,
\\
\nu_1
&
\,:=\,
dx,
\\
\nu_2
&
\,:=\,
dy.
\endaligned
\]
We want to relate these $1$-forms to the coframe on $\mathcal{M}$
introduced above. More precisely, in some (possibly more) appropriate
initial coframe, we want to determine the initial $G$-structure
for point equivalences:
\[
\big(x,y,z,z_x,z_{xx}\big)
\xrightarrow[{\rule[0pt]{40pt}{0pt}}]{}
\big(x',y',z',z_{x'}',z_{x'x'}'\big),
\]
that transfer a system $z_x = F$, $z_{xxx} = H$, to a similar
one $z_{y'} ' = F'$, $z_{x'x'x'}' = H'$, having of course the same
geometric features.

First of all, since contact forms of any fixed jet order
must be sent to contact forms of the same jet order,
there is a triangular $3 \times 3$ matrix of functions
such that, after (unwritten) pullback:
\[
\left(\!
\begin{array}{c}
\lambda'
\\
\mu_1'
\\
\mu_2'
\end{array}
\!\right)
\,=\,
\left(\!
\begin{array}{ccc}
{\sf a} & 0 & 0
\\
{\sf b}_1 & {\sf f}_1 & 0
\\
{\sf b}_2 & {\sf f}_2 & {\sf f}_3
\end{array}
\!\right)
\left(\!
\begin{array}{c}
\lambda
\\
\mu_1
\\
\mu_2
\end{array}
\!\right).
\]
Also, since the considered transformations are punctual:
\[
(x,y,z)
\,\,\longmapsto\,\,
\big(
x'(x,y,z),\,
y'(x,y,z),\,
z'(x,y,z)
\big),
\]
and since the coframes $\{\lambda, \nu_1, \nu_2\}$ and
$\{dz, dx, dy\}$ have the same span, there are functions
such that:
\leqnomode\usetagform{default}
\begin{align}
\label{non-triangular-3x3}
\left(\!
\begin{array}{c}
\lambda'
\\
\nu_1'
\\
\nu_2'
\end{array}
\!\right)
\,=\,
\left(\!
\begin{array}{ccc}
{\sf a} & 0 & 0
\\
{\sf c}_1 & {\sf h}_1 & {\sf h}_4
\\
{\sf c}_2 & {\sf h}_2 & {\sf h}_3
\end{array}
\!\right)
\left(\!
\begin{array}{c}
\lambda
\\
\nu_1
\\
\nu_2
\end{array}
\!\right).
\end{align}
But we already saw that the hypothesis that $\mathcal{M}$
is {\em also} $2$-nondegenerate with respect to variables
conducts to choose:
\[
\big\{
dx-\kaux\,dy,\ \ \
dy
\big\}
\ \ \ \ \ \ \ \ \ \ \ \ \
\text{instead of}
\ \ \ \ \ \ \ \ \ \ \ \ \
\big\{
dx,\ \ \
dy
\big\},
\]
and since there is in
({\ref{transfer-5-jet-M}})
plainly $dx = dx$, $dy = dy$, and because we know already
that from
Proposition~{\ref{Proposition-k-equal-F-zx}}:
\[
-\,\kaux
\,=\,
F_{z_x},
\]
a more natural coframe to start with later on when
running Cartan's method is (notice the modification of $\nu_1$) is
the following:
\[
\aligned
\lambda
&
\,:=\,
dz
-
z_x\,dx
-
F\,dy,
\\
\mu_1
&
\,:=\,
dz_x
-
z_{xx}\,dx
-
\Dop(F)\,dy,
\\
\mu_2
&
\,:=\,
dz_{xx}
-
H\,dx
-
\Dop_x(\Dop_x(F))\,dy,
\\
\nu_1
&
\,:=\,
dx
+
F_{z_x}\,dy,
\\
\nu_2
&
\,:=\,
dy,
\endaligned
\]
so that initial equivalences would take the {\em triangular} form:
\[
\left(\!
\begin{array}{c}
\lambda'
\\
\mu_1'
\\
\mu_2'
\\
\nu_1'
\\
\nu_2'
\end{array}
\!\right)
\,=\,
\left(\!
\begin{array}{ccccc}
{\sf a} & 0 & 0 & 0 & 0
\\
{\sf b}_1 & {\sf f}_1 & 0 & 0 & 0
\\
{\sf b}_2 & {\sf f}_2 & {\sf f}_3 & 0 & 0
\\
{\sf c}_1 & 0 & 0 & {\sf h}_1 & 0
\\
{\sf c}_2 & 0 & 0 & {\sf h}_2 & {\sf h}_3
\end{array}
\!\right)
\left(\!
\begin{array}{c}
\lambda
\\
\mu_1
\\
\mu_2
\\
\nu_1
\\
\nu_2
\end{array}
\!\right),
\]
with ${\sf f}_4 = {\sf h}_4 = 0$ disappearing. {\em This is under
the assumption that $F_{z_xz_x} \neq 0$ only!}

In fact, without assuming this, namely starting from $\nu_1 = dx$,
$\nu_2 = dy$ and from a {\em nontriangular} initial $G$-structure,
one can easily re-find~{\cite{Merker-Nurowski-2020-a}} 
that $F_{z_xz_x}$ is a relative differential invariant,
and then assume again $F_{z_xz_x} \neq 0$. 

Before launching Cartan's method, let us end up by examining a bit
what the transfer $\mathcal{M} \longrightarrow \pde_\variablessmall
(\mathcal{M})$ can tell us about
${\sf f}_4$ and ${\sf h}_4$. Indeed, we want to confirm the fact that 
$\mu_1'$ is a linear combination of $\{ \lambda, \mu_1\}$,
without $\mu_2$, a fact which the general theory of contact forms
already gave, and which we want to see as a {\em consequence}
of the triangular form of the $G$-structure
for equivalences $\mathcal{M} 
\overset{\sim}{\longrightarrow} \mathcal{M}'$.

Abbreviate the considered jet space as:
\[
J_\horsmall^2
\,:=\,
\big\{
(x,y,z,z_x,z_{xx})
\big\}.
\]
Equivalences between submanifolds
of solutions are in one-to-one correspondence
with equivalences between systems of
partial differential equations:
\[
\xymatrix{
J_\horsmall^2 
\ar@<2pt>[rr] 
\ar@<2pt>[d] 
& &
\ar@<2pt>[ll] 
{J_\horsmall'}^2
\ar@<2pt>[d] 
\\
\mathcal{M}
\ar@<2pt>[rr] 
\ar@<2pt>[u] 
& &
\ar@<2pt>[ll] 
\mathcal{M}',
\ar@<2pt>[u] 
}
\ \ \ \ \ \ \ \ \ \ \ \ \ \ \ \ \ \ \ \ \ \ \ \ \ \
\xymatrix{
\ar@<2pt>[rr] 
\big(x,y,z,z_x,z_{xx}\big)
\ar@<2pt>[d] 
& &
\ar@<2pt>[ll] 
\big(x',y',z',z_{x'}',z_{x'x'}'\big)
\ar@<2pt>[d] 
\\ 
(x,y,a,b,c)
\ar@<2pt>[rr] 
\ar@<2pt>[u] 
& &
\ar@<2pt>[ll]
(x',y',a',b',c').
\ar@<2pt>[u] 
}
\]
All steps of Cartan's method act parallely!
Then the basic contact $1$-form transfers as:
\[
\aligned
\lambda
&
\,=\,
dz-z_x\,dx-F\,dy
\\
&
\,=\,
dQ-Q_x\,dx-Q_y\,dy
\\
&
\,=\,
Q_x\,dx+Q_y\,dy+Q_a\,da+Q_b\,db+Q_c\,dc
-
Q_x\,dx-Q_y\,dy
\\
&
\,=\,
Q_a\,da+Q_b\,db+Q_c\,dc,
\endaligned
\]
hence up to a nowhere vanishing factor (absorbed anyway in the
$G$-structure matrix), we recover our $1$-form 
$\frac{Q_a}{Q_c}\, da + \frac{Q_b}{Q_c}\, db + dc$.
Next:
\[
\aligned
\mu_1
&
\,=\,
dz_x
-
z_{xx}\,dx
-
\Dop_x(F)\,dy
\\
&
\,=\,
d(Q_x)
-
Q_{xx}\,dx
-
\Dop_x(Q_y)\,dy
\\
&
\,=\,
Q_{xx}\,dx+Q_{xy}\,dy+Q_{xa}\,da+Q_{xb}\,db+Q_{xc}\,dc
-Q_{xx}\,dx-Q_{xy}\,dy
\\
&
\,=\,
Q_{xa}\,da+Q_{xb}\,db+Q_{xc}\,dc,
\endaligned
\]
and similarly:
\[
\aligned
\mu_2
&
\,=\,
dz_{xx}-H\,dx-\Dop_x(\Dop_x(F))\,dy
\\
&
\,=\,
d(Q_{xx})-Q_{xxx}\,dx-Q_{xxy}\,dy
\\
&
\,=\,
Q_{xxa}\,da+Q_{xxb}\,db+Q_{xxc}\,dc.
\endaligned
\]

In summary, we have the formulas which provide the transfer:
\[
\xymatrix{
\big\{\lambda,\,\,\mu_1,\,\,\mu_2\big\}
\ar@<2pt>[d] 
\\
\big\{dc,\,\,db,\,\,da\big\},
\ar@<2pt>[u]
}
\]
namely:
\[
\aligned
\lambda
&
\,=\,
Q_a\,da+Q_b\,db+Q_c\,dc,
\\
\mu_1
&
\,=\,
Q_{xa}\,da+Q_{xb}\,db+Q_{xc}\,dc,
\\
\mu_2
&
\,=\,
Q_{xxa}\,da+Q_{xxb}\,db+Q_{xxc}\,dc,
\endaligned
\]
with invertible determinant $\Delta(Q) \neq 0$\,\,---\,\,by
assumption!\,\,---, so that $\{dc,\, db,\,da\}$ can inversely be 
expressed in terms of $\{\lambda, \mu_1, \mu_2\}$.

If we eliminate $dc$ from lines $1$ and $2$:
\[
\aligned
\mu_1
-
\frac{Q_{xc}}{Q_c}\,
\lambda
&
\,=\,
\frac{-Q_cQ_{xa}+Q_aQ_{xc}}{Q_c}\,
\bigg\{
da
+
\frac{-Q_cQ_{xb}+Q_bQ_{xc}}{-Q_cQ_{xa}+Q_aQ_{xc}}\,
db
\bigg\}
\\
&
\,=\,
\nonzero
\cdot
\big\{
da-\laux\,db
\big\},
\endaligned
\]
we recognize the $1$-form whose kernel, in 
$T^\parameterssmall \mathcal{M}$, spans the Levi kernel
bundle $T_\kersmall^\parameterssmall \mathcal{M}$.

To finish with these considerations, remembering that
we showed that:
\[
\aligned
Q_{a'}'da'+Q_{b'}'db'+Q_{c'}dc'
&
\,\,\in\,\,
\Span\,
\big\{
Q_a\,da+Q_b\,db+Q_c\,dc
\big\},
\\
da'-\laux'\,db'
&
\,\,\in\,\,
\Span\,
\big\{
Q_a\,da+Q_b\,db+Q_c\,dc,\ \ \
da-\laux\,db
\big\},
\endaligned
\]
we deduce that on the $\pde$ side:
\[
\aligned
\lambda'
&
\,\,\in\,\,
\Span\,
\{\lambda\},
\\
\mu_1'
&
\,\,\in\,\,
\Span\,
\big\{
\lambda,\,\,\mu_1
\big\},
\endaligned
\]
which means, as predicted by the property that
the jet order is preserved by pullbacks of contact forms,
that ${\sf f}_4 = 0$, {\em always}.

\smallskip

Without assuming ${\sf h}_4 = 0$, we may launch Cartan's equivalence
method for $\pde$ systems of the form $z_y = F$, $z_{xxx} = H$, with
$11 + 1$ independent group variables, and with the lifted coframe:
\[
\left(\!
\begin{array}{c}
\pmb{\lambda}
\\
\pmb{\mu_1}
\\
\pmb{\mu_2}
\\
\pmb{\nu_1}
\\
\pmb{\nu_2}
\end{array}
\!\right)
\,:=\,
\left(\!
\begin{array}{ccccc}
{\sf a} & 0 & 0 & 0 & 0
\\
{\sf b}_1 & {\sf f}_1 & 0 & 0 & 0
\\
{\sf b}_2 & {\sf f}_2 & {\sf f}_3 & 0 & 0
\\
{\sf c}_1 & 0 & 0 & {\sf h}_1 & {\sf h}_4
\\
{\sf c}_2 & 0 & 0 & {\sf h}_2 & {\sf h}_3
\end{array}
\!\right)
\left(\!
\begin{array}{c}
\lambda
\\
\mu_1
\\
\mu_2
\\
\nu_1
\\
\nu_2
\end{array}
\!\right).
\] 
This preliminary paper may now stop, because 
starting from the
basic point reached here, 
{\em advanced and non-straightforward} 
Cartan-type computations are conducted
in~{\cite{Merker-Nurowski-2020-a, 
Merker-Nurowski-2020-b, Merker-Nurowski-2021}}.



\vfill\end{document}

%% file: macros.tex
\newtheorem{Theorem}[equation]{Theorem}

\newtheorem{Proposition}[equation]{Proposition}

\newtheorem{Lemma}[equation]{Lemma}

\newtheorem{Corollary}[equation]{Corollary}
\newtheorem{Assertion}[equation]{Assertion}

\newtheorem{Observation}[equation]{Observation}

\newtheorem{FunctionalRelations}[equation]{Functional Relations}

\theoremstyle{definition}

\newtheorem{Definition-Notation}[equation]{D\'efinition-Notation}
\newtheorem{Definition}[equation]{Definition}

\newtheorem{Terminology}[equation]{Terminology}

\newtheorem{Example}[equation]{Example}

\newtheorem{Question}[equation]{Question}

\newcommand{\C}{\mathbb{C}}

\newcommand{\K}{\mathbb{K}}

\newcommand{\N}{\mathbb{N}}

\newcommand{\R}{\mathbb{R}}

\newcommand{\MM}{\text{\sc m}}
\newcommand{\NN}{\text{\sc n}}

\newcommand{\TT}{\text{\sc t}}

\newcommand{\kaux}{{\text{\usefont{T1}{qcs}{m}{sl}k}}}
\newcommand{\laux}{{\text{\usefont{T1}{qcs}{m}{sl}l}}}

\newcommand{\zaux}{{\text{\usefont{T1}{qcs}{m}{sl}z}}}

\newcommand{\Faux}{{\text{\usefont{T1}{qcs}{m}{sl}F}}}
\newcommand{\Gaux}{{\text{\usefont{T1}{qcs}{m}{sl}G}}}

\definecolor{blue}{cmyk}{1.,1.,0.,0.63}
\definecolor{red}{cmyk}{0.,1.,1.,0.63}
\definecolor{green}{cmyk}{1.,0.,1.,0.63}
\definecolor{black}{cmyk}{1.,1.,1.,1.}

\newcommand{\green}{\textcolor{green}}
\newcommand{\red}{\textcolor{red}}

\makeatletter
\renewcommand{\@fnsymbol}[1]
{\ensuremath{\ifcase#1\or $*$\or $**$\or $***$\or $****$\or $*****$
\else\@ctrerr\fi}}
\makeatother

\newcommand{\HEAD}[2]{%
\pagestyle{fancy}
\fancyhead[RO]{\tiny\sf\thepage}
\fancyhead[CO]{{\tiny\sf #1}}
\fancyhead[LE]{\tiny\sf\thepage}
\fancyhead[CE]{{\tiny\sf #2}}
\fancyfoot{}}

\numberwithin{equation}{section}

\newcommand{\Section}[1]{
\renewcommand{\thesection}{\bf\arabic{section}}
\section{#1}
\renewcommand{\thesection}{\arabic{section}}}

\newcommand{\style}[1]{{\sf #1}}

\newcommand{\stylesmall}[1]{{\sf #1}}

\newcommand{\Aut}{\style{Aut}}

\newcommand{\codim}{\style{codim}}

\newcommand{\constant}{\style{constant}}

\newcommand{\CRdim}{\style{CRdim}}

\renewcommand{\det}{\style{det}}

\newcommand{\Diff}{\style{Diff}}

\renewcommand{\dim}{\style{dim}}

\newcommand{\genrank}{\style{genrank}}

\newcommand{\horsmall}{\stylesmall{hor}}

\newcommand{\Id}{\style{Id}}
\newcommand{\Idsmall}{\stylesmall{Id}}

\newcommand{\Jac}{\style{Jac}}

\newcommand{\Ker}{\style{Ker}}

\newcommand{\kersmall}{\stylesmall{ker}}

\newcommand{\Levi}{\style{Levi}}

\renewcommand{\lim}{\style{lim}}

\renewcommand{\min}{\style{min}}

\renewcommand{\mod}{\style{mod}}

\newcommand{\nonzero}{\style{nonzero}}

\newcommand{\parameters}{\style{par}}
\newcommand{\parameterssmall}{\stylesmall{par}}

\newcommand{\Polynomial}{\style{Polynomial}}

\newcommand{\rank}{\style{rank}}

\newcommand{\reducedsmall}{\stylesmall{reduced}}

\newcommand{\same}{\style{same}}

\newcommand{\Span}{\style{Span}}

\newcommand{\variables}{\style{var}}
\newcommand{\variablessmall}{\stylesmall{var}}

\newcommand{\bignorm}{\big\vert\!\big\vert}

\newcommand{\Dop}{{\sf D}}

\newcommand{\Hall}{\Hall}

\newcommand{\isqrt}{{\scriptstyle{\sqrt{-1}}}}

\newcommand{\pde}{\text{\sc pde}}

\newcommand{\smallbullet}{{\scriptscriptstyle{\bullet}}}

\newcommand{\smallsum}[1]{
\underset{#1}{\raisebox{1pt}{$\sum$\,}}
}

\newcommand{\vf}{\vfill\end{document}}

\newcommand{\zero}[1]{\underline{#1}_{\red{\circ}}}

\makeatletter
\def\hlinethick#1{%
\noalign{\ifnum0=`}\fi\hrule \@height #1 \futurelet
\reserved@a\@xhline}
\makeatother

\newlength{\myskip}


%% file: double-projection.pdf_t
\begin{picture}(0,0)%
\includegraphics{double-projection.pdf}%
\end{picture}%
\setlength{\unitlength}{4144sp}%
\begingroup\makeatletter\ifx\SetFigFont\undefined%
\gdef\SetFigFont#1#2#3#4#5{%
  \reset@font\fontsize{#1}{#2pt}%
  \fontfamily{#3}\fontseries{#4}\fontshape{#5}%
  \selectfont}%
\fi\endgroup%
\begin{picture}(3760,3577)(1453,-3391)
\put(1468,-538){\makebox(0,0)[lb]{\smash{{\SetFigFont{8}{9.6}{\familydefault}{\mddefault}{\updefault}{\color[rgb]{0,0,0}$(a,b)$}%
}}}}
\put(1659,-2656){\makebox(0,0)[lb]{\smash{{\SetFigFont{8}{9.6}{\familydefault}{\mddefault}{\updefault}{\color[rgb]{0,0,0}$\K_\parameterssmall^\MM$}%
}}}}
\put(5198,-3244){\makebox(0,0)[lb]{\smash{{\SetFigFont{8}{9.6}{\familydefault}{\mddefault}{\updefault}{\color[rgb]{0,0,0}$\K_\variablessmall^\NN$}%
}}}}
\put(4348,-3338){\makebox(0,0)[lb]{\smash{{\SetFigFont{8}{9.6}{\familydefault}{\mddefault}{\updefault}{\color[rgb]{0,0,0}$(x,y)$}%
}}}}
\put(3843,-523){\makebox(0,0)[lb]{\smash{{\SetFigFont{8}{9.6}{\familydefault}{\mddefault}{\updefault}{\color[rgb]{0,0,0}$\mathcal{Q}_{a,b}$}%
}}}}
\put(4713,-234){\makebox(0,0)[lb]{\smash{{\SetFigFont{8}{9.6}{\familydefault}{\mddefault}{\updefault}{\color[rgb]{0,0,0}$\mathcal{M}$}%
}}}}
\put(5071,  1){\makebox(0,0)[lb]{\smash{{\SetFigFont{8}{9.6}{\familydefault}{\mddefault}{\updefault}{\color[rgb]{0,0,0}$\K^\NN \times \K^\MM$}%
}}}}
\put(4397,-1145){\makebox(0,0)[lb]{\smash{{\SetFigFont{8}{9.6}{\familydefault}{\mddefault}{\updefault}{\color[rgb]{0,0,0}$\mathcal{P}_{x,y}$}%
}}}}
\put(3858,-2745){\makebox(0,0)[lb]{\smash{{\SetFigFont{8}{9.6}{\familydefault}{\mddefault}{\updefault}{\color[rgb]{0,0,0}$\pi_\variablessmall$}%
}}}}
\put(2217,-1290){\makebox(0,0)[lb]{\smash{{\SetFigFont{8}{9.6}{\familydefault}{\mddefault}{\updefault}{\color[rgb]{0,0,0}$\pi_\parameterssmall$}%
}}}}
\end{picture}%

%% file: par-var-equivalences.pdf_t
\begin{picture}(0,0)%
\includegraphics{par-var-equivalences.pdf}%
\end{picture}%
\setlength{\unitlength}{4144sp}%
\begingroup\makeatletter\ifx\SetFigFont\undefined%
\gdef\SetFigFont#1#2#3#4#5{%
  \reset@font\fontsize{#1}{#2pt}%
  \fontfamily{#3}\fontseries{#4}\fontshape{#5}%
  \selectfont}%
\fi\endgroup%
\begin{picture}(7518,3043)(1644,-3307)
\put(1659,-2656){\makebox(0,0)[lb]{\smash{{\SetFigFont{8}{9.6}{\familydefault}{\mddefault}{\updefault}{\color[rgb]{0,0,0}$\K_\parameterssmall^\MM$}%
}}}}
\put(4276,-691){\makebox(0,0)[lb]{\smash{{\SetFigFont{8}{9.6}{\familydefault}{\mddefault}{\updefault}{\color[rgb]{0,0,0}$\mathcal{M}$}%
}}}}
\put(2189,-1499){\makebox(0,0)[lb]{\smash{{\SetFigFont{8}{9.6}{\familydefault}{\mddefault}{\updefault}{\color[rgb]{0,0,0}$\pi_\parameterssmall$}%
}}}}
\put(3631,-2752){\makebox(0,0)[lb]{\smash{{\SetFigFont{8}{9.6}{\familydefault}{\mddefault}{\updefault}{\color[rgb]{0,0,0}$\pi_\variablessmall$}%
}}}}
\put(6537,-691){\rotatebox{360.0}{\makebox(0,0)[rb]{\smash{{\SetFigFont{8}{9.6}{\familydefault}{\mddefault}{\updefault}{\color[rgb]{0,0,0}$\mathcal{M}'$}%
}}}}}
\put(4748,-3247){\makebox(0,0)[lb]{\smash{{\SetFigFont{8}{9.6}{\familydefault}{\mddefault}{\updefault}{\color[rgb]{0,0,0}$\K_\variablessmall^\NN$}%
}}}}
\put(6066,-3240){\rotatebox{360.0}{\makebox(0,0)[rb]{\smash{{\SetFigFont{8}{9.6}{\familydefault}{\mddefault}{\updefault}{\color[rgb]{0,0,0}${\K'}_\variablessmall^\NN$}%
}}}}}
\put(8621,-1521){\rotatebox{360.0}{\makebox(0,0)[rb]{\smash{{\SetFigFont{8}{9.6}{\familydefault}{\mddefault}{\updefault}{\color[rgb]{0,0,0}$\pi_\parameterssmall'$}%
}}}}}
\put(9147,-2662){\rotatebox{360.0}{\makebox(0,0)[rb]{\smash{{\SetFigFont{8}{9.6}{\familydefault}{\mddefault}{\updefault}{\color[rgb]{0,0,0}${\K'}_\parameterssmall^\MM$}%
}}}}}
\put(7441,-2729){\rotatebox{360.0}{\makebox(0,0)[rb]{\smash{{\SetFigFont{8}{9.6}{\familydefault}{\mddefault}{\updefault}{\color[rgb]{0,0,0}$\pi_\variablessmall'$}%
}}}}}
\put(5216,-1360){\makebox(0,0)[lb]{\smash{{\SetFigFont{8}{9.6}{\familydefault}{\mddefault}{\updefault}{\color[rgb]{0,0,0}$(F,\Phi)$}%
}}}}
\end{picture}%

%% file: 0-1-3-0-stairs.pdf_t
\begin{picture}(0,0)%
\includegraphics{0-1-3-0-stairs.pdf}%
\end{picture}%
\setlength{\unitlength}{4144sp}%
\begingroup\makeatletter\ifx\SetFigFont\undefined%
\gdef\SetFigFont#1#2#3#4#5{%
  \reset@font\fontsize{#1}{#2pt}%
  \fontfamily{#3}\fontseries{#4}\fontshape{#5}%
  \selectfont}%
\fi\endgroup%
\begin{picture}(4254,2749)(1384,-3023)
\put(1399,-2115){\makebox(0,0)[lb]{\smash{{\SetFigFont{8}{9.6}{\familydefault}{\mddefault}{\updefault}{\color[rgb]{0,0,0}\red{$(0,1)$}}%
}}}}
\put(1656,-2950){\makebox(0,0)[lb]{\smash{{\SetFigFont{8}{9.6}{\familydefault}{\mddefault}{\updefault}{\color[rgb]{0,0,0}$(0,0)$}%
}}}}
\put(2329,-2959){\makebox(0,0)[lb]{\smash{{\SetFigFont{8}{9.6}{\familydefault}{\mddefault}{\updefault}{\color[rgb]{0,0,0}$(1,0)$}%
}}}}
\put(3012,-2958){\makebox(0,0)[lb]{\smash{{\SetFigFont{8}{9.6}{\familydefault}{\mddefault}{\updefault}{\color[rgb]{0,0,0}$(2,0)$}%
}}}}
\put(3668,-2961){\makebox(0,0)[lb]{\smash{{\SetFigFont{8}{9.6}{\familydefault}{\mddefault}{\updefault}{\color[rgb]{0,0,0}\red{$(3,0)$}}%
}}}}
\end{picture}%